\documentclass[10pt]{article}

\PassOptionsToPackage{x11names}{xcolor}

\usepackage{tikz}
\usetikzlibrary{arrows,snakes,positioning,patterns} 

\usepackage{amsmath, amssymb,hyperref}

\usepackage[text={6.75in,9in},a4paper,portrait]{geometry}

\newtheorem{thm}{Theorem}
\newtheorem{lem}[thm]{Lemma}
\newtheorem{prop}[thm]{Proposition}
\newtheorem{cor}[thm]{Corollary}
\def\dem{\noindent\textbf{Proof}. }
\def\qed{{\ifhmode\unskip\nobreak\hfil\penalty50 \hskip1em \else\nobreak\fi
   \mbox{}\nobreak\hfil\qedtext%
   \parfillskip=0pt \finalhyphendemerits=0 \par\mbox{}\vspace{0pt}}}

 \def\comment#1{}

\def\qedtext{\ensuremath{\square}}
\RequirePackage{amsfonts}
\DeclareMathSymbol{\square}       {\mathord}{AMSa}{"03}

\def\com#1{}

\newcommand{\Z}{\mathbb{Z}}
\newcommand{\F}{\mathbb{F}}
\newcommand{\R}{\mathbb{R}}

\begin{document}

\begin{center}

{  \LARGE  An Algorithmic Approach to Algebraic  and Dynamical \vspace{0.2cm}\\ Cancellations associated to a Spectral Sequence 
}

\vspace{0.5cm}

M.A. Bertolim 
\hspace{0.3cm} D.V.S. Lima\footnote{Supported by FAPESP under grant 2010/08579-0.}  
\hspace{0.3cm}
M.P. Mello 
\hspace{0.3cm}K.A. de Rezende\footnote{Partially supported by CNPq under grant 302592/2010-5 and by FAPESP under grant 2012/18780-0.}
\hspace{0.3cm}  M. R. da
Silveira\footnote{Partially supported by FAPESP under grant 2012/18780-0.}    

\end{center}

\begin{abstract}
In this article we study algorithms that arise in both  topological and dynamical settings, namely, the Spectral Sequence Sweeping Algorithm (SSSA)  and the Row Cancellation Algorithm (RCA) for a filtered Morse chain complex on a manifold $M^{n}$. Both algorithms have as input a connection matrix and the results obtained in this article make it possible to establish a correspondence between the algebraic cancellations in SSSA and the dynamical cancellations  in RCA. 
\end{abstract}

{\small
}

\vspace{-.5cm}

\section*{Introduction}

Algorithms for computing spectral sequences are well known in the literature and have been implemented by several authors \cite{Romero:2006}.
For example, in the field of Computational Topology, see \cite{Edelsbrunner:2010},  there has been much interest in algorithms which solve problems anywhere from the reconstruction  of surfaces in computer graphics  to  the treatment of noise in data input.  
For example, to gather information about a surface  using a computer,  it is necessary to make use of a combinatorial representation. This information can be provided by making use of simplicial complexes. Data structures are then constructed with the purpose of  storing cell complex information. To retrieve topological connectivity information from this data, homology and spectral sequences are  used. 
%

We  approach spectral sequences from a dynamical systems point of view. Our primary interest resides in relating qualitative aspects of dynamical systems which can be coded  algebraically in a chain complex description, to topological aspects of the phase space and the asymptotic behavior of stable and unstable manifolds. The study of the latter, in a one-parameter family of flows, using a homotopical tool such as the Conley index, permits a deeper understanding of bifurcation behavior in a continuation, which is reflected, for example, in birth and death of critical points. There are many techniques that may be used to achieve such an endeavor. The underlying approach used here has been to bridge the algebraic-topological and dynamical realms by use of spectral sequences. 

In order to understand this more fully, let $f:M\to \mathbb{R}$ be a Morse function on a closed Riemannian $n$-manifold $M$ and  consider the {\it Morse chain complex} $(C,\Delta)$ where the Morse chain group $C=\{C_{k}(f)\}$ is defined as the free abelian group   generated by the critical points of $f$ and graded by their  Morse indices, i.e.,
\[C_k(f)=\bigoplus_{x\in {\rm Crit}_k(f)}\mathbb{Z}[x]\]
where $\mbox{Crit}_k(f)$ denotes the set of index $k$ critical points of $f$. In this case, the differential $\Delta:C\to C$ is a collection of homomorphisms $\Delta_k:C_k(f)\to C_{k-1}(f)$ which are  defined on a generator $x\in Crit_{k}(f)$ by 
\[\Delta_k (x)=\sum_{y\in {\rm Crit}_{k-1}(f)} n(x,y)[y],\]
where $n(x,y)$ denotes the intersection number of $x$ and $y$.

A {\it spectral sequence} $E=(E^r,d^{r})_{r\geq 0}$ is a sequence, such that $E^r$ is a { bigraded $R$-module} over a principal ideal domain  $R$, i.e.,  an indexed collection of  $R$-modules $E^{r}_{p,q}$ for all pair of integers $p$ and $q$; and $d^r$ is a {\it differential} of bidegree $(-r,r-1)$, i.e., it is a collection of homomorphisms $d^r:E^{r}_{p,q} \rightarrow E^{r}_{p-r,q+r-1}$ for all $p$ and $q$ such that ${d^r\circ d^{r}} = 0$. Moreover, for all $r\geq 0$ there exists an isomorphism $H(E^r)\approx E^{r+1}$ where $$H_{p,q}(E^r)=\frac{\mbox{Ker}\,
d^r:E^{r}_{p,q}\to E^r_{p-r,q+r-1}}{\mbox{Im}\, d^r:E^{r}_{p+r,
q-r+1}\to E^r_{p,q}}$$ is the homology module. For more details see \cite{Spanier:1966}.
 
{One way to visualize a spectral sequence is as a book  consisting  of pages such that the $r$-th page corresponds to the bigraded module $E^{r}$.   On each page there are homomorphisms between the modules forming sequences of chain complexes.  Moreover, the homology module of the $r$-th page is precisely the $(r+1)$-th  page.  See Figure~\ref{fig:spectralsequencepages}. }

\begin{figure}[h!]
\centering
\includegraphics[scale=0.8]{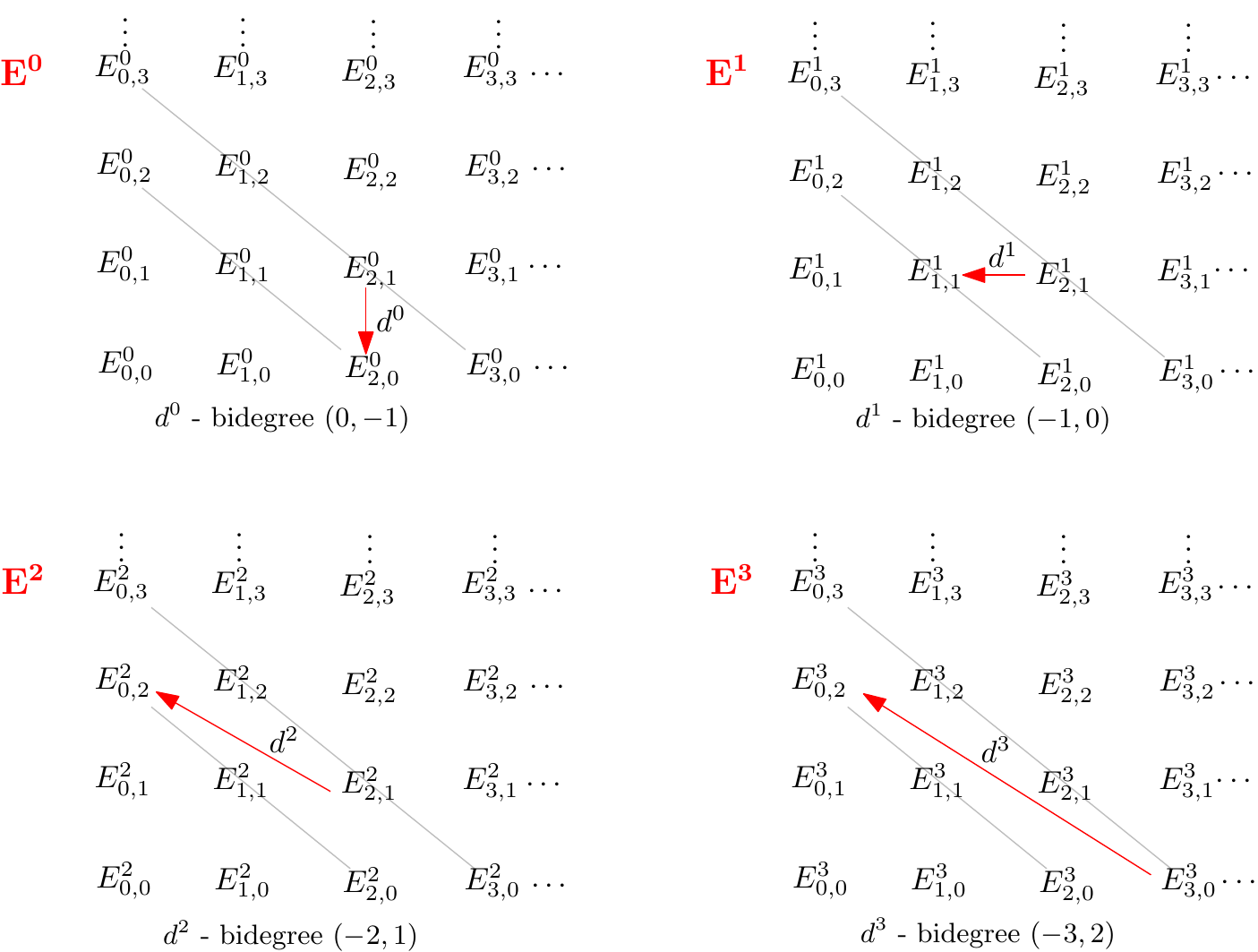}
\caption{Spectral sequence pages.}\label{fig:spectralsequencepages}
\end{figure}

In this paper, one considers a specific filtration in a Morse complex $(C,\Delta)$ defined by  $f$. Given a  finest  Morse decomposition $\{ M(p) \mid p\in P=\{1,\ldots, m \}, \  m = \# Crit(f)  \}$ such that there are distinct critical values $c_p$ with $f^{-1}(c_p)\supset M(p)$,  we can define a filtration on $M$  by $$\{F_{p-1}\}_{p=1}^{m}=\{f^{-1}(-\infty ,c_p+\epsilon)\}_{p=1}^{m}.$$ 
Let  $F=F_{p}C$ the induced filtration in $(C,\Delta)$.
Since for each $p\in P$ there is only one singularity in $F_p\setminus F_{p-1}$ the filtration $F$ is called a finest filtration.

With this in mind, we consider a filtered Morse chain complex $(C,\Delta)$ and calculate its spectral sequence. In this case the differential $\Delta$ is a connection matrix as proved in \cite{Salamon:1990}. Nonzero entries of this matrix imply the existence of connecting orbits, hence connection matrices provide information concerning global dynamics.
In this article, with a view towards application, we develop algorithms that model the spectral sequence of a filtered Morse chain complex  in the hope of coming across a relationship between algebraic cancellation of modules and dynamical cancellation of critical points. 
%
%
These algorithms  are intimately related to the Spectral Sequence Sweeping Algorithm {(SSSA)} developed in \cite{Bertolim:2013,Cornea:2010,Rezende:2010,Franzosa:2013}, and are introduced with the  intention of allowing one to gain a deeper understanding of the role of the differentials.

It was proved in \cite{Cornea:2010} that the Spectral Sequence Sweeping Algorithm applied to a connection matrix $\Delta$ generates a sequence of connection matrices $\Delta^r$ which keeps track of the differentials $d^r$ that cause the algebraic cancellations of modules $E^r$ of the spectral sequence, for $r\geq 0$. Of course, this immediately raises the question of understanding up to what point  algebraic cancellations in $(E^r,d^r)$ determine dynamical cancellations in a parameterized family of flows.

In \cite{Bertolim:2013}, a parameterized family of Morse flows\footnote{We mean by Morse flows those that arise from a gradient of a Morse function and satisfy the transversality condition. } on surfaces that undergo bifurcations was associated to this sequence of connection matrices $\Delta^{r}$. In other words, for the 2-dimensional case, the algebraic cancellations that occur as one ``turns the pages'' $E^r$ of the spectral sequence correspond to the dynamical cancellations that may occur in this parameterized family of flows. 
More specifically, the differentials of the spectral sequence that cause algebraic cancellations as $r$ increases are in 1-to-1 correspondence with the primary pivots of $\Delta^r$ determined by the Spectral Sequence Sweeping Algorithm. 
It is interesting to note that the connection matrices in the 2-dimensional setting are totally unimodular (TU). 

This instigated our inquiry into  the higher dimensional case under the hypothesis that the connection matrix $\Delta$ be totally unimodular. In this case, one can assert that for a filtered  Morse chain complex $(C,\Delta)$   on a closed simple connected manifold  $M$ of dimension $m>5$, with $\Delta$  a TU connection matrix, the algebraic cancellations of the modules of the spectral sequences  correspond to dynamical cancellations of consecutive critical points determined by the Row Cancellation Algorithm (RCA).

This result was first obtained in \cite{Bertolim:2013} with the proviso that these algorithms were compatible. Among other ingredients of this proof, two propositions were needed which are proved herein. The first being that the primary pivots were $\pm 1 $, and the second that  the primary pivots in both algorithms coincide in position and value, see  Corollaries \ref{cor:unitpivotsurface} and \ref{cor:smale=incr}. In attempting to prove this in the 2-dimensional case, we actually obtained  proofs in  dimension $n$, which are described in Theorems \ref{thm:unitpivotTU} and \ref{thm:smale=incr}  and enabled us to prove Theorem~\ref{ordering-cancellation-TU}. Hence, these constitute  the main results in this article, the first two are more algebraic in nature whereas  the last has a more dynamical flavor.

The main dynamical result in this article is the following.

\begin{thm} \label{ordering-cancellation-TU} Let $(C,\Delta)$ be a Morse chain complex associated to a  Morse-Smale function $f$ on a closed simply connected manifold $M$ of dimension $m>5$. Suppose that  $\Delta$ is totally unimodular. Let $(E^r,d^r)_{r\geq 0}$ be the associated spectral sequence for the finest filtration $F=\{F_pC\}$ defined by $f$. The algebraic cancellations of the modules $E^r$ of the spectral sequence determined by the SSSA are in one-to-one correspondence with  dynamical cancellations  that occur in Morse flows on $M$ associated to $(C,\Delta)$ determined by the RCA.
\end{thm}

\begin{cor} Let $(C,\Delta)$ be a Morse chain complex associated to a  Morse-Smale function $f$ on a closed simply connected manifold $M$ of dimension $m>5$. Suppose that $\Delta$ is a TU connection matrix then admits a perfect Morse function.
\end{cor}



 Several questions are in order here.  Of course,  the question of considering this result without the TU hypothesis for arbitrary connection matrices in dimension $n$ remains open.

The basis for this theorem is threefold. In  \cite{Cornea:2010}, it was proved that the primary pivots established by the SSSA  correspond to  differentials of the spectral sequence, which are responsible for  algebraic cancellations therein, by Theorem \ref{thm:unitpivotTU}. Secondly, in this article, we develop a row cancellation algorithm which reflects  dynamical cancellations. Finally, Theorem~\ref{thm:smale=incr} asserts that the primary pivots of both algorithms are the same. Hence, algebraic and dynamical cancellations coincide.

\begin{figure}[h!]
\centering
\includegraphics[width=10cm]{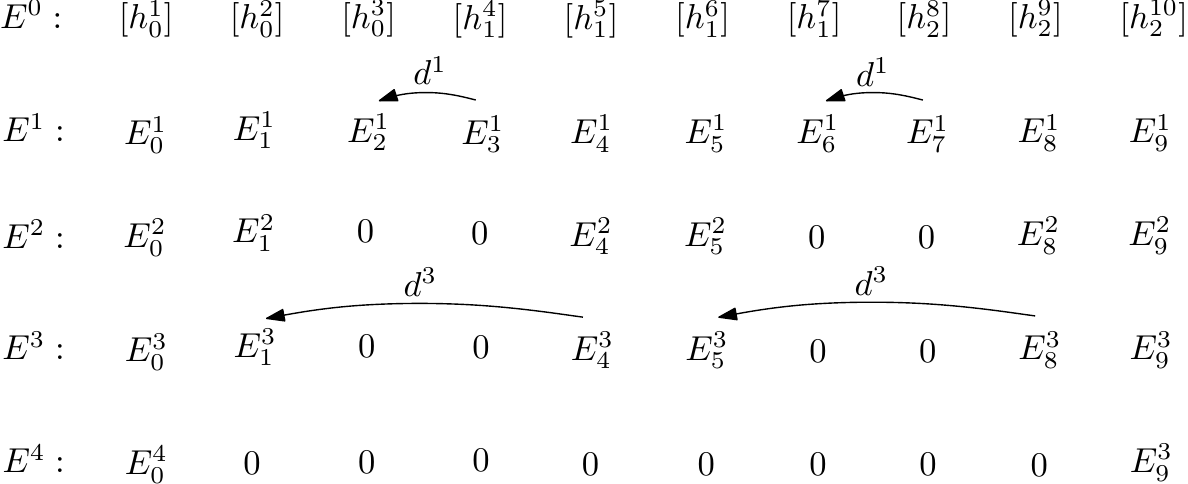}
\caption{Spectral sequence of a filtered Morse chain complex.}\label{fig:intro2}

\end{figure}

The sweeping algorithms developed in \cite{Cornea:2010,Rezende:2010} concerned connection matrices somewhat ordered in structure. This ordering is presently relaxed, and Section~\ref{sweeping} is devoted to adapting the sweeping algorithms to this general class of matrices considered here, as well as showing that the previously established properties are valid for this class. In Section~\ref{sec:TUsweeping}, we explore additional properties of the sweeping algorithm that occur when the input is restricted to TU connection matrices. Section~\ref{sec:rowcancellation} introduces a new algorithm, the row cancellation algorithm, which is better suited to a dynamical interpretation. Section~\ref{surfaces} builds the bridge between the algebraic cancellations of the sweeping algorithms and the dynamical cancellations of the row cancellation algorithm. In subsection \ref{subsec:dynamicTU}, Theorem \ref{ordering-cancellation-TU} is proved. With a view towards applications, we consider, in Subsection \ref{subsec:Smale}, the $2$-dimensional case. In fact, the restriction to surfaces of the row cancellation algorithm gives rise to Smale's Cancellation Sweeping Algorithm, which is yet another beautiful example linking dynamics and algebra, which is explored in \cite{Bertolim:2013}. See Figures  \ref{fig:intro1} and \ref{fig:intro2}.

\begin{figure}[h!]
\centering
\includegraphics[width=14cm]{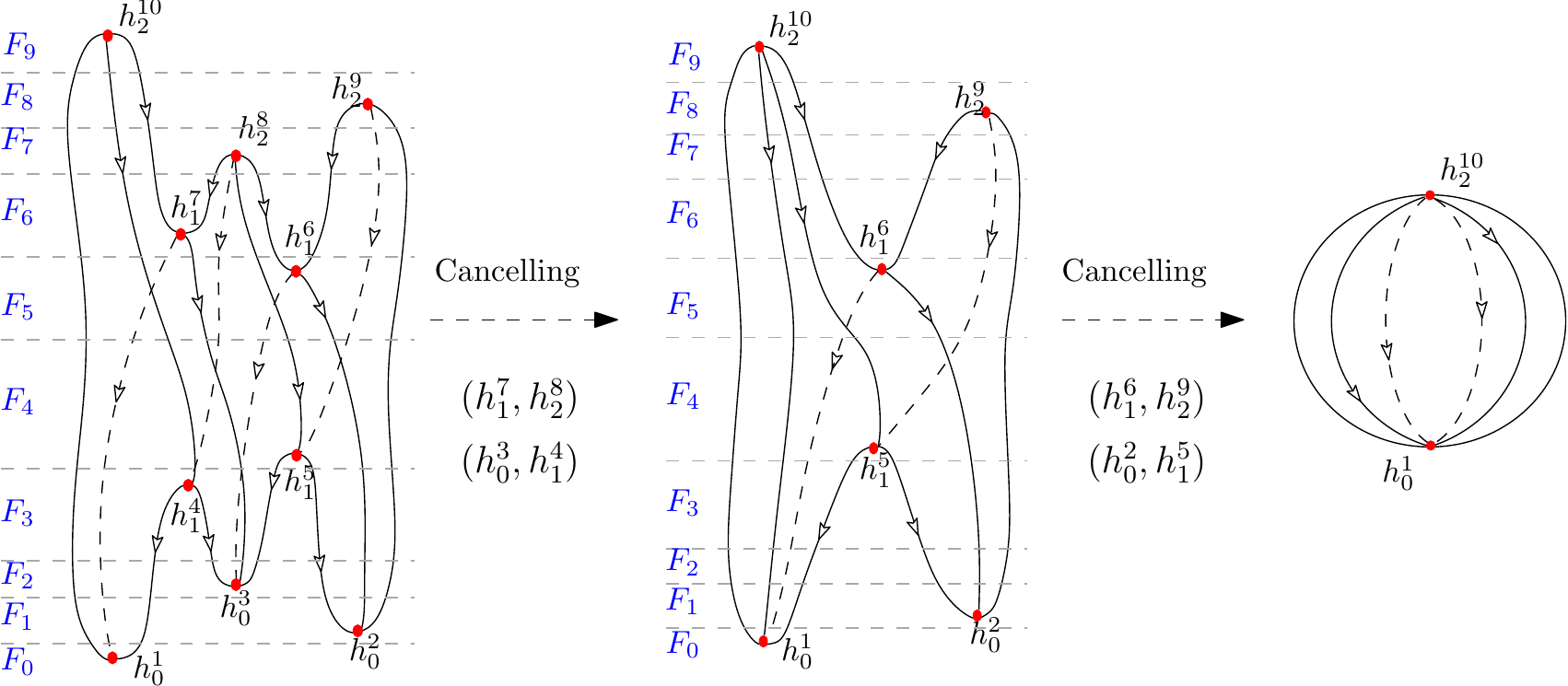}\\  \vspace{-0.3cm}
\caption{Continuation via cancellation of critical points.}\label{fig:intro1}

\end{figure}

\section{Sweeping algorithms}\label{sweeping}

In this section we present several sweepings algorithm for connection matrices. We adopt the following notation, namely, that the columns and rows of the order $m$ connection matrix $\Delta$ are partitioned into the subsets $J_0$, \ldots, $J_b$, columns and rows in $J_k$ are associated with elementary chains of index $k$. A connection matrix is grouped whenever $J_i$ is a set of consecutive integers, for all $i$, with the entries in $J_{k-1}\times J_k$ located above the diagonal, and ungrouped otherwise.

The motivation behind working with ungrouped connection matrices lies in the fact that the data of a simplicial complex is generally presented in a disorderly fashion, i.e., $k-$ simplices unordered with respect to a given filtration. Hence, columns representing index $k$-critical points appear with mixed indices from left to right. See Figure \ref{fig:scatteredmatrix}.

Notation regarding matrices is defined in Table~\ref{Table:notation}.
\begin{table}[htbp]
\centerline{\begin{tabular}{l|l}
$A_{i\textbf{\large .}}$ & $i$-th row of matrix $A$\\
$A_{\textbf{\large .}j}$ & $j$-th column of matrix $A$\\
$A_{I\textbf{\large .}}$ & submatrix of $A$ with entries\footnotemark\  $a_{ij}$ such that $i\in I$\\
$A_{IJ}$ & submatrix of $A$ with entries $a_{ij}$ such that $i\in I$ and $j\in J$,\\
 & where $I$ (resp., $J$) is a nonempty subset of the set of \\
 & row indices (resp., column indices)\\
 $A^\ell$ & $\ell$-th matrix in a sequence, non negative superscripts\\
 & \textbf{do not} denote exponents\\
 $(A^\ell)^{-1}$ & the inverse of matrix $A^\ell$
\end{tabular}}
\caption{Notation adopted for (sub)matrices.}\label{Table:notation}
\end{table}
\bigskip

\footnotetext{When there is no danger of ambiguity, the comma between row and column indices is omitted.}

The connection matrix $\Delta$ is an upper triangular nilpotent matrix with zero entries except, possibly, for the blocks $\Delta_{J_{k-1},J_k}$, for $k=1,\ldots, b$. In~\cite{Rezende:2010}, the connection matrices considered were grouped, i.e., all blocks were situated strictly above the main diagonal, and thus each subset in the partition contained consecutive indices. For example, Figure~\ref{fig:scatteredmatrix} shows two possible configurations for matrices with indices partitioned into four subsets. Both matrices have index subsets of same cardinality (3, 5, 2 and 2), but the one on the left has consecutive indices within each subset, whereas the one on the right does not. Positions in the first block are indicated with vertical lines, in the second with horizontal lines and in the third with diagonal lines. Shaded areas indicate positions that may have nonzero entries, the \emph{allowable sparsity pattern} of the connection matrix, namely positions in the set $\cup_{k=0}^{b-1} J_k\times J_{k+1}$ that lie strictly above the diagonal. Notice that, since the matrix must be upper triangular, the scattering reduces the number of positions eligible for nonzero entries. Thus, in the example on the left of Figure~\ref{fig:scatteredmatrix}, the connection matrix may have a total of $29$ nonzero entries, but the number of nonzero entries of the matrix on the right is at most $17$. We call a connection matrix with blocks above the main diagonal a \emph{grouped connection matrix}, to distinguish it from the (general) connection matrix  exemplified on the right of Figure~\ref{fig:scatteredmatrix}. Of course, the grouped connection matrix is a special type of connection matrix. The set of consecutive integers $\{i, i+1,\ldots, i+k\}$ is denoted $i..i+k$.

\begin{figure}[htbp]
\centerline{%
\begin{tikzpicture}
\begin{scope}[x=0.3cm,y=0.3cm,xscale=1.4,yscale=1.4]
\draw (1.5,12) node[above]{$\overbrace{\rule{1.2cm}{0pt}}^{J_0}$};
\draw (5.5,12) node[above]{$\overbrace{\rule{1.9cm}{0pt}}^{J_1}$};
\draw (9,12) node[above]{$\overbrace{\rule{.4cm}{0pt}}^{J_2}$};
\draw (11,12) node[above]{$\overbrace{\rule{.4cm}{0pt}}^{J_3}$};
\draw[color=DarkSlateGray4,fill=DarkSlateGray3,fill opacity=0.05] (0,0)--(0,12)--(12,12)--(12,0)--cycle;
\draw[color=DarkSlateGray4,fill=DarkSlateGray3,fill opacity=0.3] (3,12)--++(5,0)--++(0,-8)--++(4,0)--++(0,-2)--++(-2,0)--++(0,7)--++(-7,0)--cycle;
\pattern [pattern=vertical lines,pattern color=DarkSlateGray4] (3,12)--++(5,0)--++(0,-3)--++(-5,0)--cycle;
\pattern [pattern=horizontal lines,pattern color=DarkSlateGray4] (8,9)--++(2,0)--++(0,-5)--++(-2,0)--cycle;
\pattern [pattern=north east lines,pattern color=DarkSlateGray4] (10,4)--++(2,0)--++(0,-2)--++(-2,0)--cycle;
\draw (.1,10.5) node[rotate=90,above]{$\overbrace{\rule{1.2cm}{0pt}}$}
(-.8,10.5) node[left]{\scriptsize$J_0$}
(.1,6.5) node[rotate=90,above]{$\overbrace{\rule{1.9cm}{0pt}}$}
(-.8,6.5) node[left]{\scriptsize$J_1$}
(.1,3) node[rotate=90,above]{$\overbrace{\rule{.4cm}{0pt}}$}
(-.8,3) node[left]{\scriptsize$J_2$}
(.1,1) node[rotate=90,above]{$\overbrace{\rule{.4cm}{0pt}}$}
(-.8,1) node[left]{\scriptsize$J_3$};
\foreach \x in {0,...,11}
\draw[shift={(\x,-\x)},color= DarkSlateGray3]
(0,12)--(1,12)--(1,11)--(0,11)--cycle;
\draw (0,-2) node[right]{\parbox{6cm}{$J_0=1..3$, $J_1=4..8$,\\ $J_2=9..10$, $J_3=11..12$}};
\end{scope}
\begin{scope}[xshift=8cm,x=0.3cm,y=0.3cm,xscale=1.4,yscale=1.4]
\draw[color=DarkSlateGray4,fill=DarkSlateGray3,fill opacity=0.05] (0,0)--(0,12)--(12,12)--(12,0)--cycle;
\foreach \x in {0,...,11}
\draw[shift={(\x,-\x)},color= DarkSlateGray3]
(0,12)--(1,12)--(1,11)--(0,11)--cycle;
\foreach \i / \j in {6/2,6/4,12/2,12/4,12/7}
\fill [DarkSlateGray4,opacity=0.3]  (\i-1,13-\j) rectangle (\i,12-\j);
\foreach \j in {2,4,7}
\fill [DarkSlateGray4,opacity=0.3]  (9-1,13-\j) rectangle (10,12-\j);
\foreach \i  in {1,6,12}
	\foreach \j in {2,4,7}
	{\draw[color=DarkSlateGray4] (\i-1,13-\j) rectangle (\i,12-\j);
	\pattern [pattern=vertical lines,pattern color=DarkSlateGray4] (\i-1,13-\j) rectangle (\i,12-\j);
	}
\foreach \j in {2,4,7}
	{\draw[color=DarkSlateGray4] (9-1,13-\j) rectangle (10,12-\j);
	\pattern [pattern=vertical lines,pattern color=DarkSlateGray4] (9-1,13-\j) rectangle (10,12-\j);
	}
\foreach \i / \j in {3/1,8/1,8/6}
\fill [DarkSlateGray4,opacity=0.3]  (\i-1,13-\j) rectangle (\i,12-\j);
\foreach \i  in {3,8}
	\foreach \j in {1,6,12}
	{\draw[color=DarkSlateGray4] (\i-1,13-\j) rectangle (\i,12-\j);
	\pattern [pattern=horizontal lines,pattern color=DarkSlateGray4] (\i-1,13-\j) rectangle (\i,12-\j);
	}
\foreach \j in {3,8}
	{\draw[color=DarkSlateGray4] (\j-1,13-9) rectangle (\j,12-10);
	\pattern [pattern=horizontal lines,pattern color=DarkSlateGray4] (\j-1,13-9) rectangle (\j,12-10);
	}
\foreach \i / \j in {5/3,11/3,11/8}
\fill [DarkSlateGray4,opacity=0.3]  (\i-1,13-\j) rectangle (\i,12-\j);
\foreach \i  in {5,11}
	\foreach \j in {3,8}
	{\draw[color=DarkSlateGray4] (\i-1,13-\j) rectangle (\i,12-\j);
	\pattern [pattern=north east lines,pattern color=DarkSlateGray4] (\i-1,13-\j) rectangle (\i,12-\j);
	}
\foreach \i in {1,2,3,...,12}
{\draw (-.1,12-\i+.5) node[left]{\footnotesize\i}
			(\i+.2,12.5) node[left]{\footnotesize\i};}
\draw (0,-2) node[right]{\parbox{6.5cm}{$J_0=\{2,4,7\}$, $J_1=\{1,6,9,10,12\}$, $J_2=\{3,8\}$, $J_3=\{5,11\}$}};
\end{scope}
\end{tikzpicture}%
}
\caption{Two possible configurations of matrices with column/row sets partitioned into 4 subsets.}\label{fig:scatteredmatrix}
\end{figure}

The sweeping algorithm over $\Z$ was introduced in \cite{Cornea:2010} and further explored in \cite{Rezende:2010}. It was stated in terms of grouped connection matrices, and we rewrite it below, making the necessary slight notation changes to encompass the general case.\\

\textbf{\large Sweeping Algorithm over $\Z$}

\begin{description}
\item[\textbf{Input:}] nilpotent $m\times m$ upper triangular matrix $\Delta$ with column/row partition $J_0$, $J_1$, \ldots, $J_b$.
\item[\textbf{Initialization Step:}]\mbox{}\\
 $\begin{array}{@{}l}
\left[\begin{tabular}{l}
 $r=0$\\
 $\Delta^r=\Delta$\\
 $P^r=I$ ($m\times m$ identity matrix) \end{tabular}\right.\end{array}$
\item[\textbf{Iterative Step:}] (Repeated until all diagonals parallel and to the right of the main diagonal have been swept)\\
$\begin{array}{@{}l}
\left[\begin{tabular}{l}
\textbf{Matrix $\Delta$ update}\\
\begin{tabular}{@{\hspace{.5cm}}l}
$r\leftarrow r+1$ \\
$\Delta^r = (P^{r-1})^{-1} \Delta^0P^{r-1}$\\
\end{tabular}\end{tabular}\right.\\
\end{array}$
\\[5pt]
$\begin{array}{@{}l}
\left[\begin{tabular}{l}
\textbf{Markup}\\
\begin{tabular}{@{\hspace{.5cm}}l}
Sweep entries of $\Delta^r$ in the $r$-th diagonal: \\
\textbf{\textsf{If}} $\Delta^r_{j-r,j}\neq 0$ \textbf{\textsf{and}} $\Delta^r_{\textbf{\large .},j}$ does not contain a primary pivot\\
\rule{.5cm}{0pt}\textbf{\textsf{Then If}} $\Delta^r_{j-r\textbf{\large .}}$ contains a primary pivot\\
\rule{2cm}{0pt}\textbf{\textsf{Then}} temporarily mark $\Delta^r_{j-r,j}$ as a change-of-basis pivot\footnotemark\\
\rule{2cm}{0pt}\textbf{\textsf{Else}} permanently mark $\Delta^r_{j-r,j}$ as a primary pivot\\
\end{tabular}\end{tabular}\right.\\
\end{array}$
\\[5pt]
$\begin{array}{@{}l}
\left[\begin{tabular}{l}
\textbf{Matrix $P^r$ computation}\\
\begin{tabular}{@{\hspace{.5cm}}l}
$P^r \leftarrow P^{r-1}$\\
\textbf{\textsf{For each}} change-of-basis pivot $\Delta^r_{j-r,j}$,
update the $j$-th column of $P^r$ as follows\\
\begin{tabular}{@{\hspace{.5cm}}l}
Let $k$ be the index of the chain associated with column $j$ of $\Delta$.\\
Let $I=J_{k-1}\cap\{j-r,\ldots, m\}$, $J=J_k\cap\{1,\ldots, j\}$, $c=|J|$\\
Let $x^*\in \Z^c$ be an optimal solution to \\%
$\begin{array}{@{\hspace{3cm}}crcl} \min & x_c\\
\mbox{subject to} &\Delta_{IJ}x &= & 0\\
&x_c&\geq &1\\
&x&\in&\Z^c\end{array}$\\
$P^r_{Jj}\leftarrow x^*$
\end{tabular}
\end{tabular}\end{tabular}\right.\\
\end{array}$

\item[\textbf{Final Step:}]\mbox{} \\
$\begin{array}{@{}l}\left[
\begin{tabular}{l}
\textbf{Matrix $\Delta$ update}\\
\begin{tabular}{@{\hspace{.5cm}}l}
$r\leftarrow r+1$ \\
$\Delta^r = (P^{r-1})^{-1} \Delta^0 P^{r-1}$\\
\end{tabular}
\end{tabular}\right.
\end{array}$

\item[\textbf{Output:}] $(\Delta^0, \ldots, \Delta^m)$ and $(P^0,\ldots, P^{m-1})$
\end{description}
\footnotetext{Temporary marks are erased at the end of the iterative step.}

The connection matrix can be thought of as the matrix of a linear operator with respect to the basis $h=(h^1_{i_1},\ldots, h^{m}_{i_m})$, where $i_j$ denotes the index of the chain associated with column $j$. The sweeping algorithm produces a sequence of similar matrices $\Delta^0(=\Delta)$, $\Delta^1$, \ldots, $\Delta^m$. The basis associated with $\Delta^r$ is $\sigma^r=(\sigma^{1,r}_{i_1},\ldots, \sigma^{m,r}_{i_m})$. Initially, we have $\sigma^0=h$. At iteration $r$ of the sweeping algorithm, the entries on the $r$-th diagonal of $\Delta^r$ are swept. Zero entries are left unmarked. A nonzero entry at position $(j-r,j)$ is left unmarked if there is a primary pivot below it, is marked temporarily as a change-of-basis pivot if there is a primary pivot to its left, and is otherwise permanently marked as a primary pivot. If $\Delta^r_{j-r, j}$ is a change-of-basis entry, then the element $\sigma^{r+1}_j$ in the new basis $\sigma^{r+1}$ is an integer linear combination of the elements of $h$ with the same chain index, say $k$, and associated with columns to the left of, and including, column $j$, constructed so that the following two conditions are satisfied: (1) the entries on and below position $(j-r,j)$ in $\Delta^{r+1}$ are zero, and (2) the coefficient of $h^j_k$ in this linear combination, called the \emph{leading coefficient}, is the smallest positive integer that allows (1) to happen.

The correctness of the sweeping algorithm over $\Z$ rests on Proposition~\ref{prop:CorrectnessAlgorithmOverZ}, an extension of Proposition 8 of \cite{Rezende:2010} that encompasses both general and grouped connection matrices. Before proceeding, we introduce a technical lemma that facilitates the extension of the results previously developed for grouped connection matrices.
 
 \begin{lem}\label{lema:tecnico} Let $N$ be an upper triangular $m\times m$ matrix with nonzero diagonal entries and column/row partition $J_0$, \ldots, $J_b$, such that $n_{ij}=0$ if $(i,j)\notin J_k\times J_k$, for  $k=0,\ldots,b$. Then $N$ is invertible, its inverse is upper triangular and also has zero entries outside the set $\cup_{k=0}^b J_k\times J_k$. Furthermore,
 \[ (N_{J_kJ_k})^{-1} = (N^{-1})_{J_kJ_k},\quad \mbox{ for }k=0,\ldots, b. \]
 Additionally, if $j\in J_k$, for some $k\in\{0,\ldots, b\}$, $I=J_k\cap\{j,\ldots,m\}$ and $J=J_k\cap\{1,\ldots, j\}$, $\overline{I}$ and $\overline{J}$ denote their respective complements with respect to $\{1,\ldots,m\}$, then
\[(N^{-1})_{I\overline{I}} = 0\quad\mbox{and}\quad (N^{-1})_{\overline{J} J} = 0.\]
 \end{lem}

 \dem If $N$ is upper triangular, its determinant is the product of the diagonal entries. Since these entries are nonzero by hypothesis, $N$ is nonsingular and its inverse is upper triangular. Let the column/row partition be given by
\[\begin{array}{rcl}
J_0&=&\{i_1,\ldots,i_{|J_0|}\},\quad\mbox{where }i_1<\cdots<i_{|J_0|}\\
J_1&=&\{i_{|J_0|+1},\ldots,i_{|J_0|+|J_1|}\},\quad\mbox{where }i_{|J_0|+1}<\cdots<i_{|J_0|+|J_1|},\\
&\vdots&\\
J_b&=&\{i_{|J_0|+\cdots+|J_{b-1}|+1},\ldots,i_{|J_0|+\cdots+|J_b|}\},\quad\mbox{where }i_{|J_0|+\cdots+|J_{b-1}|+1}<\cdots<i_{|J_0|+\cdots+|J_{b}|},
\end{array}\]
and let $p\colon \{1,\ldots, m\}\to\{1,\ldots,m\}$  be the permutation that ``orders'' the columns according to the subsets, while maintaining the order inside each subset, that is,
\[\begin{array}{rcl}
p(i_1) &=&1,\\
p(i_2) &=& 2,\\
&\vdots&\\
p(i_{|J_0|+\cdots+|J_{b}|}) &=& m.
\end{array}\]
This definition implies $p(J_k)$ is the set of consecutive integers $|J_0|+\cdots+|J_{k-1}|+1..|J_0|+\cdots +|J_k|$.
Let $Q$ be the permutation matrix that embodies the column interchanges defined by $p$:
\[ Q_{\textbf{\large .}j} = e_{i_j}, \quad \mbox{for }j=1,\ldots,m,\]
where the $m$-column vector $e_{i_j}$ is zero except for its $i_j$-th entry, which is one. Then
\begin{eqnarray*} 
(Q^{-1}N Q)_{rc} &=& (Q^T N Q)_{rc}\\
&=& e^T_{i_r} N e_{i_c}\\
&=& n_{i_r i_c},\end{eqnarray*}
which implies
\[(Q^{-1}NQ)_{p(J_k)\, p(J_k)} = N_{J_k J_k}.\]
Thus $Q$ has the effect of grouping entries in each block. Furthermore, if $i_r$ and $i_c$ belong to the same subset in the column/row partition, say $J_k$, then $i_r<i_c$ (resp., $i_r=i_c$) if and only if $p(i_r)<p(i_c)$ (resp., $p(i_r)=p(i_c)$), so $Q^T N Q$ has a block diagonal structure, with upper triangular diagonal blocks $(Q^T N Q)_{p(J_k) \,p(J_k)}$, $k=1,\ldots, b$, each with nonzero diagonal elements. This implies that the inverse of $Q^T N Q$ has the same structure, where each diagonal block is replaced by its inverse. Entries of $(Q^T N Q)^{-1}$ outside these blocks are zero, therefore entries in $N^{-1}$ outside $\cup_{k=1}^b J_{k-1}\times J_k$ are zero. In particular, this implies that, for $k=0,\ldots, b$,
\begin{eqnarray*}
(N_{J_k\, J_k})^{-1} &=& ((Q^T N Q)_{p(J_k)\,p(J_k)})^{-1}\\
&=&(Q^T N Q)^{-1}_{p(J_k) \,p(J_k)}\\
&=& (Q^T N^{-1} Q)_{p(J_k) \,p(J_k)}\\
&=& (N^{-1})_{J_k\, J_k}.
\end{eqnarray*}

Additionally, we have that,
\[(N^{-1})_{I \overline{I}} =(Q^T N Q)^{-1}_{p(I) \,p(\overline{I})}.\]
The block diagonal structure of $Q^T N Q$ implies that $(Q^T N Q)_{p(I),p(\overline{J}_k)}=0$. Let $L=J_k\cap\{1$,\ldots, $j-1\}$. Then $(Q^T N Q)^{-1}_{p(I)\, p(L)}=0$ because the block $(Q^T N Q)^{-1}_{p(J_k)\, p(J_k)}$ is upper triangular. The result 
\[(N^{-1})_{I\overline{I}} = 0\]
follows from the fact $\overline{I}=L\cup\overline{J}_k$. The remaining result follows analogously.  \qed

\begin{prop}\label{prop:CorrectnessAlgorithmOverZ} Let $(\Delta^0, \Delta^1,\ldots)$ and $(P^0,P^1,\ldots)$ be the sequence of connection and change-of-basis matrices, respectively, produced by the application of the Sweeping Algorithm over $\Z$ to the connection matrix $\Delta\in \Z^{m\times m}$ with column/row partition $J_0$, \ldots, $J_b$. Then, for $r\geq 1$, we have
\begin{itemize}
\item[(i)] $\Delta^r$ and $\Delta^0 P^{r-1}$ is compliant with the allowable sparsity pattern of $\Delta$;
\item[(ii)] the nonzero entries of $\Delta^0 P^{r-1}$ strictly below the $r$-th diagonal  are located on or above a unique primary pivot position, entries in primary pivot positions are nonzero;
\item[(iii)] the nonzero entries of $\Delta^r$ strictly below the $r$-th diagonal  are either primary pivots (always nonzero) or are above a unique primary pivot;
\item[(iv)] each linear integer program formulated in the matrix $P^{r-1}$ computation step has an optimal solution;
\item[(v)] the change-of-basis matrices $P^{r-1}$ have the following upper triangular structure: the diagonal elements are nonzero, the positions of the off-diagonal nonzero elements are contained in $(\cup_{i=1}^b J_i\times J_i) \cap \{(i,j)\in \{1,\ldots, m\}^2\mid i< j\}$.\end{itemize}

\end{prop}

\dem For $r=1$, items (i)--(v) are either trivially true, since $P^0=I$, which implies $\Delta^1=\Delta^0=\Delta$, whose diagonal is zero, or vacuously true. The proof is by induction.

Assume by induction that the items are true for fixed arbitrary $r$, greater than $1$. First we show that (iv) and (v) are satisfied for $r+1$. Consider the sweeping of the $r$-th diagonal and let the entry in position $(j-r,j)$ be marked as a change-of-basis pivot. Let $k$ be the chain index associated with column $j$. Let $I=J_{k-1}\cap \{j-r,\ldots,m\}$, $J=J_k\cap \{1,\ldots,j\}$, $c=|J|$.  Using Lemma~\ref{lema:tecnico} and the facts that the induction hypothesis that (v) is satisfied for $r$, we have that
\begin{eqnarray*}\Delta^r_{IJ} &=& (P^{r-1})^{-1}_{I\textbf{\large .}} \Delta^0 P^{r-1}_{\textbf{\large .}J}\\
&=&(P^{r-1})^{-1}_{II} \Delta^0_{IJ}P^{r-1}_{JJ}.
\end{eqnarray*}

We will construct a column vector $y\in\R^c$ such that $\Delta^r_{IJ}y=0$. Let $p$ be the column containing the primary pivot entry to the left of $\Delta^r_{j-r,j}$. Then, by induction hypothesis (i), $j-r\in I$ and $p\in J$. Additionally, by induction hypothesis (iii), the entries below row $j-r$ in columns $p$ and $j$ of $\Delta^r$ are zero. Now let $y$ be a zero vector except for the entries corresponding to columns $\Delta^r_{Ip}$, whose value is $-\Delta^r_{j-r,j}$, and $\Delta^r_{Ij}$ (the $c$-th, or last, entry of $y$), whose value is $\Delta^r_{j-r,p}$. Then
\[\Delta^r_{IJ} y=0.\] 

The integrality of $P^{r-1}$ and $\Delta^0$ imply that $y$ is rational, and thus also $P^{r-1}_{JJ}y$. The facts that $P^{r-1}_{JJ}$ is upper triangular with nonzero diagonal entries and that $y_c\neq0$ imply that $(P^{r-1}_{JJ}y)_c\neq 0$. Thus there is a suitable multiple $\bar{x}$ of $y$ that is an integral solution to 
\[\left\{\begin{array}{ll} (P^{r-1})^{-1}_{II} \Delta^0_{IJ} x=0,\\x_c\geq 1,\end{array}\right.\]
and thus an integral solution to 
\[\left\{\begin{array}{ll} \Delta_{IJ} x=0,\\x_c\geq 1.\end{array}\right.\]
Therefore the linear integral problem in the matrix $P^r$ computation step is feasible. Its objective is also bounded by construction, since $x_c\geq 1$. We may conclude, see \cite{Nemhauser:1988}, that this linear integral problem has an optimal solution $x^*$. Hence (iv) is  for $P^r$. Furthermore, the facts that $I\subset J_{k-1}$ and $J\subset J_k$, and the assignment rule $P^r_{Jj}\leftarrow x^*$ imply that (v) is satisfied for $P^r$.

Now consider (i) and (ii) for $\Delta^0 P^r$. By induction, (i) and (ii) are satisfied for $\Delta^0P^{r-1}$. By construction of $P^r$, the only columns that change from $\Delta^0 P^{r-1}$ to $\Delta^0 P^r$ are the ones containing change-of-basis entries in $\Delta^r$. The corresponding columns in $P^r$ are built via the solution of the linear programs in the matrix $P^r$ computation step. This guarantees that the entries on and below the change-of-basis position in the product $\Delta^0P^r$ are zero. Since the change-of-basis position is on the $r$-th diagonal, it is strictly below the $(r+1)$-th diagonal, so (ii) is satisfied for these columns. Also, since we are replacing the column containing the change-of-basis entry with a linear combination of this column and other columns to its left, all with same chain index, the allowable sparsity pattern of $\Delta$ is preserved, so (i) is also satisfied for these columns. Now consider a column $s$ that was not changed. Trivially, (i) is satisfied in this case. To establish (ii), we consider the three possibilities for $(\Delta^0P^r)_{\{s-r,\ldots,m\},s}=(\Delta^0 P^{r-1})_{\{s-r,\ldots,m\},s}$. If $(\Delta^0 P^{r-1})_{\{s-r+1,\ldots,m\},s}\neq 0$, (ii) is valid by induction. If $(\Delta^0P^{r-1})_{\{s-r,\ldots,m\},s}=0$, (ii) is trivially satisfied. The remaining possibility is $(\Delta^0P^{r-1})_{s-r,s}\neq 0$ and $(\Delta^0P^{r-1})_{\{s-r+1,\ldots,m\},s}=0$. In this case, using induction and (v),
\[\Delta^r_{s-r,s} =(P^{r-1})^{-1}_{s-r,\{s-r,\ldots,m\}}(\Delta^0P^{r-1})_{\{s-r,\ldots,m\},s}= (P^{r-1})^{-1}_{s-r,s-r}(\Delta^0P^{r-1})_{s-r,s} \neq 0,\]
whereas, for $i>s-r$,
\[\Delta^r_{i,s} =(P^{r-1})^{-1}_{i,\{i,\ldots,m\}}(\Delta^0P^{r-1})_{\{i,\ldots,m\},s}= 0.\]
But these facts imply that $\Delta^r_{s-r,s}$ must have been marked as a primary pivot when the $r$-th diagonal was swept. Thus the unique nonzero element of $\Delta^0P^r$ in column $s$, strictly below the $(r+1)$-th diagonal, is located on a primary pivot position and (ii) is satisfied.

Next we show (i) and (iii) are true for $\Delta^{r+1}$. Given the (established above) upper triangular structure of $P^r$, shared by its inverse, if $i\in J_k$, then
\[\Delta^{r+1}_{i\textbf{\large .}} = \sum_{i'\in J_k\cap\{i+1,\ldots,m\}}P^r_{i, i'} (\Delta^0P^r)_{i'\textbf{\large .}}.\] 
Then the fact that $\Delta^0P^r$ satisfies (i), imply $\Delta^{r+1}$ also satisfies. Furthermore, the fact that $\Delta^0P^r$ satisfies (ii) implies that the nonzero entries strictly below the $(r+1)$-th diagonal are located on or above a primary pivot position. Now suppose position $(i,\ell)$, strictly below the $(r+1)$-th diagonal, is marked as a primary pivot. Then, using (ii), 
\[\Delta^{r+1}_{i\ell} = \sum_{i'\in J_k\cap\{i+1,\ldots,m\}}P^r_{i\, i'} (\Delta^0P^r)_{i'\ell}= P^r_{i\,i} (\Delta^0P^r)_{i\ell}\neq 0,\]
so the primary pivot entries in $\Delta^{r+1}$ are nonzero. Finally, algorithm rules allow at most one primary pivot per column, so the uniqueness in (iii) is easily guaranteed. \qed

The following corollary addresses the configuration of the last matrix obtained by the Sweeping Algorithm over $\Z$ and is a simple consequence of Proposition~\ref{prop:CorrectnessAlgorithmOverZ}.

\begin{cor}\label{cor:zeropattern3} Let $\Delta^m$ be the last matrix produced by the application of the Sweeping Algorithm over $\Z$ to the connection matrix $\Delta\in \Z^{m\times m}$. Then the primary pivot entries are nonzero and each nonzero entry is located above a unique primary pivot.
\end{cor}

As in the special case of grouped connection, see \cite{Rezende:2010}, the configuration established in Corollary~\ref{cor:zeropattern3} of the last connection matrix in the sequence produced by the algorithm leads to the complementary relation between columns and rows expressed in the next proposition. The proof is a straightforward adaption of the proof in \cite{Rezende:2010}.

\begin{prop}\label{prop:compl} Let $\Delta^m$ be the last matrix produced by the application of the Sweeping Algorithm over $\Z$ to the connection matrix $\Delta\in \Z^{m\times m}$. If the $j$-th column of $\Delta^m$ is nonzero, then its $j$-th row is null, or, equivalently,
\begin{equation}\label{eq:complF}
\Delta^m_{\textbf{\large .}j}\; \Delta^m_{j\textbf{\large .}}= 0, \quad \mbox{for all }j.
\end{equation}
\end{prop}

\dem Equation (\ref{eq:complF}) is trivial when $\Delta^m_{\textbf{\large .}j}$ is a zero column. So suppose $\Delta^m_{\mbox{\large\bf.}j}\neq 0$. By the inherited allowable sparsity pattern established in Proposition \ref{prop:CorrectnessAlgorithmOverZ}, there exists $s$ such that $j\in J_s$. By Corollary~\ref{cor:zeropattern3}, the nonzero columns of $\Delta^m$ are precisely the columns containing primary pivots. Label the primary pivots of columns in $J_s$ in increasing order of row index: if $\Delta^m_{i_1j_1}$, \ldots, $\Delta^m_{i_a j_a}$ are the primary pivots in columns in $J_s$, then $i_1<i_2<\cdots<i_a$. Thus, $j_1$, $j_2$, \ldots, $j_a$ are the nonzero columns in $\Delta^m_{\textbf{\large .}J_s}$ and $j\in \{j_1,\ldots, j_a\}$. Furthermore, $\Delta^m_{i_a j_a}$ is the unique nonzero entry of row $\Delta^m_{i_a \textbf{\large.}}$, row $i_{a-1}$ has a nonzero entry in column $j_{a-1}$ and may have another one in column $j_a$, and so on.

The fact that $\Delta^m$ is nilpotent implies that
\[ 0 = \Delta^m_{i_a\textbf{\large .}} \Delta^m_{\textbf{\large .}j'} = \Delta^m_{i_a j_a}\Delta^m_{j_a j'},\]
for fixed arbitrary $j'$.
Using the fact that the primary pivot entry is nonzero, we conclude that $\Delta^m_{j_a j'}=0$, for all $j'$. Repeating the argument for $i_{a-1}$ and using the fact that the $j_a$-th row of $\Delta^m$ is null, we establish that its $j_{a-1}$-th row is null. The nullity of rows $j_{a-2}$, \ldots, $j_1$ of $\Delta^m$ follow analogously. Therefore, we conclude that $\Delta^m_{\mbox{\bf\large .}j}=0$.\qed

If the connection matrix has entries in a field $\F$, the sweeping algorithm can be easily adapted by letting the variables in the minimization problem in the \emph{Matrix $P^r$ computation} step have entries in $\F$, and not in $\Z$. Let us call this version the \emph{Accumulated Sweeping Algorithm over $\F$}. The adjective ``accumulated'' refers to matrix $P^r$, which accumulates all information regarding the successive basis changes along iterations $1$, \ldots, $r$. The adaptation of Proposition \ref{prop:CorrectnessAlgorithmOverZ} to this new algorithm is presented below.

\begin{prop}\label{prop:CorrectnessAccumulatedAlgorithmOverF} Let $(\Delta^0, \Delta^1,\ldots)$ and $(P^0,P^1,\ldots)$ be the sequence of connection and change-of-basis matrices, respectively, produced by the application of the Accumulated Sweeping Algorithm over $\F$ to the connection matrix $\Delta\in \F^{m\times m}$ with column/row partition $J_0$, \ldots, $J_b$. Then, for $r\geq 1$, we have
\begin{itemize}
\item[(i)] $\Delta^r$ and $\Delta^0 P^{r-1}$ is compliant with the allowable sparsity pattern of $\Delta$;
\item[(ii)] the nonzero entries of $\Delta^0 P^{r-1}$ strictly below the $r$-th diagonal  are located on or above a unique primary pivot position, entries in primary pivot positions are nonzero;
\item[(iii)] the nonzero entries of $\Delta^r$ strictly below the $r$-th diagonal  are either primary pivots (always nonzero) or are above a unique primary pivot;
\item[(iv)] the optimization program formulated in the matrix $P^{r}$ computation step corresponding to change-of-basis pivot $\Delta^r_{j-r,j}$, where $j\in J_k$, has an optimal solution $x=P^{r-1}_{JJ}y$, where $y\in \F^c$ is given by
\begin{equation}\label{eq:optimalxF}
y_s = \left\{\begin{array}{ll}-\dfrac{\Delta^r_{j-r,j}}{\Delta^r_{j-r,p}},\quad&\mbox{if }s=p',\\
1,&\mbox{if }s=c,\\
0,&\mbox{otherwise,}\end{array}\right.
\end{equation}
where $p'$ is the column of $\Delta^r_{IJ}$ that contains $\Delta^r_{j-r,p}$,  the primary pivot to the left of the change-of-basis pivot $\Delta^r_{j-r,j}$;
\item[(v)] the change-of-basis matrices $P^{r-1}$ have the following upper triangular structure: the diagonal elements are equal to $1$, the positions of the remaining nonzero elements are contained in $(\cup_{i=1}^b J_i\times J_i) \cap \{(i,j)\in \{1,\ldots, m\}^2\mid i<j\}$.\end{itemize}

\end{prop}

\dem The only items that need to be proved are (iv) and (v). For the remaining ones,  the proof of Proposition \ref{prop:CorrectnessAlgorithmOverZ} is valid.

When $r=1$ item (iv) is vacuously true and (v) is trivially true, since $P^1$ is the identity matrix. Assume they are true for fixed arbitrary $r$, greater than $1$. Consider the case $r+1$, that is, matrix $P^r$ will be constructed based on the sweeping of the $r$-th diagonal of $\Delta^r$. Suppose $\Delta^r_{j-r,j}$ is a change-of-basis pivot and $\Delta^r_{j-r,p}$ is the primary pivot entry to its left. Suppose $j\in J_k$ and let $I$ and $J$ be as defined in the algorithm. Let $p'$ be the column of $\Delta^r_{IJ}$ that contains the primary pivot entry $\Delta^r_{j-r,p}$. By (iii)  and the rules of the algorithm, we have
\[(\Delta^r_{IJ})_{\textbf{\large .}p'} = \left[\begin{array}{c} \Delta^r_{j-r,p}\\ 0 \\\vdots \\ 0\end{array}\right]  \ \ \  \text{and} \ \ \    (\Delta^r_{IJ})_{\textbf{\large .}c} = \left[\begin{array}{c} \Delta^r_{j-r,j}\\ 0 \\\vdots \\ 0\end{array}\right].\]
Furthermore, by the induction hypotheses, the vector $y$ in (iv) is well defined and belongs to $\F^c$. Therefore
$\Delta^r_{IJ} y = 0.$

Let 
$x=P^{r-1}_{JJ}y.$
Since $\Delta^r_{IJ}=(P^{r-1})^{-1}_{II} \Delta^0_{IJ} P^{r-1}_{JJ}$, we have that $x$ satisfies $\Delta_{IJ} x = 0$. Furthermore, by the induction hypothesis (v),
$x_c = (P^{r-1}_{JJ})_{c\textbf{\large .}} y=y_c=1,$
which proves (iv) and (v), given that $P^r$ is constructed from $P^{r-1}$ by the replacement
$P^r_{Jj}=x,$
for each change-of-basis column $j$.\qed

Proposition~\ref{prop:CorrectnessAccumulatedAlgorithmOverF} implies we can simplify the matrix $P^r$ computation step in the Accumulated Sweeping Algorithm over $\F$, by using the ready-made optimal solution $x=P^{r-1}_{JJ}y$ as described in its item (iv). Now notice that, if $x=P^{r-1}_{JJ}y$, for $y$ described in item (iv) of the proposition, then
\[P^r_{Jj}=-\frac{\Delta^r_{j-r,j}}{\Delta^r_{j-r,p}} P^{r-1}_{Jp} + P^{r-1}_{Jj}.\]
Using simple algebra, we have that\\[-60pt]
\[P^r_{JJ} = P^{r-1}_{JJ}\;%
\setlength{\unitlength}{1.7pt}
\begin{picture}(50,46)(0,23)
\put(0,0){\line(0,1){44}}
\multiput(0,0)(0,44){2}{\line(1,0){2}}
\put(48,0){\line(0,1){44}}
\multiput(48,0)(0,44){2}{\line(-1,0){2}}
\put(4,2){%
\put(0,10){\line(1,0){30}}
\put(0,10){\line(0,1){30}}
\put(0,40){\line(1,0){40}}
\put(30,0){\line(0,1){40}}
\put(40,0){\line(0,1){40}}
\put(30,0){\line(1,0){10}}
\put(15,25){\makebox(0,0){$I$}}
\put(15,5){\makebox(0,0){$0~~\cdots~~0$}}
\put(35,20){\makebox(0,0){$y$}}}
\end{picture}.\]

\vspace{40pt}

\noindent Since Proposition~\ref{prop:CorrectnessAccumulatedAlgorithmOverF}(v) implies that entries in columns $P^{r-1}_{\textbf{\large .}p}$ and $P^{r-1}_{\textbf{\large .}j}$ outside rows in $J$ are zero, all updates relative to the sweeping of the $r$-th diagonal may be represented by one single matrix $T^r$, that coincides with the identity matrix, except for columns corresponding to change-of-basis columns. For each such column $j$, if $\Delta^r_{j-r,p}$ is the primary pivot to the left of the change-of-basis pivot $\Delta^r_{j-r,j}$, the matrix $T^r$ contains the entry $-\Delta^r_{j-r,j}/\Delta^r_{j-r,p}$ in position $(p,j)$. Summarizing, this matrix is given by
\begin{equation}\label{eq:defTr}
T^r_{ij}=\left\{\begin{array}{ll}
1, & \mbox{if }i=j,\\
-\dfrac{\Delta^r_{j-r,j}}{\Delta^r_{j-r,p}},\quad& \mbox{if }i=p, \Delta^r_{j-r,j} \mbox{ is a change-of-basis pivot and }\\
&\Delta^r_{j-r,p} \mbox{ is the primary pivot to its left,}\\
0,&\mbox{ow.}\end{array}\right.\end{equation}
With this definition, the accumulated change-of-basis $P^r$ satisfies
\begin{equation}\label{eq:Pr} P^r = P^{r-1} T^r.\end{equation}
Given that 
\begin{equation}\label{eq:Dr}
\Delta^{r+1} = (P^r)^{-1} \Delta^0 P^r = (P^{r-1} T^r)^{-1} \Delta^0 P^{r-1} T^r = (T^r)^{-1} \Delta^r T^r,\end{equation}
the algorithm may be simplified by dealing exclusively with the $T^r$ \emph{transition matrices}, thus avoiding the product in (\ref{eq:Pr}). To distinguish this simplified version from the previous one, we name it the \emph{Incremental Sweeping Algorithm over $\F$}, in the sense that in this case we keep a matrix that performs an incremental change, from $\Delta^{r}$ to $\Delta^{r+1}$, say, instead of a matrix that accumulates all changes, taking $\Delta^0$ directly to $\Delta^{r+1}$. This algorithm, inserted below for completeness, was introduced in~\cite{Rezende:2010}. Although it was defined and studied with grouped connection matrices in mind, we've seen that it may be equally applied to general connection matrices. In particular, the items of Proposition~\ref{prop:CorrectnessAccumulatedAlgorithmOverF} that concern the sequence $\Delta^0$, $\Delta^1$, \ldots, are satisfied by the Incremental Sweeping Algorithm over $\F$. 

The incremental way of updating the connection matrix is just a computational detail. The relevant difference between the accumulated and the incremental algorithms is that, in the former, any optimal solution to the optimization problem may be used, while, in the latter, a specific optimal solution is employed.\\

\textbf{\large Incremental Sweeping Algorithm over $\F$}

\begin{description}
\item[\textbf{Input:}] nilpotent upper triangular matrix $\Delta\in \F^{m\times m}$ with column/row partition $J_0$, $J_1$, \ldots, $J_b$.
\item[\textbf{Initialization Step:}]\mbox{}\\
 $\begin{array}{@{}l}
\left[\begin{tabular}{l}
 $r=0$\\
 $\Delta^r=\Delta$\\
 $T^r=I$ ($m\times m$ identity matrix) \end{tabular}\right.\end{array}$

\item[\textbf{Iterative Step:}] (Repeated until all diagonals parallel and to the right of the main diagonal have been swept)\\
$\begin{array}{@{}l}
\left[\begin{tabular}{l}
\textbf{Matrix $\Delta$ update}\\
\begin{tabular}{@{\hspace{.5cm}}l}
$r\leftarrow r+1$ \\
$\Delta^r = (T^{r-1})^{-1} \Delta^{r-1}T^{r-1}$\\
\end{tabular}\end{tabular}\right.\\
\end{array}$
\\[5pt]
$\begin{array}{@{}l}
\left[\begin{tabular}{l}
\textbf{Markup}\\
\begin{tabular}{@{\hspace{.5cm}}l}
Sweep entries of $\Delta^r$ in the $r$-th diagonal: \\
\textbf{\textsf{If}} $\Delta^r_{j-r,j}\neq 0$ \textbf{\textsf{and}} $\Delta^r_{\textbf{\large .},j}$ does not contain a primary pivot\\
\rule{.5cm}{0pt}\textbf{\textsf{Then If}} $\Delta^r_{j-r\textbf{\large .}}$ contains a primary pivot\\
\rule{2cm}{0pt}\textbf{\textsf{Then}} temporarily mark $\Delta^r_{j-r,j}$ as a change-of-basis pivot\footnotemark
\\
\rule{2cm}{0pt}\textbf{\textsf{Else}} permanently mark $\Delta^r_{j-r,j}$ as a primary pivot\\
\end{tabular}\end{tabular}\right.\\
\end{array}$
\\[5pt]
$\begin{array}{@{}l}
\left[\begin{tabular}{l}
\textbf{Matrix $T^r$ construction}\\
\begin{tabular}{@{\hspace{.5cm}}l}
$T^r \leftarrow I$\\
\textbf{\textsf{For each}} change-of-basis pivot $\Delta^r_{j-r,j}$,
change the $j$-th column of $T^r$ as follows\\
\begin{tabular}{@{\hspace{.5cm}}l}
Let $p$ be such that $\Delta^r_{j-r,p}$ is a primary pivot\\
$T^r_{pj}\leftarrow -\Delta^r_{j-r,j}/\Delta^r_{j-r,p}$
\end{tabular}
\end{tabular}\end{tabular}\right.\\
\end{array}$

\item[\textbf{Final Step:}]\mbox{} \\
$\begin{array}{@{}l}\left[
\begin{tabular}{l}
\textbf{Matrix $\Delta$ update}\\
\begin{tabular}{@{\hspace{.5cm}}l}
$r\leftarrow r+1$ \\
$\Delta^r = (T^{r-1})^{-1} \Delta^{r-1} T^{r-1}$\\
\end{tabular}
\end{tabular}\right.
\end{array}$

\item[\textbf{Output:}] $(\Delta^0, \ldots, \Delta^m)$ and $(T^0,\ldots, T^{m-1})$

\end{description}
\footnotetext{Temporary marks are erased at the end of the iterative step.}

Let $\mathcal{P}^r$ be the set of column indices of primary pivots of $\Delta^m$ on or below the $r$-th diagonal. Then $\mathcal{P}^0=\emptyset$ and $\mathcal{P}^r\backslash \mathcal{P}^{r-1}$ is the set of column indices of primary pivots on the $r$-th diagonal, for $r=1, \ldots, m-1$.  Proposition~\ref{prop:CorrectnessAccumulatedAlgorithmOverF}, the definition of $T^r$ given in (\ref{eq:defTr}), the relationship between the change-of-basis matrices given in (\ref{eq:Pr}) and the definition of the Incremental Sweeping Algorithm over $\F$ imply the following result.

\begin{prop}\label{prop:CorrectnessAlgorithmOverF} Let $(\Delta^0(=\Delta),\Delta^1,\ldots)$ and $(T^0,T^1,\ldots)$ be the sequence of connection and change-of-basis matrices produced by the Incremental Sweeping Algorithm over $\F$, where $\Delta$ is a connection matrix in $\F^{m\times m}$, with column/row partition $J_0$, \ldots, $J_b$. Then, for $r=1,\ldots, m$, we have
\begin{itemize}
\item[(i)] $\Delta^r$ is compliant with the allowable sparsity pattern of $\Delta$;
\item[(ii)] the nonzero entries of $\Delta^r$ strictly below the $r$-th diagonal are either primary pivots or are above a unique primary pivot;
\item[(iii)] primary pivot entries of $\Delta^r$ are nonzero;
\item[(iv)] the transition matrices $T^r$ have the following upper triangular structure: the diagonal elements are equal to $1$, each change-of-basis column has two nonzero entries and the positions of the off-diagonal nonzero elements are contained in \linebreak
$(\cup_{i=1}^b ((\mathcal{P}^{r-1}\backslash\mathcal{P}^{r-2})\cap J_i)\times ((\mathcal{P}^{r-1}\backslash\mathcal{P}^{r-2})\cap J_i)\cap \{(i,j)\in \{1,\ldots,m\}^2\mid i< j\}$.
\end{itemize}
\end{prop}

The next corollary and proposition follow from Proposition~\ref{prop:CorrectnessAccumulatedAlgorithmOverF}(ii) and mimic analogous results for the Sweeping Algorithm over $\Z$.

\begin{cor}\label{cor:zeropattern4} Let $\Delta^m$ be the last matrix produced by the application of the Incremental Sweeping Algorithm over $\F$ to the connection matrix $\Delta\in \F^{m\times m}$. Then the primary pivot entries are nonzero and each nonzero entry is located above a unique primary pivot.
\end{cor}

\begin{prop}\label{prop:complF} Let $\Delta^m$ be the last matrix produced by the application of the Incremental Sweeping Algorithm over $\F$ to the connection matrix $\Delta\in \F^{m\times m}$.
If the $j$-th column of $\Delta^m$ is nonzero, then its $j$-th row is null, or, equivalently,
\begin{equation}
\Delta^m_{\textbf{\large .}j}\; \Delta^m_{j\textbf{\large .}}= 0, \quad \mbox{for all }j.\label{eq:complZ}\end{equation}
\end{prop}

The change of basis that effects the transition from $\Delta^r$ to $\Delta^{r+1}= (T^r)^{-1}\Delta^r T^r$ is constructed so as to zero out in $\Delta^{r+1}$ the change-of-basis entries, without disturbing the (already null) entries below the change-of-basis pivots. If $\{j_1,\ldots, j_{t_r}\}$ denotes the set of distinct column indices that contain change-of-basis entries in $\Delta^r$ and $\{p_1,\ldots, p_{t_r}\}$ the column indices such that $(j_s-r,p_s)$ is the position of the primary pivot in the same row as the change-of basis pivot $\Delta^r_{j_s-r,j_s}$, for $s\in\{1,\ldots, t_r\}$, then the  matrix constructed in the iterative step, that effects this basis change, may be expressed as 
\begin{equation}\label{eq:matrixTr}T^r = I - \sum_{s=1}^{t_r} \frac{\Delta^r_{j_s-r,j_s}}{\Delta^r_{j_s-r,p_s}} U^{p_sj_s} =  I - \sum_{s=1}^{t_r}\alpha_{p_s j_s} U^{p_sj_s},\end{equation}
where $U^{ij}$ is the $m\times m$ matrix unit whose entries are all $0$ except in position $(i,j)$, where it is $1$. These matrix units have the property, see \cite{Lam:1999},
\begin{equation}\label{eq:matrixunitproperty}
U^{ij}U^{k\ell} = \delta_{jk}U^{i\ell}.\end{equation}

Algebraically, the post-multiplication of $\Delta^r$ by $T^r$ in (\ref{eq:matrixTr}) consists of $t_r$ elementary column operations on the change-of-basis columns of $\Delta^r$. Of the three possible elementary column operations on a column $j$, we use only one in this work: ``add to column $j$ a multiple of another column''. The property (\ref{eq:matrixunitproperty}) and the fact that no column of $\Delta^r$ may contain both a primary pivot mark and a change-of-basis pivot imply that
\begin{equation}\label{eq:Tr=product}
T^r=(I-\alpha_{p_1j_1}U^{p_1j_1})\cdots(I-\alpha_{p_{t_r} j_{t_r}}U^{p_{t_r}j_{t_r}}),\end{equation}
and the order of the matrix product on the right-hand-side of (\ref{eq:Tr=product}) is irrelevant. Furthermore, as already noted in \cite{Rezende:2010}, the inverse of $T^r$ is given by
\begin{equation}\label{eq:Trinverse}
(T^r)^{-1} = I+\sum_{s=1}^t \alpha_{p_sj_s}U^{p_sj_s}=
(I+\alpha_{p_1j_1}U^{p_1j_1})\cdots(I+\alpha_{p_{t_r} j_{t_r} }U^{p_{t_r} j_{t_r}}).
\end{equation}
This follows from (\ref{eq:matrixTr}), (\ref{eq:matrixunitproperty}) and the fact that a column cannot simultaneously contain a primary pivot and a change-of-basis pivot. In keeping with its counterpart, the pre-multiplication of $\Delta^r T^r$ by $(T^r)^{-1}$ consists of $t_r$ elementary row operations. Amongst the ones available, the only row operation on row $i$ considered herein is of the type ``add to row $i$ a multiple of another row''. This discussion implies the following lemma.

\begin{lem}\label{lema:columnrowchanges} Suppose the Incremental Sweeping Algorithm over $\F$ is applied to the connection matrix $\Delta\in \F^{m\times m}$, with column/row partition $J_0$, \ldots, $J_b$. Then the update from $\Delta^r$ to $\Delta^{r+1}$ may be accomplished blockwise, according to the individual updates
\begin{equation}\label{eq:blockupdateIncr}
\Delta^r_{J_{k-1}J_k} = (T^{r-1}_{J_{k-1}J_{k-1}})^{-1}\Delta^{r-1}_{J_{k-1}J_k}T^{r-1}_{J_kJ_k},\quad\mbox{for }k=1,\ldots,b.\end{equation}
Consequently, only columns containing change-of-basis pivots are subjected to elementary column operations, and only rows with the same index as the column of the primary pivots used for canceling out the change-of-basis pivots suffer elementary row operations.\end{lem}

\dem Items (i) and (iv) of Proposition \ref{prop:CorrectnessAlgorithmOverF} imply that the connection matrix update step may be restricted to the entries in $\cup_{k=1}^b J_{k-1}\times J_k$, since entries outside these positions are zero. The update for the individual blocks in (\ref{eq:blockupdateIncre}) is obtained using Lemma~\ref{lema:tecnico}. Column operations are due to the post-multiplication, and by construction of $T^{r-1}_{J_kJ_k}$, affect only columns of $J_k$ that contain change-of-basis pivots. Using (\ref{eq:Trinverse}), the pre-multiplication by $(T^{r-1}_{J_{k-1}J_{k-1}})^{-1}$ affects only rows with same index as the columns that contain the primary pivots used for canceling out the change-of-basis pivots.\qed

Let $\Delta\in \F^{m\times m}$ be a connection matrix, with column/row partition $J_0$, \ldots, $J_b$, to which the Incremental Sweeping Algorithm over $\F$ is applied. Suppose there is a change in entry in position $(i,j) \in J_{k-1}\times J_k$ from $\Delta^r$ to $\Delta^{r+1}$. Lemma~\ref{lema:columnrowchanges} implies this may be due to an elementary column operation on column $j$, due to a change-of-basis mark on entry $(j-r,j)$ on this column, and/or to an elementary row operation, due to a primary pivot in column $i\in J_{k-1}$ that is being used to cancel out a change-of-basis pivot in position $(q-r,q)$ in some column $q \in J_{k-1}$, where $q>i$. Let $\mathcal{J}_k$ be the set of columns in $J_k$ that contain primary pivot entries in $\Delta^m$, for $k=1,\ldots b$. Let $\overline{J}_k=J_k\backslash \mathcal{J}_k$, for all $k$. Corollary~\ref{cor:zeropattern4} and Proposition~\ref{prop:complF} imply that rows of $\Delta^m$ in $\mathcal{J}_k$ must be zero, for all $k$,  and thus cannot contain primary pivots. Furthermore, since change-of-basis pivots must be located to the right of and on the same row as primary pivots, there can be no change-of-basis pivots on rows in $\mathcal{J}_k$, for all $k$. Putting together these facts leads us to the last ingredient needed to reformulate the Incremental Sweeping Algorithm over $\F$.

\begin{lem}\label{lema:independence}
Let $\Delta \in \F^{m\times m}$ be a connection matrix with column/row partition $J_0$, \ldots, $J_b$. Consider the sweeping of the $r$-th diagonal of $\Delta^r$ in the Incremental Sweeping Algorithm over $\F$. The markings on entries in columns belonging to $J_k$ and the construction of \/$T^r_{J_k J_k}$ are completely determined by the values of the change-of-basis and primary pivots in $\Delta^r_{\overline{J}_{k-1}J_k}$.
\end{lem}

\dem Consider a fixed generic entry $\Delta^r_{j-r,j}$, where $j\in J_k$, eligible for receiving a mark. Since zero entries are not marked, we may suppose $\Delta^r_{j-r,j}\neq 0$ and, using Proposition~\ref{prop:CorrectnessAlgorithmOverF} (i), conclude that $j-r\in J_{k-1}$. Furthermore, rows in $\mathcal{J}_{k-1}$ cannot contain primary or change-of-basis pivots, so we may assume $j-r\in \overline{J}_{k-1}$. 

Nonzero entry $\Delta^r_{j-r,j}$ will not be marked if there exists a primary pivot below it, say in position $(i,j)$, $i>j-r$, which would imply that $i\in \overline{J}_{k-1}$, which means this primary pivot is in $\Delta^r_{\overline{J}_{k-1}J_k}$.

Now suppose column $j$ does not contain a primary pivot. If $\Delta^r_{j-r\textbf{\large .}}$ contains a primary pivot, say in position $(j-r,p)$, then, Proposition~\ref{prop:CorrectnessAlgorithmOverF} (iii), $\Delta^r_{j-r,p}\neq0$, and, by item (i), $p\in J_k$, so this primary pivot is located in $\Delta^r_{\overline{J}_{k-1}J_k}$. In this case, $\Delta^r_{j-r,j}$ will be marked as a change-of-basis pivot. Otherwise, it will be marked as a primary pivot. Notice that, in this case $j-r\in \overline{J}_{k-1}$, since rows in $\mathcal{J}_{k-1}$ cannot contain primary pivots.

Consider the construction of $T^r_{J_k J_k}$. The diagonal is always initialized with ones. Column $T^r_{\textbf{\large .}j}$ will have another nonzero entry only if $\Delta^r_{j-r,j}$ is a change-of-basis pivot. Suppose this is the case. Then, as we've established above, for some $p\in J_k$, there is a primary pivot $\Delta^r_{j-r,p}$, with $j-r\in \overline{J}_{k-1}$ and $p\in J_k$. This will give rise to the nonzero entry $-\Delta^r_{j-r,j}/\Delta^r_{j-r,p}$ in position $(p,j)$ of $T^r$.\qed

Lemma~\ref{lema:independence} implies the sweeping could be done sequentially blockwise as follows. Let $\Delta\in\F^{m\times m}$ be a connection matrix with column/row partition $J_0$, \ldots, $J_b$. During the blockwise sequential procedure, all connection matrices have the same column/row partition as $\Delta$. At the first step, one would apply the sweeping algorithm to a connection matrix that coincides with $\Delta$ in positions $J_0\times J_1$, and is zero otherwise. The transition submatrices $T^{1,r}_{J_1\,J_1}$, for $r=1$, \ldots, $m$, produced during the execution of the algorithm are recorded and will be used in the next step. At step $k$, for $k=1$, \ldots, $b$, we consider a connection matrix that coincides with $\Delta$ in positions $J_{k-1}\times J_k$, and is zero otherwise. The diagonals of this matrix are swept as in the Incremental Sweeping Algorithm over $\F$, and the same rules are used for marking up the entries, as well as building the transition matrix, which will contain nonzero off-diagonal entries only in positions $J_k\times J_k$, by virtue of the construction of the connection matrix used as input. The remaining principal submatrices $T^{k,r}_{J_s\,J_s}$, for $s\neq k$, are identity matrices, for all $r$. In the update step, only entries in positions $J_{k-1}\times J_k$ need be updated. The block update formula (\ref{eq:blockupdateIncr}) is employed, with the difference that the submatrices $T^{k-1,r}_{J_{k-1}\, J_{k-1}}$ produced in the previous step are used,  instead of the ones produced during the markup in the current step. Lemma~\ref{lema:independence} implies that submatrix $T^{k,r}_{J_k\,J_k}$ produced in step $k$ of this sequential procedure coincide with submatrix $T^r_{j_k\,J_k}$ constructed in the application of the Incremental Sweeping Algorithm over $\F$ to $\Delta$. Therefore the evolution of entries in positions $J_{k-1}\times J_k$ in step $k$ coincides with the evolution of entries of $\Delta$ in the same positions, during the application of the algorithm to $\Delta$. The sequence $\Delta^1$, \ldots, $\Delta^m$ can be obtained summing the sequences produced by sweeping the 1-block matrices as described.

Nevertheless, the procedure described above is still not simple enough for our purposes. Our aim is to forgo the pre-multiplication operation, and that of course comes at a price. Proposition~\ref{prop:complF} and Lemma~\ref{lema:columnrowchanges} imply that the only rows that change in the pre-multiplication are the ones that are equal to zero at the end of the sweeping. Furthermore, Lemma~\ref{lema:independence} implies that the entries in these rows do not contribute to the blocks $T^{k,r}_{J_kJ_k}$. Therefore, the sequential sweeping can be even further simplified, if we give up keeping track of the evolution of the rows in $\cup_{k=1}^b\mathcal{J}_k$. Instead of storing the blocks $T^{k-1,r}_{J_{k-1}\, J_{k-1}}$, one need only record the primary pivot columns of the last connection matrix in step $k-1$ and zero the rows with same index before applying the sweeping in the next step. The algorithm is stated below. Notice that in this sequential version, there are no elementary row operations, only elementary column operations in each execution of the Incremental Sweeping Algorithm over $\F$, since all blocks of connection matrix $\Delta(k)$ are zero, except for block $\Delta(k)_{J_{k-1},J_k}$, for $k=1,\ldots, b$.\\

\textbf{Block Sequential Sweeping Algorithm over $\F$}\\[-20pt]
\begin{description}
\item[\textbf{Input:}] nilpotent $m\times m$ upper triangular matrix $\Delta$ with column/row partition $J_0$, $J_1$, \ldots, $J_b$.

\item[\textbf{Initialization Step:}] $\mathcal{J}_0=\emptyset$

\item[\textbf{Iterative Step:}] \mbox{}\\
\textbf{For $k=1,\ldots,b$ do}\mbox{}\\
$\begin{array}{@{}l}
\left[\begin{tabular}{l}
Let $\Delta(k)$ be the matrix obtained from $\Delta$ by zeroing rows in $\mathcal{J}_{k-1}$ and entries\\ \rule{1cm}{0pt}outside positions in $J_{k-1}\times J_k$\\
Apply the Incremental Sweeping Algorithm over $\F$ to $\Delta(k)$, with column/row\\
\rule{1cm}{0pt}partition $J_0$,  \ldots, $J_b$, obtaining $\Delta(k)^m$\\
$\mathcal{J}_k =$ indices of columns of $\Delta(k)^m$ containing primary pivots
\end{tabular}
\right.\end{array}$
\item[\textbf{Output:}] $(\Delta(k)^0, \ldots, \Delta(k)^m)$ and $(T(k)^0,\ldots, T(k)^{m-1})$, for $k=1,\ldots, b$
\end{description}

\begin{thm}[Uncoupling]\label{thm:uncoupling} Let $\Delta\in\F^{m\times m}$ be a connection matrix with row/column partition $J_0$, \ldots, $J_b$. Let $\Delta^m$ be the matrix produced by the Incremental Sweeping Algorithm over $\F$ applied to $\Delta$, and let $\Delta(k)^m$, for $k=1,\ldots, b$, be the matrices obtained in the Block Sequential Sweeping Algorithm over $\F$ applied to $\Delta$. Then $\Delta^m_{J_{k-1}J_k} = \Delta(k)^m_{J_{k-1}J_k}$, for $k=1,\ldots, b$ and the collection of change-of-basis and primary pivots encountered in the application of the Incremental Sweeping Algorithm over $\F$ to $\Delta(k)$, for $k=1, \ldots, b$, coincides with the change-of-basis and primary pivots found when it is applied to $\Delta$.
\end{thm}

\dem Let  $T(k)^0$, $T(k)^1$, \ldots, $T(k)^m$ be the transition matrices constructed when applying the Incremental Sweeping Algorithm over $\F$ to $\Delta(k)$. Given that $\Delta(k)$ has at most one nonzero block, $\Delta(k)_{J_{k-1}J_k}$, Proposition~\ref{prop:CorrectnessAlgorithmOverF}, Lemma~\ref{lema:tecnico} and induction imply that  $T(k)^r_{J_iJ_i}$ is an identity matrix, for all $i\neq k$. Now, using Lemma~\ref{lema:columnrowchanges}, the update of $\Delta(k)$ is reduced to the following update
\[\Delta(k)^r_{J_{k-1}J_k} = (T(k)^{r-1}_{J_{k-1}J_{k-1}} )^{-1}\Delta(k)^{r-1}_{J_{k-1}J_k} T(k)^{r-1}_{J_kJ_k} = \Delta(k)^{r-1}_{J_{k-1}J_k} T(k)^{r-1}_{J_kJ_k},\]
for all $r\geq 1$. Since the update of $\Delta(k)$ involves only the post-multiplication step, only elementary column operations are performed.

Lemma~\ref{lema:independence} implies that $\Delta(1)^m_{J_0 J_1}=\Delta^m_{J_0 J_1}$ and the change-of-basis and primary pivots marked during the application of the algorithm to $\Delta(1)$ coincide with the ones marked in columns in $J_1$ when the algorithm is applied to $\Delta$. Assume by induction that $\Delta(k-1)^m_{J_{k-2} J_{k-1}}=\Delta^m_{J_{k-2} J_{k-1}}$ and change-of-basis and primary pivots of $\Delta(k-1)$ agree with the ones in columns in $J_{k-1}$ marked when $\Delta$ is swept.  Also by induction, $\mathcal{J}_{k-1} = N_{k-1}$. Proposition~\ref{prop:complF} implies rows of $\Delta^m$ in $\mathcal{J}_{k-1}$ are zero. By construction, rows of $\Delta(k)$ in $\mathcal{J}_{k-1}$ are also zero, and, since $\Delta(k)$ suffers only elementary column operations during the application of the Incremental Sweeping over $\F$, these rows are not changed. So $\Delta^m_{\mathcal{J}_{k-1}J_k} = \Delta(k)^m_{\mathcal{J}_{k-1}J_k}$.

By Lemma~\ref{lema:columnrowchanges}, entries in the rows of $\Delta$ in $\overline{J}_{k-1}$ are not subjected to elementary row operations during the application of the Incremental Sweeping Algorithm over $\F$ thereto. Furthermore, all change-of-basis and primary pivots in columns in $J_k$ occur in positions in $\overline{J}_{k-1}\times J_k$. Hence changes to entries in these rows are only due to elementary column operations, as also happens when the Incremental Sweeping is applied to $\Delta(k)$. Since $\Delta_{\overline{J}_{k-1}J_k}=\Delta(k)_{\overline{J}_{k-1}J_k}$ and the elementary column operations on columns in $J_k$ are solely dependent on the entries in this submatrix, it follows that the change-of-basis and primary pivots marked in columns in $J_k$ during execution of the Incremental Sweeping Algorithm over $\F$ to both $\Delta$ and $\Delta(k)$ coincide, and so $\Delta^m_{\overline{J}_{k-1}J_k}=\Delta(k)^m_{\overline{J}_{k-1}J_k}$. By induction, the statements are true for $k=1,\ldots, b$.\qed 

\comment{The Uncoupling Theorem~\ref{thm:uncoupling} implies that we may restrict our attention to connection matrices containing at most one nonzero block when studying the Incremental Sweeping Algorithm over $\F$, if we accept to miss the evolution of rows in $\cup_{k=1}^b \mathcal{J}_k$, which we know will end up zero and will not contain neither primary nor change-of-basis pivots. To ease the discussion that follows, we henceforth call this special case the \emph{1-Block Incremental Sweeping Algorithm over $\F$}.}

{The Uncoupling Theorem~\ref{thm:uncoupling} will lead to a significant simplification of the task of studying the effects of the application of the Incremental Sweeping Algorithm over $\F$ to a connection matrix with row/column partition $J_0$, \ldots, $J_b$, as long as one is content to forgo the observation of the evolution of rows in $\cup_{k=1}^b \mathcal{J}_k$ , which are known to end up zero and will not contain neither primary nor change-of-basis pivots. To ease the discussion that follows, we henceforth call the special case of the Incremental Sweeping Algorithm over $\F$, where the input is restricted to connection matrices with row/column partition consisting of only two subsets, the \emph{1-Block Incremental Sweeping Algorithm over $\F$}.}

The analysis of the 1-Block Incremental Sweeping Algorithm over $\F$ is much easier than that of its more general counterpart. If the columns of the connection matrix $\Delta\in \F^m$ are partitioned into two subsets $J_0$ and $J_1$, only rows in $J_1$ are altered in the pre-multiplication by $(T^{r-1})^{-1}$, in the matrix update step. But $\Delta_{J_1}=0$ and this zero pattern is invariant under elementary column operations. So the  post-multiplication part of the update doesn't change the nullity of rows in $J_1$ and the elementary row operations performed during the pre-multiplication part of the update step involve only the zero rows in $J_1$. This implies $\Delta^r_{J_1}=0$ for all $r$, and we may eliminate the pre-multiplication part of the update step, so that only elementary column operations need be performed. 

During the execution of the 1-Block Incremental Sweeping Algorithm over $\F$, columns of $\Delta^r$ may be classified as \emph{active} (resp., \emph{passive}), if they contain (resp., do not contain) a primary pivot mark. The active columns effect change upon the passive columns. The passive columns suffer changes caused by active columns. At the beginning of the algorithm all columns are passive and before changes are allowed to happen, at least one column must become active.  Once a column reaches the active state, it doesn't leave it, since primary pivot marks are permanent. Passive columns undergo a (possibly empty) sequence of elementary column operations and either reach an active state or become zero, since Corollary~\ref{cor:zeropattern4} implies columns without primary pivots must be zero. If a column reaches an active state, it does so when the lowest nonzero entry in the column is marked as a primary pivot. The order of sweeping implies that the change-of-basis pivots that occur in a fixed column, say $j$, are marked in an upward fashion. If the entry in position $(j-r,j)$ is marked as a change-of-basis pivot, and the entry in position $(j-r,p)$ contains the primary pivot to its left, columns $p$ and $j$ exhibit a sequence of trailing of zeros from row $j-r+1$ to the last row. The elementary column operation that eliminates this change-of-basis pivot changes only the entries in rows $1$ through $j-r$, the actual operation being determined by the values of the two pivots. By construction, each operation increases the number of trailing zeros by at least one. 

In the 1-Block Incremental Sweeping Algorithm over $\F$, the passive columns dictate the cancellations, which are done only once a change-of-basis pivot is marked. We propose a reengineered version therefor, in which this role is transferred to the active columns. Once an primary pivot is identified, all cancellations it is responsible for in the 1-block Incremental Sweeping Algorithm over $\F$ are performed. To arrive at the same final matrix as in the original algorithm, the primary pivots must also be identified in a upward order, from the bottom up. The second and last important aspect for the identification of primary pivots, is the left-to-right order of the sweeping. The algorithm below incorporates both of these. In the algorithm we adopt the usual convention that if $S$ is an empty set and $A$ is a real matrix, then $A_{S\textbf{\large .}}=0$. Since we are considering connection matrices with at most one nonzero block, we may assume, without loss of generality, that the row/column partition has two subsets.\\

\textbf{Revised 1-Block Incremental Sweeping Algorithm over $\F$}\\[-20pt]
\begin{description}
\item[\textbf{Input:}] nilpotent $m\times m$ upper triangular matrix $\Delta$ with column/row partition $J_0$, $J_1$.

\item[\textbf{Initialization Step:}] $C^1=\{1,\ldots, m\}$, $\tilde{\Delta}^1=\Delta$, $t=1$.

\item[\textbf{Iterative Step:}] \mbox{}\\
\textbf{While $\tilde{\Delta}^t_{\textbf{\large .}C^t}\neq 0$ do}\mbox{}\\
$\begin{array}{@{}l}
\left[\begin{tabular}{l}
Let $i_t=\max\{i\mid \tilde{\Delta}^t_{iC^t}\neq0\}$\\ 
Let $j_t=\min\{j\in C^t\mid \tilde{\Delta}^t_{i_t j}\neq0\}$\\
Permanently mark $\Delta^t_{i_t j_t}$ as a primary pivot\\
$\left[\begin{tabular}{l}
\textbf{Update Matrix Construction}\\
\begin{tabular}{@{\hspace{.5cm}}l}
$\tilde{T}^t\leftarrow I - \displaystyle\sum_{\scriptsize\begin{array}{c}j\in C^t\\[-5pt] j>j_t\end{array}} \frac{\tilde{\Delta}^t_{i_t j}}{\tilde{\Delta}^t_{i_t j_t}} U^{j_t j}$\\
\end{tabular}\end{tabular}\right.$\\
\\[-15pt]
$\begin{array}{@{}l}
\left[\begin{tabular}{l}
\textbf{Simplified Matrix $\Delta$ update}\\
\rule{.5cm}{0pt}$\tilde{\Delta}^{t+1} = \tilde{\Delta}^t\tilde{T}^t$\\
\end{tabular}\right.\\
C^{t+1}\leftarrow C^t\backslash \{j_t\}\\
t\leftarrow t+1 \\
\end{array}$
\end{tabular}
\right.\end{array}$
\item[\textbf{Output:}] $(\tilde{\Delta}^0, \ldots)$ and $(\tilde{T}^0,\ldots)$
\end{description}

The attentive reader may realize that the Revised 1-block Incremental Sweeping Algorithm over $\F$ is an implementation of the transposed version of the Gaussian method for obtaining a row echelon matrix, minus the row swaps (which would amount to column swaps in the revised algorithm).

\begin{prop}\label{prop:RevisedPattern} Let $(i_1,j_1)$, \ldots, $(i_{t^*},j_{t^*})$ be the positions of primary pivots marked in the application of the Revised 1-block Incremental Sweeping Algorithm over $\F$ to the connection matrix $\Delta\in \F^{m\times m}$ with row/column partition $J_0$, $J_1$, in the order in which they were marked. Then the following are true:
\begin{itemize}
\item[(i)] once a column receives a primary pivot mark, it remains invariant until the end of the algorithm,
\item[(ii)] $j_t>i_t$, for $t=1,\ldots, t^*$,
\item[(iii)] $i_1>i_2>\cdots>i_{t^*}$,
\item[(iv)] $\tilde{\Delta}^{t+1}_{\{i_t,\dots,m\}C^{t+1}}=0$, for $t=1,\ldots,  t^*$,
\item[(v)] the number of consecutive zero entries at the bottom of each column never decreases.
\end{itemize}
\end{prop}

\dem Without loss of generality, suppose $\Delta\neq 0$, otherwise the proposition is vacuously true. Given the construction of $\tilde{T}^t$, only columns $j$ such that $j> j_t$ and $j\in C^t$ change from $\tilde{\Delta}^t$ to $\tilde{\Delta}^{t+1}$. Namely,
\begin{equation}\label{eq:update}\tilde{\Delta}^{t+1}_{\textbf{\large .} j} = \left\{\begin{array}{ll}
\tilde{\Delta}^t_{\textbf{\large .} j}  -  \frac{\tilde{\Delta}^t_{i_t j}}{\tilde{\Delta}^t_{i_t j_t}}\tilde{\Delta}^t_{\textbf{\large .} j_t}, \quad&\mbox{if }j\in C^t \mbox{ and }j> j_t,\\
\tilde{\Delta}^t_{\textbf{\large .} j} , & \mbox{otherwise},\end{array}\right.\end{equation}
so (i) is true.

The initial matrix is upper triangular and this property is preserved under the update step, since columns that changed receive the sum of the old column and a multiple of a column to its left. Hence $j_t>i_t$ for all $t$, since $\tilde{\Delta}^t_{i_t j_t}\neq 0$. Thus (ii) is also true.

By construction, the lowest row of $\tilde{\Delta}^1$ with nonzero entries is $i_1$. Formula (\ref{eq:update}) gives
\begin{equation}\label{eq:pivo1}\tilde{\Delta}_{i_1 j}^2 =\left\{\begin{array}{ll}
\tilde{\Delta}^1_{i_1 j} - \frac{\tilde{\Delta}^1_{i_1 j}}{\tilde{\Delta}^1_{i_1 j_1}}\tilde{\Delta}^1_{i_1 j_1}=0,\quad&\mbox{if }j> j_1,\\
\tilde{\Delta}^1_{i_1 j_1}, &\mbox{if }j\leq j_1.\end{array}\right.\end{equation}
Using the fact that there are no nonzero entries to the left of $\tilde{\Delta}^1_{i_1j_1}$, the only nonzero entry on or below row $i_1$ in $\tilde{\Delta}^2$ is the primary pivot in position $(i_1,j_1)$, which proves (iv) for $t=1$. Since $j_1\notin C^2$, the lowest nonzero row of $\tilde{\Delta}^2_{\textbf{\large .}C^2}$ must be higher than $i_1$, so $i_2>i_1$. 

Assume by induction that $i_1>\cdots>i_{t-1}$ and $\tilde{\Delta}^t_{\{i_{t-1},\ldots,m\} C^t}=0$. The last fact implies $i_{t-1}>i_t$, proving (iii) by induction. The definition of $j_t$ implies $\tilde{\Delta}^t_{\{i_t+1,\ldots,m\} C^t}=0$. Let $j\notin \{j_1, \ldots, j_t\}$. These zero entries won't be altered by the elementary column operations embodied in the matrix update formula (\ref{eq:update}), hence $\tilde{\Delta}^{t+1}_{\{i_t+1,\ldots,m\} C^t}=0$. When $i=i_t$, (\ref{eq:update}) implies
\[\tilde{\Delta}_{i_t j}^{t+1} = \left\{\begin{array}{ll}\tilde{\Delta}_{i_t j}^t - \frac{\tilde{\Delta}^t_{i_t j}}{\tilde{\Delta}^t_{i_t j_t}}\tilde{\Delta}^t_{i_t j_t}=0,\quad&\mbox{if }j\in C^t \mbox{ and }j>j_t,\\
\tilde{\Delta}^t_{i_t j},&\mbox{if }j\in C^t \mbox{ and } j\leq j_t.\end{array}\right.\]
Using the  that $(i_t,j_t)$ is the position of the leftmost nonzero entry of $\tilde{\Delta}^t_{i_t C^t}$, we conclude that $\tilde{\Delta}^{t+1}_{i_t j_t}$ is the only nonzero entry on row $i_t$ in the submatrix $\tilde{\Delta}^{t+1}_{\textbf{\large .}C^{t+1}}$. Therefore (iv) is true for $\tilde{\Delta}^{t+1}$. By induction, it follows that (iv) is true for all values of $t$.

Finally, to show (v), notice that a column $j$ suffers change when $j\in C^t \backslash \{j_t\}$, for some $t$. But if this is so, the entries in column $j$, in rows $i_t+1, \ldots, m$ must be zero and the update (\ref{eq:update}) will zero the entry in position $(i_t,j)$ while keeping the entries lower down at zero. Entries higher up may or not zero out, and may even change from zero to a nonzero value. But if a column changes, the number of trailing zeros in that column will always increase by at least one.\qed

\begin{cor}\label{cor:1blockend} Let $\tilde{\Delta}^{t^*+1}$ be the last matrix obtained by the application of the Revised 1-block Incremental Sweeping Algorithm to the connection matrix $\Delta \in \F^{m\times m}$ with row/column partition $J_0$, $J_1$. Then the primary pivot entries are nonzero and each nonzero entry of $\tilde{\Delta}^{t^*+1}$ lies above a primary pivot. \end{cor}

\dem Let $(i_1,j_1)$, \ldots, $(i_{t^*},j_{t^*})$ be the positions of primary pivots. A simple induction shows that $C^{t^*+1}$ is the set of indices of columns of $\tilde{\Delta}^{t^*+1}$ without primary pivots. The stopping criterium implies $\tilde{\Delta}^{t^*+1}_{\textbf{\large .} C^{t^*+1}}=0$. Finally, primary pivots, when marked, are, by the rules of the algorithm, the lowest nonzero entry of the column, and, by Proposition~\ref{prop:RevisedPattern} (i), columns do not change after receiving a primary pivot mark.\qed

The next proposition establishes the equality between the final matrices produced by the Revised 1-block Incremental Sweeping over $\F$ and the 1-block Incremental Sweeping Algorithm over $\F$. There is of course no sense in looking for equality between other matrices in the sequence produced by the algorithm, since the order of cancellation is in all likelihood quite different in the two algorithms.

\begin{prop}\label{prop:RevisedCorrect} Let $\Delta\in \F^{m\times m}$ be a connection matrix with column/row partition $J_0$, $J_1$. Let $\tilde{\Delta}^{t^*+1}$ and $\Delta^m$ be the matrices obtained by applying the Revised 1-block Incremental Sweeping Algorithm over $\F$ and the 1-Block Incremental Sweeping Algorithm over $\F$ to $\Delta$, respectively. Then $\tilde{\Delta}^{t^*+1}=\Delta^m$ and their primary pivots coincide.
\end{prop}

\dem 
If $\tilde{\Delta}^{t^*+1}=\Delta^m$, then their primary pivots coincide in position and value, by Corollaries~\ref{cor:zeropattern4} and~\ref{cor:1blockend}.

In both algorithms columns may suffer elementary column operations until they either reach zero or receive a primary pivot mark. Additionally, a column may suffer an elementary column operation only from another column with a primary pivot mark on its left and successive operations on a column may only increase the range of trailing zeros in that column, that is, the set of successive rows, ending in row $m$, containing zero entries. The order in which the primary pivots are identified probably differs between algorithms, but the important thing is that when an entry is eligible for receiving a primary pivot in either algorithm, the column it belongs to will have been subjected to the same changes in both of them. In order for these changes to be the same in both algorithms, the active columns acting on them must be the same and the changes they provoke in a fixed column must occur in the correct order, from bottom up.

The proof is by induction on the number of nonzero columns of $\Delta$. If there is but one nonzero column, the nonzero entry at the bottom of this column will be marked by both algorithms as the unique primary pivot of $\Delta$, so $t^*=1$, there will be no elementary column operations and in fact $\tilde{\Delta}^{t^*+1} = \Delta^m = \Delta$. 

Admit by induction that the last matrix of both algorithms coincide, when the number of nonzero columns is smaller than $k$. Suppose $\Delta$ has $k$ nonzero columns. Let $(i_1,j_1)$ be the first position to receive a primary pivot mark in the Revised 1-block Incremental Sweeping Algorithm over $\F$. Then the entry $\Delta^{j_1-i_1}_{i_1 j_1}$ must also receive a primary pivot mark in the Incremental Sweeping Algorithm over $\F$. To see that, note that, by definition of $(i_1,j_1)$, the entries in $\Delta_{\{i_1,\ldots,m\} \{1,\ldots, j_1-1\}}$ and $\Delta_{\{i_1+1,\ldots, m\}\{1,\ldots,m\}}$ are zero, and the number of trailing zeros can only increase. Consequently no entries in these submatrices may have been marked as a primary pivot before the sweeping of the $j_1-i_1$ diagonal, so $\Delta^{j_1-i_1}_{i_1j_1}$ is nonzero, with no primary pivot marks on its left or below it. So, in this case, $\Delta^{t^*+1}_{\textbf{\large .} j_1} = \Delta^m_{\textbf{\large .} j_1} = \Delta_{\textbf{\large .}j_1}$. Furthermore, notice that, entries in positions $(i_1,j_1+1)$, \ldots $(i_1,m)$ will be swept after this and be marked, if nonzero, as change-of-basis entries in the Incremental Sweeping Algorithm over $\F$, since rows $i_1+1$, \ldots, $m$ of $\Delta$ are zero and the entry in position $(i_1,j_1)$ has  a primary pivot. The changes possibly effected on columns due to these markings in the Incremental Sweeping Algorithm over $\F$ on the corresponding iterations are precisely the changes done in the first iteration of the Revised 1-block Incremental Sweeping Algorithm over $\F$.

Let $\Delta' \in \F^{m\times m}$ be defined as follows:
\[\Delta'_{\textbf{\large .}j} = \left\{\begin{array}{ll} 
\tilde{\Delta}^2_{\textbf{\large .}j},\quad&\mbox{if } j\neq j_1,\\
0, &\mbox{otherwise.}\end{array}\right.\]
The matrix $\Delta'$ agrees with the matrix obtained from $\Delta$ after the first iteration of the Revised 1-block Incremental Sweeping Algorithm over $\F$, except for column $j_1$, which is zero. So $\Delta'$ encompasses the changes to columns due to change-of-basis entries in row $i_1$, and has the $j_1$-th column equal to zero. Thus if we apply the Incremental Sweeping Algorithm over $\F$ to $\Delta'$, the matrix ${\Delta'}^{m}$ obtained agrees with $\Delta^m$, except for column $j_1$, which would remain zero throughout the algorithm. Analogously, $\Delta'$ encompasses changes made to $\Delta$ in the first iteration of the Revised 1-block Incremental Sweeping Algorithm over $\F$, but differs from it in column $j_1$, which has been zeroed. Thus the application of the Revised 1-block Incremental Sweeping Algorithm over $\F$ to $\Delta'$ produces a matrix $\tilde{\Delta'}^{t^*}$ whose columns coincide with the corresponding ones from $\tilde{\Delta}^{t^*+1}$, except for the $j_1$-th column. Since $\Delta'$ has at most $k-1$ nonzero columns, induction implies that ${\Delta'}^{m} = \tilde{\Delta'}^{t^*+1}$. This implies 
\[\Delta^m_{\textbf{\large .}j}= {\Delta'}^{m}_{\textbf{\large .}j}=\tilde{\Delta'}^{t^*}_{\textbf{\large .}j}=\tilde{\Delta}^{t^*+1}_{\textbf{\large .}j},\quad \mbox{for }j\neq j_1.\]
But since we already have that $\Delta^m_{\textbf{\large .}j_1}=\tilde{\Delta}^{t^*+1}_{\textbf{\large .}j_1}$, we conclude that $\Delta^m=\tilde{\Delta}^{t^*+1}$.\qed

\section{Sweeping algorithm for {TU} connection matrices}\label{sec:TUsweeping}

{A matrix is totally unimodular  if all its square submatrices have determinant $0$, $1$ or $-1$, see, for instance, \cite{Nemhauser:1988}. This property is invariant under transposition, multiplying a row or column by $0,\pm 1$, adding or removing zero rows/columns or unit rows/columns. In this section  we explore the additional facts that may be established when the input to the  Spectral Sequence Sweeping Algorithm is restricted to totally unimodular connection matrices. An important class of connection matrices satisfying this property is that associated to Morse flows on surfaces. }

The set $\{0,\pm 1\}$ is not closed under addition, so the appropriate algorithm to apply to totally unimodular connection matrices is the Sweeping Algorithm over $\Z$. Nevertheless, using the Block Sequential Sweeping Algorithm over $\F$, the Revised 1-Block Incremental Sweeping Algorithm over $\F$, Theorem~\ref{thm:uncoupling} and the total unimodularity property, we will show that the sequence of matrices and bases obtained when applying the Incremental Sweeping Algorithm over $\F$ to a TU connection matrix is compatible with the corresponding ones produced by the application of Sweeping Algorithm over $\Z$ to $\Delta$, in the following sense. The change from $\Delta^r$ to $\Delta^{r+1}$ results from replacing basis element $\sigma^r_j$, of each column $j$ containing a change-of-basis entry, with an integral linear combination of elements of $h$ with same chain index and associated with columns of index less than or equal $j$. The coefficients of this linear combination are calculated so as to zero out the entry in the change-of-basis position in $\Delta^{r+1}$, while maintaining the pattern of trailing zeros below it. Furthermore, amongst the integral linear combinations that accomplish this, one must choose one with the smallest possible positive leading coefficient. Even with this condition, there may be more than one optimal integral linear combination. The bases constructed by the Incremental Algorithm over $\F$ trivially satisfy the conditions pertaining the zero patterns. We show that, for every $r$ and $j$, the leading coefficient in the integral linear combination that is $\sigma^r_j$, produced by the application of the Incremental Algorithm over $\F$ to $\Delta$, is $1$, so it satisfies the optimality criterium, and the choice is thus compatible with the rules of the Sweeping Algorithm over $\Z$. 

In the context of linear programming, the pivot operation on a nonzero entry $a_{pq}$, called the \emph{pivot}, of an $m\times n$ matrix $A$ is defined as the following set of elementary row operations on $A$: 1) add to row $i$, for $i\neq p$, $-a_{iq}/a_{pq}$ times row $p$, 2) divide row $p$ by $a_{pq}$. These operations may be expressed as the pre-multiplication of $A$ by the $m\times m$ matrix $B$ that coincides with the identity matrix except for its $p$-th column, which is given by 
\[b_{ip} = \left\{\begin{array}{ll} -\dfrac{a_{iq}}{a_{pq}},\quad&\mbox{if }i\neq p,\\
\dfrac{1}{a_{pq}},&\mbox{otherwise,}\end{array}\right.\]
see, for instance, \cite{Bazaraa:2010}. The pivot operation assigns 1 to entry in position $(p,q)$ and cancels entries on top and below position $(p,q)$, so that column $q$ of $BA$ is equal to $e_p$, the $p$-th canonical basis vector of $\R^m$. It follows from Lemma 9.2.2 of \cite[p. 190]{Truemper:1998}, that, if $A$ is totally unimodular, then $BA$ also is. Now consider a variant of linear programming pivoting, where the $m\times m$ matrix $\widetilde{B}$ coincides with $B$ except for entry in position $(p,p)$, which is equal to $1$. This means $(\widetilde{B}A)_{i\textbf{\large .}} = (BA)_{i\textbf{\large .}}$, for $i\neq p$, and $(\widetilde{B}A)_{p\textbf{\large .}} = a_{pq} (BA)_{p\textbf{\large .}}$. If $A$ is totally unimodular, then $a_{pq}=\pm 1$, so $\widetilde{B}A$ may be obtained from $BA$ by at most a change in sign of row $p$ of $BA$. But then, $\widetilde{B}A$ is also totally unimodular.

\begin{prop}\label{prop:unitarypivots}
Let $\Delta\in\{0,\pm1\}^{m\times m}$ be a totally unimodular connection matrix with column/row partition $J_0$, $J_1$. If we apply the Incremental Sweeping Algorithm over $\F$ to $\Delta$, then all primary pivots value either $1$ or $-1$.\end{prop}

\dem First we analyze the application of the Revised 1-Block Sweeping Algorithm over $\F$ to $\Delta$. Let $(i_1,j_1)$, \ldots, $(i_{t^*},j_{t^*})$ be the positions of the primary pivots marked. 

We claim that $\tilde{\Delta}^t_{\textbf{\large .}C^t}$ is totally unimodular, for $t=1,\ldots, t^*+1$. This is trivially true for $t=1$, by hypothesis. Assume it is true for $t$. Since
\[(\tilde{\Delta}^t T^t)_{\textbf{\large .}j} = \left\{\begin{array}{ll}
\tilde{\Delta}^t_{\textbf{\large .}j},&\mbox{if }j\notin C^t,\\
\tilde{\Delta}^t_{\textbf{\large .}j}, &\mbox{if }j\in C^t \mbox{ and }j\leq j_t,\\
\tilde{\Delta}^t_{\textbf{\large .}j}-\dfrac{\tilde{\Delta}^t_{i_t j}}{\tilde{\Delta}^t_{i_t j_t}} \tilde{\Delta}^t_{\textbf{\large .}j_t},&\mbox{if }j\in C^t\mbox{ and } j> j_t,\end{array}\right.\]
the marking of an entry in position $(i_t,j_t)$ of $\tilde{\Delta}^t_{\textbf{\large .} C^t}$ as a primary pivot implies the cancellation, in $\tilde{\Delta}^{t+1}$, of entries on row $i_t$, in columns in $C^t$ other than $j_t$ (although the update matrix construction provides the cancellation of entries to the right of column $j_t$, entries to its left are zero, since $\tilde{\Delta}^t_{i_tj_t}$ is the leftmost nonzero entry in row $i_t$ of $\tilde{\Delta}^t_{\textbf{\large .}C^t}$).  Thus $\tilde{\Delta}^{t+1}_{\textbf{\large.}C^t}=(\tilde{\Delta}^tT^t)_{\textbf{\large.}C^t} = \tilde{\Delta}^t_{\textbf{\large .}C^t} T^t_{C^t C^t}$. Notice that the cancellations are achieved by adding to column $j\in C^t$, $j\neq j_t$, the appropriate multiple of column $j_t$. But this is simply a transposed version of the variant of the linear programming pivoting described above, 
with $\tilde{\Delta}^t_{\textbf{\large .}C^t} = A^T$ and $T^t_{C^t C^t} = \widetilde{B}^T$. Therefore, if $\tilde{\Delta}^t_{\textbf{\large .}C^t}$ is totally unimodular, then $\tilde{\Delta}^{t+1}_{\textbf{\large.}C^t}$ is also totally unimodular. Since this property is, by definition, inherited by  submatrices, and $C^{t+1}\subset C^t$, we conclude $\tilde{\Delta}^{t+1}_{\textbf{\large .}C^{t+1}}$ is also totally unimodular. By induction, $\tilde{\Delta}^t_{\textbf{\large .}C^t}$ is totally unimodular for all $t$.

This implies $\tilde{\Delta}^t_{\textbf{\large .}C^t}\in \{0,\pm 1\}^{m\times |C^t|}$ for all $t$. Hence all primary pivots marked in the application of the Revised 1-Block Sweeping Algorithm over $\F$ to $\Delta$ value $1$ or $-1$ when marked, since they are, by choice, nonzero. Finally, they do not change once marked, see Proposition~\ref{prop:RevisedPattern} (i).

Proposition~\ref{prop:RevisedCorrect} implies $\Delta^{t^*+1}=\Delta^m$, the last matrix produced by the application of the Incremental Sweeping Algorithm over $\F$ to $\Delta$. Therefore, by Corollaries~\ref{cor:zeropattern4} and~\ref{cor:1blockend}, their primary pivots coincide in position and value.
\qed

\begin{figure}[htbp]\small
$\Delta^0=$\raisebox{-.17\textwidth}{\includegraphics[height=.4\textwidth,viewport=15 0 200 170,clip]{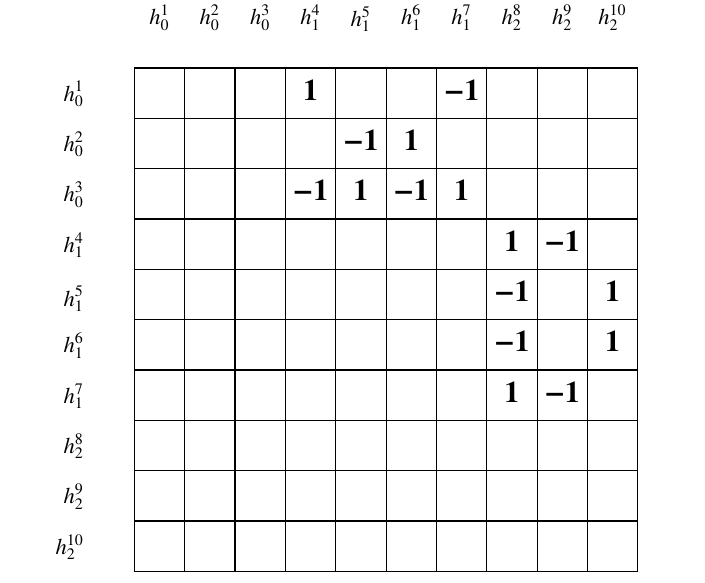}}%
\hfill$\Delta^3=$\kern-.2cm\raisebox{-.17\textwidth}{\includegraphics[height=.4\textwidth]{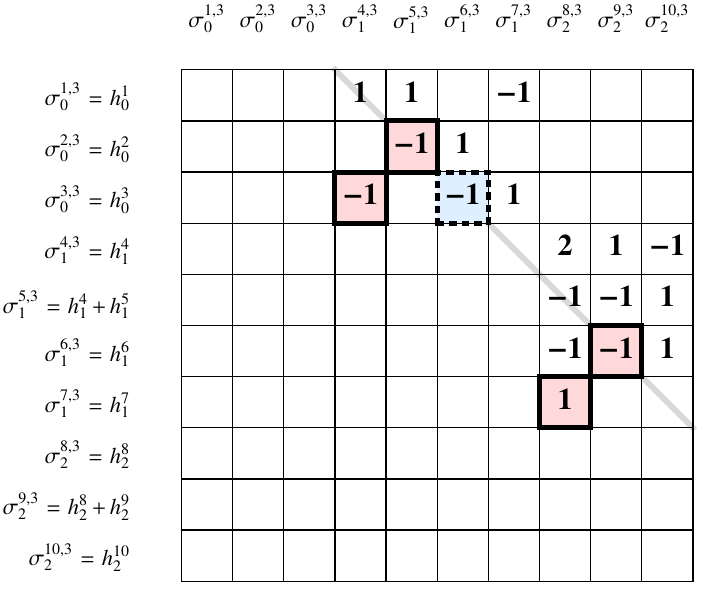}}%
\caption{$\Delta^3$ obtained by Incremental Sweeping Algorithm over $\F$ applied to $\Delta^0$.}\label{fig:surfaceexample}
\end{figure}

The elements are in place to establish that the primary pivots obtained when applying the Incremental Sweeping Algorithm over $\F$ to a TU connection matrix are unitary, see Theorem~\ref{thm:unitpivotTU}. This justifies its application, in place of the Sweeping Algorithm over $\Z$, to TU connection matrices, see Proposition~\ref{prop:leadingunitary} and Corollary~\ref{cor:finsteadofz}, even though the matrices produced during the sweeping may contain entries that would not have an inverse in $\Z$. An example is shown in Figure~\ref{fig:surfaceexample}. The highlighted entries in position $(3,4)$ and $(7,8)$ have been marked as primary pivots during the sweeping of the first diagonal. No entry in the second diagonal is a primary pivot. When the third diagonal is swept, in $\Delta^3$, entries in $(2,5)$ and $(6,9)$ are marked as primary pivots, entry in $(3,6)$ as a change-of-basis pivot. Notice that $\Delta^3_{4,8}=2$. Since column $4$ contains a primary pivot, row $4$ will be zero at the end of the algorithm.

\begin{thm}[Primary pivots for TU connection matrices]\label{thm:unitpivotTU} Let $\Delta$ be a TU connection matrix. Then the primary pivots obtained when applying the Incremental Sweeping Algorithm over $\F$ thereto value $\pm 1$.
\end{thm}

\dem Suppose $\Delta$ has column/row partition $J_0$, \ldots, $J_b$. The Uncoupling Theorem~\ref{thm:uncoupling} justifies the application of the Block Sequential Sweeping Algorithm over $\F$ to $\Delta$. If $\Delta$ is TU, then $\Delta_{J_{k-1}J_k}$ is TU, for $k=1$, \ldots, $b$,  and this property is maintained if one adds zero rows and/or columns to a matrix. Thus each $\Delta(k)$, for $k=1$, \ldots, $b$, in the application of the Block Sequential Algorithm over $\F$ to $\Delta$ is TU. Then the result follows from Proposition~\ref{prop:unitarypivots}.\qed

\begin{prop}\label{prop:leadingunitary}
Suppose that all the primary pivots obtained when applying the Incremental Sweeping Algorithm over $\F$ to an integral connection matrix $\Delta$ are $\pm1$. Let $\sigma^r$ be the basis associated with the $r$-th matrix in the sequence produced by the algorithm, $\Delta^r$. Then each basis element $\sigma^r_j$, for all $r$ and $j$, is an integral linear combination of the elements of $h$ associated with columns to the left of, and including, column $j$, whose leading coefficient is $1$.
\end{prop}

\dem First of all, from the definition of $T^r$ in the Incremental Sweeping Algorithm over $\F$ and the hypotheses that $\Delta$ is integral and the primary pivots value $\pm 1$, we obtain that $T^1$ and $(T^1)^{-1}$, given by (\ref{eq:Trinverse}), are also integral. Then $\Delta^2$ is integral and, applying induction, we may conclude that $\Delta^r$ is integral, for all $r$.

Initially $\sigma^1=h$, so the property is true for $r=1$. Suppose it is true for $\sigma^r$. If $\Delta^r_{\textbf{\large .}j}$ does not contain a change-of-basis, then $\sigma^{r+1}_j=\sigma^r_j$ and the result is true by the induction hypothesis. Suppose it does contain. It is sufficient to consider one change-of-basis pivot, in a generic fixed position, say $(j-r,j)$. Proposition~\ref{prop:unitarypivots} implies that the primary pivot to its left, say in position $(j-r,p)$,  is $\pm 1$. Then, by the rules of the algorithm, 
\begin{eqnarray*}
\sigma^{r+1}_j &=& \sigma^r_j - \frac{\Delta^r_{j-r,j}}{\Delta^r_{j-r,p}}  \sigma^r_p \\
&=&  \sigma^r_j \pm \Delta^r_{j-r,j}\,\sigma^r_p\\
&=& h^j_k + \sum_{j'<j}c^{r,j}_{j'} h^{j'}_k \pm \Delta^r_{j-r,j}\sum_{j'<p} c^{r,p}_{j'} h^{j'}_k,\end{eqnarray*}
where the induction hypothesis is used in the last equality. This proves the result for $r+1$, since $p<j$, $\Delta^r_{j-r,j}$ and, by induction, all coefficients $c^{r\textbf{\large .}}_{\textbf{\large .}}$ are integer. Using induction, the result if true for all $r$.\qed

\begin{cor}\label{cor:finsteadofz} Suppose that all the primary pivots obtained when applying the Incremental Sweeping Algorithm over $\F$ to an integral connection matrix $\Delta$ are $\pm1$. Then the sequence of bases associated with the matrices $\Delta^1$, $\Delta^2$, \ldots, is compatible with the application of the Sweeping Algorithm over $Z$ to $\Delta$. In particular, the sequence produced by the application of the Incremental Sweeping Algorithm over $\F$ to a TU connection matrix is compatible with the application of the Sweeping Algorithm over $\Z$ thereto.
\end{cor}

\dem The first result follows from Proposition~\ref{prop:leadingunitary}. The second is a special case, since a surface connection matrix is integral and, by Theorem~\ref{thm:unitpivotTU}, the application of the Incremental Sweeping Algorithm over $\F$ thereto produces unitary primary pivots.\qed

\section{Row cancellation algorithm}\label{sec:rowcancellation}

In the Revised 1-Block Incremental Sweeping Algorithm over $\F$ there are no change-of-basis pivots, each primary pivot ``takes the lead'' in the sense that all cancellations it would be responsible for in the Incremental Sweeping Algorithm over $\F$ are performed at once, as soon as it is marked. To achieve the same end matrix as the Incremental Sweeping Algorithm over $\F$, the sweeping order is altered and only one primary pivot is marked per iteration. 

This inspired the alteration of the Incremental Sweeping Algorithm over $\F$ to have all primary pivots have this commandeering role. In the new algorithm we keep the same order of sweeping along the diagonals and the criteria for marking an entry as a primary pivot. But the upper triangular unit-diagonal transition matrix is calculated so that in the next matrix all entries to the right of the primary pivots are zeroed by means of elementary column operations using exclusively the primary pivot columns. Although the computation of the transition matrix is focused only on this post-multiplication aspect, it will be shown that the pre-multiplication by its inverse does not interfere with the cancellations achieved in the post-multiplication. 

The new algorithm, called \emph{Row Cancellation Algorithm over $\F$}, is developed for generic connection matrices with entries in $\F$. Nevertheless, we will see that it can be applied to surface connection matrices, and, in this special case, it is possible to give a topological/dynamic interpretation to its workings. 

The next proposition shows that the transition matrix with the properties described and with the desired effect upon the connection matrix is unique and gives a constructive formula therefor.

\begin{prop}\label{prop:transitionmatrixSmale}
Let $A\in \F^{m\times m}$ and $P=\{j_1,\ldots, j_t\}$ be a subset of column indices, such that $j_1< \cdots < j_t$ and $SP=\{j_1-r,\ldots, j_k-r\}\subseteq 1..m$, for some positive integral $r$. Suppose $a_{j_s-r,j_s}\neq 0$, while entries in column $j_s$ below $a_{j_s-r,j_s}$ are zero, for $s=1$, \ldots, $t$. Then there exists a unique upper triangular unit-diagonal matrix $T\in\F^{m\times m}$ such that \begin{itemize}
\item[(i)] the entries of $T$, strictly above the diagonal and outside rows with indices in $P$, are zero;
\item[(ii)] $(AT)_{j_s-r,j_s+1..m}=0$, for $s=1$, \ldots, $t$.\end{itemize} 
This matrix may be expressed as the product
\[T^1\cdots T^t,\]
where
\[(T^s)_{i\mbox{\bf\large.}}=\left\{\begin{array}{ll} e_i,&\mbox{if }i\neq j_s,\\
e_{j_s} - \dfrac{1}{a_{j_s-r,j_s}} (0 \cdots 0~a_{j_s-r,j_s+1} \cdots a_{j_s-r,m}),\quad &\mbox{ow.}\end{array}\right.\] 
\end{prop}

\dem First we show that $T=T^1\cdots T^t$ satisfies the conditions and then show its unicity. By construction, $T^s$ is a unit-diagonal upper triangular matrix satisfying (i), for $s=1$, \ldots, $t$. The post-multiplication of a matrix, say $B$, by $T^s$ has the effect of adding to column $j$ of $B$ a multiple of column $j_s$ of $B$, for $j>j_s$. If $B$ is unit-diagonal upper triangular and satisfies (i), these properties will be preserved by this operation. Thus, by induction, $T^1 \cdots T^t$ is unit-diagonal, upper triangular and satisfies (i).

By construction, columns $1$ through $j_1$ of $A$ do not change with the post-multiplication by $T^1$. The entries strictly to the right of position $(j_1-r,j_1)$ are zeroed out, since
\[(AT^1)_{j_1-r,j_1+k} = a_{j_1-r,j_1} \cdot\left(-\frac{1}{a_{j_1-r,j_1}} a_{j_1-r,j_1+k}\right) + a_{j_1-r,j_1+k}=0,\]
for $k=1$, \ldots, $m-j_1$. Because entries of $A$ below $a_{j_1-r,j_1}$ are zero by hypothesis, rows of $A$ strictly below $j_1-r$ do not change with the post-multiplication by $T^1$. The fact that $j_1$ is the smallest element of $P$ implies that $(AT^1)_{j_s-r..m,j_s}=A_{j_s-r..m,j_s}$, for $s=1$, \ldots, $t$, and $(AT^1)_{j_1-r+1..m,j}=A_{j_1-r+1..m,j}$, for $j> j_1$.

Assume by induction that $(AT^1\cdots T^s)_{j-r..m,j} = A_{j-r..m,j}$, for $j=j_1$, \ldots, $j_t$, \linebreak $(AT^1\cdots T^s)_{j_s-r+1..m,j} = A_{j_s-r+1..m,j}$, for $j>j_s$, and that entries of $AT^1\cdots T^s$ strictly to the right of positions $(j_1-r,j_1)$, \ldots, $(j_s-r,j_s)$ are zero. The post-multiplication of $AT^1\cdots T^s$ by $T^{s+1}$ has the effect of adding to column $(AT^1 \cdots T^s)_{\mbox{\bf\large .}j}$ a multiple of column $(AT^1\cdots T^s)_{\mbox{\bf\large .}j_{s+1}}$, for $j>j_{s+1}$. Thus columns $1$ through $j_{s+1}$ of $(AT^1\cdots T^s)$ do not change with the post-multiplication by $T^{s+1}$ and, in particular,  $(AT^1\cdots T^{s+1})_{j-r..m,j} = A_{j-r..m,j}$, for $j=j_1$, \ldots, $j_{s+1}$. If $j>j_{s+1}$,
\begin{eqnarray*} 
(AT^1\cdots T^sT^{s+1})_{\mbox{\bf\large.}j} &=& (AT^1\cdots T^s)T^{s+1}_{\mbox{\bf\large .} j}\\
 & = & (AT^1\cdots T^s)_{\mbox{\bf\large.}j} - \frac{1}{a_{j_{s+1}-r,j_{s+1}}} a_{j_{s+1}-r,j} (AT^1\cdots T^s)_{\mbox{\bf\large.}j_{s+1}} .\end{eqnarray*}

Then, for $j>j_{s+1}$, the induction hypothesis implies
\[  (AT^1\cdots T^sT^{s+1})_{j_k-r,j} \!=\! \left\{\begin{array}{@{}l@{}l}
0 - \dfrac{1}{a_{j_{s+1}-r,j_{s+1}}} a_{j_{s+1}-r,j} \cdot 0= 0, &\mbox{if }k=1, \ldots, s,\\
a_{j_{s+1}-r,j} - \dfrac{1}{a_{j_{s+1}-r,j_{s+1}}} a_{j_{s+1}-r,j} \cdot a_{j_{s+1}-r,j_{s+1}}= 0,~ &\mbox{if }k=s+1.\end{array}\right.\]
So the entries to the right of $(j_1-r,j_1)$, \ldots, $(j_{s+1}-r,j_{s+1})$ are zero in $AT^1 \cdots T^{s+1}$.

The fact that $(AT^1\cdots T^s)_{j_{s+1}-r+1..m,j_{s+1}} = A_{j_{s+1}-r+1..m,j_{s+1}} = 0$ implies that  $(AT^1 \cdots T^{s+1})_{j_{s+1}-r+1..m,j} = (AT^1 \cdots T^s)_{j_{s+1}-r+1..m,j}=A_{j_{s+1}-r+1..m,j}$, for $j>j_{s+1}$. In particular, $(AT^1 \cdots T^{s+1})_{j-r..m,j} = A_{j-r..m,j}$, for $j=j_{s+2}$,  \ldots, $j_t$, since $j_1<\cdots < j_t$. Applying induction, we obtain (ii).

To show unicity, we translate into equations the conditions $T$ must satisfy. Pick $j\in \{1,\ldots, m\}$ such that $j_s< j\leq j_{s+1}$. Column $j$ of $T$ has entries already determined from its general description: $t_{jj}=1$ and $t_{ij}=0$, for $i\notin \{j_1,j_2,\ldots, j_s, j\}$. The set of conditions in (ii) imply the remaining entries must satisfy the linear system
\[A_{SP_s\,P_s}T_{P_s,j} + A_{SP_s,j}=0,\]
where $SP_s=\{j_1-r,\ldots, j_s-r\}$ and $P_s=\{j_1,\ldots,j_s\}$. By hypothesis, $A_{SP_s\,P_s}$ is upper triangular with nonzero diagonal entries. Thus the linear system above has a unique solution, which must coincide with the corresponding entries of $T^1\cdots T^t$.\qed

The Row Cancellation Algorithm over $\F$ is presented below. Its properties will be the subject of the following proposition. Notice that, unlike the Incremental Sweeping Algorithm over $\F$, this one does not need a final step. This follows from the fact that the only marks assigned during the sweeping of a diagonal are primary pivots. Even if a primary pivot mark is assigned to the only entry in diagonal $m-1$, it will not provoke any changes to the connection matrix, since there are no entries to its right to cancel. Consequently, no transition matrix need be computed if $r=m-1$.\\

\textbf{\large Row Cancellation Algorithm over $\F$}

\begin{description}
\item[\textbf{Input:}] nilpotent upper triangular matrix $\Delta\in \F^{m\times m}$ with column/row partition $J_0$, $J_1$, \ldots, $J_b$.
\item[\textbf{Initialization Step:}]\mbox{}\\
 $\begin{array}{@{}l}
\left[\begin{tabular}{l}
 $r=0$\\
 $\widetilde{\Delta}^r=\Delta$\\
 $\widetilde{T}^r=I$ ($m\times m$ identity matrix) \end{tabular}\right.\end{array}$

\item[\textbf{Iterative Step:}] (Repeated until all diagonals parallel and to the right of the main diagonal have been swept)\\
$\begin{array}{@{}l}
\left[\begin{tabular}{l}
\textbf{Matrix $\Delta$ update}\\
\begin{tabular}{@{\hspace{.5cm}}l}
$r\leftarrow r+1$ \\
$\widetilde{\Delta}^r = (\widetilde{T}^{r-1})^{-1} \widetilde{\Delta}^{r-1}\widetilde{T}^{r-1}$\\
\end{tabular}\end{tabular}\right.\\
\end{array}$
\\[5pt]
$\begin{array}{@{}l}
\left[\begin{tabular}{l}
\textbf{Markup}\\
\begin{tabular}{@{\hspace{.5cm}}l}
Sweep entries of $\widetilde{\Delta}^r$ in the $r$-th diagonal: \\
\textbf{\textsf{If}} $\widetilde{\Delta}^r_{j-r,j}\neq 0$ \textbf{\textsf{and}} $\widetilde{\Delta}^r_{\textbf{\large .},j}$ does not contain a primary pivot\\
\rule{.5cm}{0pt}\textbf{\textsf{Then}} permanently mark $\widetilde{\Delta}^r_{j-r,j}$ as a primary pivot\\
\end{tabular}\end{tabular}\right.\\
\end{array}$
\\[5pt]
$\begin{array}{@{}l}
\left[\begin{tabular}{l}
\textbf{Matrix $\widetilde{T}^r$ construction}\\
\begin{tabular}{@{\hspace{.5cm}}l}
\textbf{\textsf{If}} diagonal $r$ has no primary pivots or $r=m-1$\\
\rule{.5cm}{0pt}\textbf{\textsf{Then}} $\widetilde{T}^r=I$\\
\rule{.5cm}{0pt}\textbf{\textsf{Else}} Let $j_1<\cdots<j_{t_r}$ be the indices of columns with primary pivots in\\
\rule{1cm}{0pt}diagonal $r$\\ 
\rule{1cm}{0pt}$\widetilde{T}^r = \widetilde{T}^{r,1}\cdots \widetilde{T}^{r,t_r}$, where\\
\rule{1cm}{0pt}$(\widetilde{T}^{r,s})_{i\mbox{\bf\large.}}=\left\{\begin{array}{ll} e_i,&\mbox{if }i\neq j_s,\\
e_{j_s} - \dfrac{1}{\widetilde{\Delta}^r_{j_s-r,j_s}} (0 \cdots 0~\widetilde{\Delta}^r_{j_s-r,j_s+1} \cdots \widetilde{\Delta}^r_{j_s-r,m}),\quad &\mbox{ow.}\end{array}\right.$
\end{tabular}\end{tabular}\right.\\
\end{array}$


\item[\textbf{Output:}] $(\widetilde{\Delta}^0,\ldots, \widetilde{\Delta}^{m-1})$ and $(\widetilde{T}^0,\ldots, \widetilde{T}^{m-2})$

\end{description}

To keep up the parallel with the Incremental Sweeping Algorithm over $\F$, we denote by $\widetilde{\mathcal{P}}^r$ the set of column indices of primary pivots on or below diagonal $r$ in $\widetilde{\Delta}^{m-1}$. Notice that the off-diagonal nonzero entries of $\widetilde{T}^r$ belong to rows with indices equal to the column indices of the primary pivots in diagonal $r$. Using induction, it is straightforward that the off-diagonal nonzero entries in $\widetilde{T}^1\cdots \widetilde{T}^r$ belong to rows in $\widetilde{\mathcal{P}}^r$. When $\widetilde{\Delta}^{r+1}$ is obtained in the update step, we will see in Proposition \ref{prop:SmaleOverF} that the entries to the right of the primary pivots on the $r$-th diagonal are zeroed out. Here the parallel breaks, because the changes from $\Delta^r$ to $\Delta^{r+1}$ in the Incremental Sweeping Algorithm over $\F$ arise from the cancellation of change-of-basis entries on the $r$-th diagonal, using primary pivots identified in previous iterations, and thus belonging to diagonals $1$, \ldots, $r-1$. Thus the off-diagonal nonzero entries of $T^r$ belong to $\mathcal{P}^{r-1}\backslash\mathcal{P}^{r-2}$ and, by induction, the off-diagonal nonzero entries of $T^1\cdots T^r$ belong to $\mathcal{P}^{r-1}$. Finally, the condition that the candidate for primary pivot do not have a primary pivot to its left is not required in this algorithm, due to Proposition ~\ref{prop:SmaleOverF} (vii).

The following proposition mimics the results obtained for the Sweeping Algorithm over $\Z$ and the Incremental Sweeping Algorithm over $\F$. 

\begin{prop}\label{prop:SmaleOverF} Let $(\widetilde{\Delta}^0(=\Delta), \widetilde{\Delta}^1, \ldots, \widetilde{\Delta}^{m-1})$ and $(\widetilde{T}^0=I,\widetilde{T}^1,\ldots,\widetilde{T}^{m-2})$ be the sequences of connection and transition matrices produced by the application of the Row Cancellation Algorithm over $\F$ to the connection matrix $\Delta\in \F^{m\times m}$ with column/row partition $J_0$, \ldots, $J_b$. Then, for $r=0, 1, \ldots, m-1$, we have:
\begin{itemize}
\item[(i)] $\widetilde{\Delta}^r$ is compliant with the allowable sparsity pattern of $\Delta$;
\item[(ii)]\label{prop:zerobelowppsmale} the nonzero entries of $\widetilde{\Delta}^r$ strictly below the $r$-th diagonal are either primary pivots (always nonzero) or lie above a unique primary pivot;
\item[(iii)] if $\widetilde{\Delta}^s_{p-r,p}$ is a primary pivot, then $\widetilde{\Delta}^r_{p\textbf{\large.}}$ does not contain a primary pivot, for $s\leq r$; 
\item[(iv)] if $\widetilde{\Delta}^r_{p-r,p}$ is marked as a primary pivot, then $\widetilde{\Delta}^{r+1}_{p-r,p+1..m}=0$; 
\item[(v)] if $\widetilde{\Delta}^r_{p-r,p}$ is marked as a primary pivot, then $\widetilde{\Delta}^{s}_{p\textbf{\large .}}=0$, for $s\geq r+1$.
\end{itemize}
Additionallly,
\begin{itemize}
\item[(vi)] If $\widetilde{\Delta}^{m-1}_{ij}$ is a primary pivot, then $\widetilde{\Delta}_{\textbf{\large .}i}$ has no primary pivot;
\item[(vii)] if $\widetilde{\Delta}^r_{p-r,p}$ is marked as a primary pivot, then $\widetilde{\Delta}^{s}_{p-r,p+1..m}=0$, for $s\geq r+1$; 
\item[(viii)] each row of $\widetilde{\Delta}^{m-1}$ may contain at most one primary pivot.
\end{itemize}

\end{prop}

\dem Assertion (i) is true for $r=0$. Suppose it is true for arbitrary fixed $r$. If there are no primary pivots on the $r$-th diagonal, it trivially holds for $r+1$, since $\widetilde{\Delta}^{r+1}=\widetilde{\Delta}^r$ in this case.

Now suppose $(1<)j_1<\cdots < j_{t_r}$ are the column indices of primary pivots in diagonal $r$. Then the transition matrix $\widetilde{T}^r$ is the product $\widetilde{T}^r=\widetilde{T}^{r,1}\cdots \widetilde{T}^{r,t_r}$, where 
\begin{equation}\label{eq:Trs}
(\widetilde{T}^{r,s})_{i\mbox{\bf\large.}}=\left\{\begin{array}{ll} e_i,&\mbox{if }i\neq j_s,\\
e_{j_s} - \dfrac{1}{\widetilde{\Delta}^r_{j_s-r,j_s}} (0 \cdots 0~\widetilde{\Delta}^r_{j_s-r,j_s+1} \cdots \widetilde{\Delta}^r_{j_s-r,m}),\quad &\mbox{ow.}\end{array}\right.
\end{equation}
The inverse of $\widetilde{T}^{r,s}$ is given by
\begin{equation}\label{eq:Trs-1}
(\widetilde{T}^{r,s})_{i\mbox{\bf\large.}}^{-1}=\left\{\begin{array}{ll} e_i,&\mbox{if }i\neq j_s,\\
e_{j_s} + \dfrac{1}{\widetilde{\Delta}^r_{j_s-r,j_s}} (0 \cdots 0~\widetilde{\Delta}^r_{j_s-r,j_s+1} \cdots \widetilde{\Delta}^r_{j_s-r,m}),\quad &\mbox{ow.}\end{array}\right.
\end{equation}

The connection matrix update may be computed in the following nested fashion:
\[\begin{array}{r@{\hspace{.3cm}}l@{\hspace{2cm}}r@{\hspace{.3cm}}l}
1\mbox{st step:}&\widetilde{\Delta}^r \widetilde{T}^{r,1}=\widetilde{\Delta}^{r,0} \widetilde{T}^{r,1}=\widetilde{\Delta}^{r,1}
&2\mbox{nd step:}&(\widetilde{T}^{r,1})^{-1}\widetilde{\Delta}^{r,1}=\widetilde{\Delta}^{r,2}\\
3\mbox{rd step:}& \widetilde{\Delta}^{r,2}\widetilde{T}^{r,2}=\widetilde{\Delta}^{r,3}
&4\mbox{th step:} &(\widetilde{T}^{r,2})^{-1}\widetilde{\Delta}^{r,3}=\widetilde{\Delta}^{r,4}\\
&\vdots&\end{array}\]
By induction, the matrix $\widetilde{\Delta}^r=\widetilde{\Delta}^{r,0}$ satisfies (i). We will use induction in the nested product to prove that the last matrix $\widetilde{\Delta}^{r,2t_r}=\widetilde{\Delta}^{r+1}$ also does. Suppose $\widetilde{\Delta}^{r,2s-2}$ satisfies (i). In the post-multiplication by $\widetilde{T}^{r,s}$, to column $j_s+j$ of $\widetilde{\Delta}^{r,2s-2}$ we add  $-\widetilde{\Delta}^r_{j_s-r,j_s+j}/\widetilde{\Delta}^r_{j_s-r,j_s}$ times column $j_s$, for $j=1,\ldots, m-j_s$. This coefficient is nonzero only if columns $j_s$ and $j_s+j$ belong to the same subset of the column partition. Since column $j_s$ is to the left of column $j_s+j$, and $\widetilde{\Delta}^{r,2s-2}$ is compliant with the sparsity pattern of $\Delta$, this property is inherited by the product $\widetilde{\Delta}^{r,2s-2}\widetilde{T}^{r,s}$, so $\widetilde{\Delta}^{r,2s-1}$ satisfies (i). Using (\ref{eq:Trs-1}), in the pre-multiplication by $(\widetilde{T}^{r,s})^{-1}$, to row $j_s$ of $\widetilde{\Delta}^{r,2s-1}$ we add $\widetilde{\Delta}^r_{j_s-r,j_s+j}/\widetilde{\Delta}^r_{j_s-r,j_s}$ times row $j_s+j$, for $j=1,\ldots, m-j_s$. Again the multiplication coefficient is nonzero only if the two rows belong to the same subset of the row partition. Given that row $j_s$ lies above row $j_s+j$ and that $\widetilde{\Delta}^{r,2s-1}$ is compliant with $\Delta$'s allowable sparsity pattern, it follows that $\widetilde{\Delta}^{r,2s}=(\widetilde{T}^{r,s})^{-1}\widetilde{\Delta}^{r,2s-1}$ is compliant with the allowable sparsity pattern of $\Delta$. By induction, we conclude that $\widetilde{\Delta}^{r+1}$ satisfies (i). By induction, (i) is satisfied for all $r$.

Let $f$ be the first nonzero diagonal of $\Delta$ and $(1<)j_1<\cdots<j_{t_f}$ be the columns of the nonzero entries. Then $\Delta=\widetilde{\Delta}^0=\cdots=\widetilde{\Delta}^f$. Given the construction of $\widetilde{T}^r$ and Proposition~\ref{prop:transitionmatrixSmale}, we have that $(\widetilde{\Delta}^r\widetilde{T}^r)_{j-r,j+1..m}=0$, for $j\in\{j_1,\ldots,j_{t_r}\}$.

Since there are no primary pivots to the left of diagonal $f$, (iii) is trivially valid for $s<f$. The pre-multiplication of $\widetilde{\Delta}^f\widetilde{T}^f$ by $(\widetilde{T}^r)^{-1}$ may be obtained by the sequence of pre-multiplications: by $(\widetilde{T}^{r,1})^{-1}$, by $(\widetilde{T}^{f,2})^{-1}$, etc. Therefore, given the structure of $(\widetilde{T}^{f,j})^{-1}$ in (\ref{eq:Trs-1}), only rows $j_1$, \ldots, $j_{t_f}$, are affected.
Suppose there exists $j_i\in\{j_1,\ldots,j_{t_f}\}$ such that $\widetilde{\Delta}^f_{j_i\textbf{\large .}}$ contains a primary pivot. Since entries to the left and entries below a primary pivot in $\widetilde{\Delta}^f$ are zero, and primary pivots are nonzero, we have
\[
\widetilde{\Delta}^f_{j_i-f\textbf{\large.}}\widetilde{\Delta}^f_{\textbf{\large.}j_i+f} =
\widetilde{\Delta}^f_{j_i-f,j_i}\widetilde{\Delta}^f_{j_i,j_i+f}
\neq 0,\]
contradicting the nilpotency of $\widetilde{\Delta}^f$. Therefore (iii) is also true for $s=f$.

Given the construction of $\widetilde{T}^f$ and Proposition \ref{prop:transitionmatrixSmale}, we have that
\[(\widetilde{\Delta}^f\widetilde{T}^f)_{j-f,j}=\widetilde{\Delta}^f_{j-f,j}\neq 0,\quad \mbox{for }j\in\{j_1,\ldots, j_{t_f}\}\]
and
\[(\widetilde{\Delta}^f\widetilde{T}^f)_{j-f,j+1..m}=0,\quad \mbox{for }j\in\{j_1,\ldots, j_{t_f}\}.\]
Recall that the pre-multiplication by $(\widetilde{T}^f)^{-1}$ only affects rows $j_1,\ldots, j_{t_f}$ of $\widetilde{\Delta}^f\widetilde{T}^f$, none of which contains a primary pivot. This implies that
\begin{equation}\label{eq:ppnonzero}
\widetilde{\Delta}^{f+1}_{j-f,j}=((\widetilde{T}^f)^{-1}\widetilde{\Delta}^f\widetilde{T}^f)_{j-f,j}=\widetilde{\Delta}^f_{j-f,j}\neq 0,\quad \mbox{for }j\in\{j_1,\ldots, j_{t_f}\}
\end{equation}
and
\begin{equation}\label{eq:rightppzero}
((\widetilde{T}^f)^{-1}\widetilde{\Delta}^f\widetilde{T}^f)_{j-f,j+1..m}=(\widetilde{\Delta}^f\widetilde{T}^f)_{j-f,j+1..m}=0,\quad \mbox{for }j\in\{j_1,\ldots, j_{t_f}\}.\end{equation}
Hence (iv) is valid for $r=f$.

The fact that entries to the left of primary pivots are zero, (\ref{eq:rightppzero}) and the nilpotency of $\widetilde{\Delta}^{f+1}$ imply that
\begin{equation}\label{eq:pptimesrowzero}
0=\widetilde{\Delta}^{f+1}_{j-f \textbf{\large .}}\widetilde{\Delta}^{f+1}=\widetilde{\Delta}^{f+1}_{j-f,j} \widetilde{\Delta}^{f+1}_{j\textbf{\large .}},\quad \mbox{for } j\in\{j_1,\ldots, j_{t_f}\}.\end{equation}
Using (\ref{eq:ppnonzero}) and (\ref{eq:pptimesrowzero}), we have that
\begin{equation}\label{eq:nullrow}
\widetilde{\Delta}^{f+1}_{j\textbf{\large .}}=0,\quad \mbox{for }j\in\{j_1,\ldots, j_{t_f}\},\end{equation}
which partly proves (v), for $r=f$, when $s=f+1$.

The post-multiplication part of the update, in the transition from $\widetilde{\Delta}^f$ to $\widetilde{\Delta}^{f+1}$, may only change entries strictly to the right of diagonal $f$, since entries below a primary pivot are zero. From (\ref{eq:nullrow}), we see that the pre-multiplication part of the update affects exclusively rows $j_1$, \ldots, $j_{t_f}$, which are zeroed, and so does not create new nonzero entries on or below the $f$-th diagonal. Furthermore, since rows $j_1$, \ldots, $j_{t_f}$ do not contain primary pivots by (iii), the primary pivots remain nonzero in $\widetilde{\Delta}^{f+1}$. We conclude that (ii) true for $r=f+1$.

Elementary column operations do not change a zero row, so to establish the rest of (v) we must show that the elementary row operations done in the pre-multiplication part of future updates do not change rows $j_1$, \ldots, $j_{t_f}$. Rows that change in an update are precisely those whose indices coincide with the columns that receive primary pivot marks in that iteration. But a column may contain at most one primary pivot. Therefore rows $j_1$, \ldots, $j_{t_f}$ will remain invariant for the rest of the execution of the algorithm and  we conclude that (v) is true for $r=f$. 

Suppose (ii)--(v) are true for $\widetilde{\Delta}^0$, \ldots, $\widetilde{\Delta}^r$. If all the nonzero entries of the $r$-th diagonal of $\widetilde{\Delta}^r$ lie above primary pivots, then this diagonal receives no primary pivot marks, $\widetilde{\Delta}^{r+1}=\widetilde{\Delta}^r$, and (ii)--(v) are trivially true for $r+1$.

Now consider the remaining case, where there are nonzero entries on the $r$-th diagonal of $\widetilde{\Delta}^r$ with no primary pivots below. These will be marked as primary pivots and (ii) implies all entries below these primary pivots of the $r$-th diagonal in $\widetilde{\Delta}^r$ are zero. Let $j_1<\cdots<j_{t_r}$ be the indices of the columns containing these entries. Then, using the nilpotency of $\widetilde{\Delta}^r$ and (ii), we have
\[
0=\widetilde{\Delta}^r_{j-r\textbf{\large.}}\widetilde{\Delta}^r_{\textbf{\large.}j+r} =
\sum_{k=1}^m\widetilde{\Delta}^r_{j-r,k}\widetilde{\Delta}^r_{k,j+r}=
\sum_{k=1}^{j}\widetilde{\Delta}^r_{j-r,k}\widetilde{\Delta}^r_{k,j+r},\quad\mbox{for }j\in \{j_1,\ldots, j_{t_r}\}.\]
The nonzero entries on the left of $\widetilde{\Delta}^r_{j-r,j}$ lie strictly below the $r$-th diagonal. Using the induction hypothesis (ii), these nonzero entries lie above primary pivots. But then, using (v), we obtain
\[
0=\!\!\sum_{k\in \widetilde{\mathcal{P}}^{r-1}\cap \{i\mid i<j\}}\widetilde{\Delta}^r_{j-r,k} \widetilde{\Delta}^r_{k,j+r} + \widetilde{\Delta}^r_{j-r,j}\widetilde{\Delta}^r_{j,j+r}=\widetilde{\Delta}^r_{j-r,j}\widetilde{\Delta}^r_{j,j+r}, \mbox{ for }j\in \{j_1,\ldots, j_{t_r}\}.\]
Thus $\widetilde{\Delta}^r_{j,j+r}=0$ and cannot be a primary pivot, for $j\in \{j_1,\ldots, j_{t_r}\}$. This implies (iii) is partly true for $r$, when $s=r$.

The proof of (iv) for fixed arbitrary $r$ is analogous to that for $r=f$. The structure of $\widetilde{T}^r$, $(\widetilde{T}^r)^{-1}$ and the recently established (iii) implies
\begin{equation}\label{eq:ppnonzero'}
\widetilde{\Delta}^{r+1}_{j-f,j}=((\widetilde{T}^r)^{-1}\widetilde{\Delta}^r\widetilde{T}^r)_{j-r,j}=\widetilde{\Delta}^f_{j-r,j}\neq 0,\quad \mbox{for }j\in\{j_1,\ldots, j_{t_r}\}
\end{equation}
and
\begin{equation}\label{eq:rightppzero'}
((\widetilde{T}^r)^{-1}\widetilde{\Delta}^r\widetilde{T}^r)_{j-r,j+1..m}=(\widetilde{\Delta}^r\widetilde{T}^r)_{j-r,j+1..m}=0,\quad \mbox{for }j\in\{j_1,\ldots, j_{t_r}\}.\end{equation}
Notice, in particular, that since the post-multiplication by $\widetilde{T}^r$ affects only entries to the right of diagonal $r$ and, by (iii) for $s=r$, the pre-multiplication does not alter row $j-r$, for $j\in \{j_1, \ldots, j_{t_r}\}$, we also have
\begin{equation}\label{eq:beforeppequal}
\widetilde{\Delta}^{r+1}_{j-r,1..j}=((\widetilde{T}^r)^{-1}\widetilde{\Delta}^r\widetilde{T}^r)_{j-r,1..j}=\widetilde{\Delta}^r_{j-r,1..j},\quad\mbox{for }j\in \{j_1,\ldots, j_{t_r}\}.
\end{equation}

The nilpotency of $\widetilde{\Delta}^{r+1}$ implies
\begin{eqnarray*}
0&=&\widetilde{\Delta}^{r+1}_{j-r\textbf{\large .}} \widetilde{\Delta}^{r+1}\\
&\stackrel{\textrm{(iv)}}{=}& \widetilde{\Delta}^{r+1}_{j-r,1..j} \widetilde{\Delta}^{r+1}_{1..j,\textbf{\large .}}\\
&\stackrel{(\ref{eq:beforeppequal})}{=}& \widetilde{\Delta}^r_{j-r,1..j} \widetilde{\Delta}^{r+1}_{1..j,\textbf{\large .}}\\
&\stackrel{\textrm{(ii)}}{=}&
\!\!\sum_{k\in \widetilde{\mathcal{P}}^{r-1}\cap \{i\mid i\leq j\}} \widetilde{\Delta}^r_{j-r,k} \widetilde{\Delta}^{r+1}_{k\textbf{\large .}}\\
&\stackrel{\textrm{(v)}}{=}& \widetilde{\Delta}^r_{j-r,j}\widetilde{\Delta}^{r+1}_{j\textbf{\large .}},\makebox(0,0)[l]{\rule{3cm}{0pt}for $j\in \{j_1,\ldots, j_{t_r}\}$.}
\end{eqnarray*}
The fact that the primary pivot entry $\widetilde{\Delta}^r_{j-r,j}$ is nonzero implies $\widetilde{\Delta}^{r+1}_{j\textbf{\large .}}=0$, for $j\in \{j_1,\ldots, j_{t_r}\}$. This means (iii) must be also true for $s<r+1$. Suppose not and let $j_i \in\{j_i,\ldots, j_{t_r}\}$ be such that $\widetilde{\Delta}^r_{j_i\,j}$ is a primary pivot, where $j=j_i+s<j_i+r$. By (ii), $\widetilde{\Delta}^r_{j_i\,j}\neq 0$ and the entries strictly below $\widetilde{\Delta}^r_{j_i\,j}$ are zero. But then the elementary column operations performed in the post-multiplication of $\widetilde{\Delta}^r$ by $\widetilde{T}^r$ will not change this entry, since it is located strictly to the left of diagonal $r$. So $\widetilde{\Delta}^{r+!}_{j_i\,j}\neq 0$. Likewise, the elementary row operations done in the pre-multiplication by $(\widetilde{T}^r)^{-1}$ cannot zero this entry, because to row $j_i$ one adds multiples of rows below it, and we arrive at a contradiction to the fact that $\widetilde{\Delta}^{r+1}_{j_i\textbf{\large .}}=0$. Thus (iii) is true for $r+1$.

The post-multiplication part of future updates will not affect the rows $j_1$, \ldots, $j_{t_r}$. The pre-multiplication part of future updates affects only rows whose indices correspond to a column that received a primary pivot mark at the future iteration. But each column may receive at most one primary pivot mark, so the rows $j_1$, \ldots, $j_{t_r}$ will not be changed in future iterations. Therefore (v) is true for $r+1$.

The post-multiplication part of the update from $\widetilde{\Delta}^r$ to $\widetilde{\Delta}^{r+1}$ may change entries only strictly to the right of diagonal $r$. Thus this operation does not change the primary pivots marked in this or previous iterations. The pre-multiplication part affects exclusively rows $j_1$, \ldots, $j_{t_r}$. These rows did not receive a primary pivot mark on diagonal $r$ by (iii). Suppose $j_i\in\{j_1,\ldots, j_{t_r}\}$ had received one in a previous iteration. Then the entries below this primary pivot would be zero by (ii). Thus adding multiples of rows below $j_i$ to row $j_i$ could not zero out this primary pivot, contradicting (v). Therefore the only rows that change in the  pre-multiplication part do not contain primary pivots. This implies that the primary pivots of $\widetilde{\Delta}^r$ do not change with the update step. Since the update step does not introduce new nonzero entries on the $r$-th diagonal, (ii) is true for $r+1$.

By induction, we conclude that (ii)--(v) are true for all $r$.

Suppose $\widetilde{\Delta}^{m-1}_{ij}$ is a primary pivot. Item (ii) implies $\widetilde{\Delta}^{m-1}_{ij}\neq 0$. Then (iv) implies that no primary pivot was marked in column $i$ during the sweeping of diagonals $j-i$, $j-i+1$, \ldots But if column $i$ had received a primary pivot in some earlier iteration, say as the diagonal $j-i-s$, where $s\geq 1$, was swept, then, by (v), $\widetilde{\Delta}^{j-i}_{i\textbf{\large .}i}=0$, contradicting the fact that $\widetilde{\Delta}^{j-i}_{ij}\neq 0$. Thus (vi) is true.

If $\widetilde{\Delta}^{j-i}_{ij}$ is marked as a primary pivot, then the algorithm's rules implies column $j$ will not be marked again in future iterations. Item (vii) is true for $s=r+1$ by item (iv). The validity of (vii) for larger values of $s$ follows from the facts that (a) elementary column operations involving columns $j+1$, \ldots, $m$ will not change this trailing zero pattern, so primary pivots marked in these columns in future iterations do not alter this part of the row, (b) primary pivots marked in future iterations, so in diagonals to the right of diagonal $j-i$,  in columns in the range $1..j-1$ are necessarily in rows above row $i$, and, by (ii), the entries below these primary pivots are zero, so elementary column operations caused by these markings will not interfere with the zeros to the right of the primary pivot in position $(i,j)$, and (c) by (vi) column $i$ will not receive a primary pivot mark, so row $i$ will not be altered by elementary row operations.

Finally, (viii) is a consequence of (vii) and the rules of the algorithm, since an entry must be nonzero to be eligible for receiving a primary pivot mark.
\qed

This time the complementary relationship between a column $p$ containing a primary pivot and row $p$ in the last matrix in the sequence produced by the Row Cancellation Algorithm over $\F$ is a straightforward consequence of Proposition~\ref{prop:SmaleOverF}.

\begin{cor}\label{cor:complementaritySmale} Let $\widetilde{\Delta}^{m-1}$ be the last connection matrix in the sequence produced by the application of the Row Cancellation Algorithm over $\F$ to the connection matrix $\Delta\in \F^{m\times m}$ with column/row partition $J_0$, \ldots, $J_b$. Then
\begin{equation}\label{eq:complementaritySmale}
\widetilde{\Delta}^{m-1}_{\textbf{\large .}j}\widetilde{\Delta}^{m-1}_{j\textbf{\large .}} =0,\quad \mbox{for all }j.\end{equation}
\end{cor}

\dem By Proposition~\ref{prop:SmaleOverF} (ii) for $r=m-1$, the only nonzero columns of $\widetilde{\Delta}^{m-1}$ are those containing primary pivots. But if column $p$ contains a primary pivot, then by Proposition~\ref{prop:SmaleOverF} (v), $\widetilde{\Delta}^{m-1}_{p\textbf{\large .}}=0$. Thus, in either case, (\ref{eq:complementaritySmale}) is verified. \qed

Another feature that is shared with the Incremental Sweeping Algorithm over $\F$ is the possibility of doing the matrix update in blocks.

\begin{lem}\label{lema:columnrowchangesSmale}
The matrix updates in the application of the Row Cancellation Algorithm over $\F$ to the connection matrix $\Delta\in \F^{m\times m}$ with column/row partition $J_0$, \ldots, $J_b$ can be done in a blockwise fashion as follows.
\begin{equation}\label{eq:blockupdate}
\widetilde{\Delta}^r_{J_{k-1}J_k} = (\widetilde{T}^{r-1}_{J_{k-1}J_{k-1}})^{-1}\widetilde{\Delta}^{r-1}_{J_{k-1}J_k}\widetilde{T}^{r-1}_{J_kJ_k},\quad\mbox{for }k=1,\ldots,b.\end{equation}
\end{lem}

\dem The result is nontrivial only if $\widetilde{T}^{r-1}$ is not equal to the identity matrix. Suppose therefore that columns $j_1$, \ldots, $j_{t_r}$ contain primary pivot marks at diagonal $r-1$. 

If $A$ and $B$ are two upper triangular, unit diagonal matrices in $\F^{m\times m}$ such that the positions of their nonzero entries is contained in $\cup_{k=0}^b J_k\times J_k$, then $AB$ is trivially upper triangular with unit diagonal. The entry of $AB$ in position $(i,j)$ is given by
\[(AB)_{ij} = \sum_{s=1}^m a_{is}b_{sj}.\]
The sum is zero unless there is some $s$ and $k$ such that $i$, $j$, and $s$ belong to $J_k$. But this is possible only if $i$ and $j$ belong to $J_k$. Thus the support of $AB$ is also contained in $\cup_{k=0}^b J_k \times J_k$. 

The expression of $\widetilde{T}^{r-1,s}$ in (\ref{eq:Trs}) implies $\widetilde{T}^{r-1,s}_{i,j}\neq 0$ only if $i=j$ or $i=j_s$ and $j$ belongs to the subset of the partition that contains $j_s$. Thus the support of $\widetilde{T}^{r-1,s}$ is contained in $\cup_{k=0}^b J_k \times J_k$. Using the argument in the previous paragraph and induction, we conclude that $\widetilde{T}^{r-1}$ is upper triangular, has unit diagonal and its support is contained in $\cup_{k=0}^b J_k \times J_k$. Given the structure of $(\widetilde{T}^{r-1,s})^{-1}$ in (\ref{eq:Trs-1}), the same is true about $(\widetilde{T}^{r-1})^{-1}$. This implies
\begin{eqnarray}
((\widetilde{T}^{r-1})^{-1} \widetilde{\Delta}^{r-1} \widetilde{T}^{r-1})_{J_{k-1}J_k}& =&
(\widetilde{T}^{r-1})^{-1}_{\rule{0pt}{7pt}J_{k-1}\textbf{\large .}} \widetilde{\Delta}^{r-1} \widetilde{T}^{r-1}_{\textbf{\large .} J_k} \nonumber\\
&=&(\widetilde{T}^{r-1})^{-1}_{J_{k-1} J_{k-1}} \widetilde{\Delta}^{r-1}_{J_{k-1} J_k}  \widetilde{T}^{r-1}_{J_k J_k}.\label{eq:hatupdate}
\end{eqnarray}

Lemma~\ref{lema:tecnico} implies $(\widetilde{T}^{r-1})^{-1}_{J_{k-1}J_{k-1}}=(\widetilde{T}^{r-1}_{J_{k-1}J_{k-1}})^{-1}$. This and (\ref{eq:hatupdate}) conclude the proof. \qed
 
The last ingredient to prove that the Row Cancellation Algorithm over $\F$ can be done sequentially is the analogue of Lemma~\ref{lema:independence} below. It concerns the application of the Row Cancellation Algorithm over $\F$ to a connection matrix $\Delta\in \F^{m\times m}$, with column/row partition $J_0$, \ldots, $J_b$. We let $\widetilde{\mathcal{J}}_k$ be the set of columns in $J_k$ that contain primary pivot entries in $\widetilde{\Delta}^{m-1}$, for $k=1,\ldots b$. Additionally, $\widetilde{\overline{J}}_k = J_k\backslash \widetilde{\mathcal{J}}_k$, for $k=1,\ldots b$.

\begin{lem}\label{lema:independenceSmale} Suppose the Row Cancellation Algorithm over $\F$ is applied to the connection matrix $\Delta \in \F^{m\times m}$ with column/row partition $J_0$, \ldots, $J_b$. Then the markings during the sweeping of the $r$-th diagonal of $\widetilde{\Delta}^r$ on entries in columns belonging to $J_k$ and the construction of \/$\widetilde{T}^r_{J_k J_k}$ are completely determined by the values of the entries in $\widetilde{\Delta}^r_{\widetilde{\overline{J}}_{k-1}J_k}$.
\end{lem}

\dem Let $j\in J_k$. Since rows in $\widetilde{\mathcal{J}}_{k-1}$ do not receive primary pivot marks, we may assume that $j-r\in \widetilde{\overline{J}}_{k-1}$. If $\widetilde{\Delta}_{j-r,j}=0$, it will not be marked. If $\widetilde{\Delta}^r_{j-r,j}\neq 0$, one must check whether there is a primary pivot below it. It there is, it must belong to a row in $\widetilde{\overline{J}}_{k-1}$. If there is not, then this entry will be marked as a primary pivot and contribute to one of the matrices, say $\widetilde{T}^{r,s}$ whose product constitute $\widetilde{T}^r$ is the identity with the addition of the vector $- 1/\widetilde{\Delta}^r_{j-r,j} (0 \cdots 0~\widetilde{\Delta}^r_{j-r,j+1} \cdots \widetilde{\Delta}_{j-r,m})$ to the $j$-th row. The entry $-\widetilde{\Delta}^r_{j-r,j+t}/\widetilde{\Delta}^r_{j-r,j}$ in this vector are nonnull only if $j+t$ also belongs to $J_k$. Since we already have that $j-r\in \widetilde{\overline{J}}_{k-1}$, the result is proved.\qed

The validity of the block update established Lemma~\ref{lema:independenceSmale} and the complementarity relationship between a column with a primary pivot and the row of same index given in Corollary~\ref{cor:complementaritySmale} give rise to a simplified version of the Row Cancellation Algorithm over $\F$. The Block Sequential Row Cancellation Algorithm over $\F$ is a straightforward adaptation of the Block Sequential Sweeping Algorithm over $\F$, where, at step $k$,  instead of applying the Incremental Sweeping Algorithm over $\F$ to $\Delta(k)$, one applies the Row Cancellation Algorithm over $\F$ thereto. The proof of the corresponding Uncoupling Theorem is a straightforward adaptation of the original one.

\begin{thm}[Row Cancellation Uncoupling]\label{thm:uncouplingSmale}Let $\Delta\in\F^{m\times m}$ be a connection matrix with row/column partition $J_0$, \ldots, $J_b$. Let $\widetilde{\Delta}^{m-1}$ be the matrix produced by the application of the Row Cancellation Algorithm over $\F$ to $\Delta$, and let $\widetilde{\Delta}(k)^{m-1}$, for $k=1,\ldots, b$, be the matrices obtained in the Block Sequential Row Cancellation Algorithm over $\F$ applied to $\Delta$. Then $\widetilde{\Delta}^{m-1}_{J_{k-1}J_k} = \widetilde{\Delta}(k)^{m-1}_{J_{k-1}J_k}$, for $k=1,\ldots, b$ and the collection of primary pivots encountered in the application of the Row Cancellation Algorithm over $\F$ to $\Delta(k)$, for $k=1, \ldots, b$, coincides with the primary pivots found when it is applied to $\Delta$.
\end{thm}

Theorem~\ref{thm:uncouplingSmale} significantly simplifies the next results, since it allows us to consider connection matrices with only one block, which means only elementary column operations need be performed in the matrix update step. This special instance of the Row Cancellation Algorithm over $\F$ will be called 1-Block Row Cancellation Algorithm over~$\F$. Notice that, although the proposition guarantees the equalities of the primary pivots up to $r=m-1$, this is sufficient, since the primary pivots of $\Delta^{m-1}$ and of $\Delta^m$ are equal.

\begin{prop}\label{prop:PPcoincidem}
Let $\Delta\in\F^{m\times m}$ be a connection matrix with row/column partition $J_0$, \ldots, $J_b$. Let $\Delta^1$, \ldots, $\Delta^m$ and $T^1$, \ldots, $T^{m-1}$ (resp., $\widetilde{\Delta}^1$, \ldots, $\widetilde{\Delta}^{m-1}$ and $\widetilde{T}^1$, \ldots, $\widetilde{T}^{m-2}$) be the matrices produced in the application of the Incremental Sweeping Algorithm over $\F$ (resp., Row Cancellation Algorithm over $\F$) to $\Delta$. Then the primary pivots of $\Delta^r$ and $\widetilde{\Delta}^r$ coincide in position and value, for $r=1,\ldots, m-1$.
\end{prop}

\dem By Theorem~\ref{thm:uncouplingSmale}, it is enough to prove the result for a one block connection matrix. In this case, only column operations are performed in either algorithm.  

The position and value equalities are trivially true for primary pivots on the first diagonal. Assume the equalities hold for those in diagonals $1$, \ldots, $r-1$. This implies, in particular, that $\mathcal{P}^{r-1}=\widetilde{\mathcal{P}}^{r-1}$. We will prove both equalities for a generic pivot on diagonal $r$.

Suppose $\Delta^r_{j^*-r,j^*}$ is marked as primary pivot in the application of the 1-Block Incremental Sweeping Algorithm over $\F$ to $\Delta$. This means that there are no primary pivots on its left or below it in $\Delta^r$. Since these entries lie below the $r$-th diagonal, the same is true for $\widetilde{\Delta}^r$. Furthermore, Proposition \ref{prop:CorrectnessAlgorithmOverF} (ii) and Proposition~\ref{prop:SmaleOverF} (ii) imply that
\begin{eqnarray}
\Delta^r_{j^*-r+1..m,j^*} &=&0, \label{eq:delrzero}\\
\widetilde{\Delta}^r_{j^*-r+1..m,j^*} &=&0. \label{eq:hatdelrzero}\end{eqnarray}

The update rules of the 1-Block Incremental Sweeping Algorithm over $\F$, imply that each column, say $j$, of each matrix in the sequence $\Delta^0$, $\Delta^1$, \ldots produced by the algorithm, is the sum of the original column plus a linear combination of primary pivot columns, i.e., columns that have already received a primary pivot mark, on its left. Each addition of a multiple of a primary pivot column aims to increase the number of trailing zeros in column $j$. Therefore, if $i_{\max}=\max\{i \mid \Delta^r_{ij}\neq 0\}$, then $\Delta^\ell_{\textbf{\large .}j}$ may be expressed as
\begin{equation}\label{eq:deltalj1}
\Delta^\ell_{\textbf{\large .}j} = \Delta_{\textbf{\large .}j} + \alpha_1\dot{\Delta}^{r_1}_{\textbf{\large .}j_1} + \cdots + \alpha_t \dot{\Delta}^{r_t}_{\textbf{\large .}j_t},
\end{equation}
where (the dot indicates that) $\dot{\Delta}^{r_i}$ is the first matrix in the sequence that contains a primary pivot mark in column $j_i$, $j_i<j$, $r_i<\ell$ and the primary pivot in column $j_i$ lies on a row strictly below $i_{\max}$, for $i=1, \ldots, t$.

Doing this recursively, after a finite number of steps we reach an expression involving only original columns, transforming (\ref{eq:deltalj1}) into 
\begin{equation}\label{eq:deltalj2}
\Delta^\ell_{\textbf{\large .}j} = \Delta_{\textbf{\large .}j} + \Delta_{\textbf{\large .}\mathcal{P}^{\ell,j}} \beta^{\ell,j},
\end{equation}
where $\mathcal{P}^{\ell,j}(\subset \mathcal{P}^{\ell-1})$ is the set of column indices with primary pivot entries in the submatrix $\Delta^\ell_{i_{\max}+1..m,1..j-1}$, that is, with entries in rows strictly below row $i_{\max}$ and in columns strictly to the left of column $j$, and $\beta^{\ell,j}$ is a column vector of appropriate dimension.

When $\ell=r$ and $j=j^*$, (\ref{eq:deltalj2}) becomes
\begin{equation}\label{eq:deltarj*}
\Delta^r_{\textbf{\large .}j^*} = \Delta_{\textbf{\large .} j^*} + \Delta_{\textbf{\large .}\mathcal{P}^{r,j^*}} \beta^{r,j^*}.\end{equation}

Applying this procedure to each column of $\Delta^r$ in $\mathcal{P}^{r,j^*}$ we obtain the following matrix equality
\begin{equation}
\Delta^r_{\textbf{\large .}\mathcal{P}^{r,j^*}} = \Delta_{\textbf{\large.}\mathcal{P}^{r,j^*}} B^{r,\mathcal{P}^{r,j^*}}, \label{eq:deltaplj1}\end{equation}
where $B^{r,\mathcal{P}^{r,j^*}}$ is a square matrix. Submatrix $\Delta^r_{\textbf{\large .}\mathcal{P}^{r,j^*}}$ has full column rank, since, by a column permutation, it can be cast in a column echelon format. Thus $B^{r,\mathcal{P}^{r,j^*}}$ is invertible and $\Delta_{\textbf{\large.}\mathcal{P}^{r,j^*}}$ also has full column rank. Additionally, since the primary pivots in $\Delta^r_{\textbf{\large .}\mathcal{P}^{r,j^*}}$ lie strictly below row $j^*-r$, we have have that 
\begin{equation}
\Delta^r_{j^*-r+1..m,\mathcal{P}^{r,j^*}} = \Delta_{j^*-r+1..m,\mathcal{P}^{r,j^*}} B^{r,\mathcal{P}^{r,j^*}}, \label{eq:deltaplj2}\end{equation}
and both $\Delta^r_{j^*-r+1..m,\mathcal{P}^{r,j^*}}$ and $\Delta_{j^*-r+1..m,\mathcal{P}^{r,j^*}}$ have full column rank.

It is easy to see that each column, say $j$, in each matrix of the sequence generated by the application of the 1-Block Row Cancellation Algorithm over $\F$ to $\Delta$, say $\widetilde{\Delta}^\ell$, is the sum of the original column plus a linear combination of primary pivot columns to its left. But in this algorithm, this linear combination may contain a primary pivot column whose primary pivot entry lies above the lowest nonzero entry of $\widetilde{\Delta}^\ell_{\textbf{\large .}j}$.

Since, by induction the primary pivot entries below diagonal $r$ in $\Delta^r$ and $\widetilde{\Delta}^r$ coincide in position and value, the column indices of primary pivot entries in $\widetilde{\Delta}^r_{j^*-r+1..m, 1..j^*-1}$ is precisely $\mathcal{P}^{r,j^*}$.  Then
\begin{equation}\label{eq:^deltarj*1}
\widetilde{\Delta}^r_{\textbf{\large .}j^*} = \Delta_{\textbf{\large .}j^*} + 
\sum_{j\in \mathcal{P}^{r,j^*}}\alpha_j \dot{\widetilde{\Delta}}^{r_j}_{\textbf{\large .}j} +
\sum_{\scriptsize\begin{array}{c}
j\in \mathcal{P}^{r-1}\backslash \mathcal{P}^{r,j^*} \\[-5pt]
j<j^*\end{array}} \alpha_j \dot{\widetilde{\Delta}}^{r_j}_{\textbf{\large .}j},
\end{equation}
 where, in both sums, $\dot{\widetilde{\Delta}}^{r_j}$ denotes the first matrix in the sequence that contains a primary pivot in column $j$. The second summation in the right-hand-side of (\ref{eq:^deltarj*1}) aggregates the contribution to $\widetilde{\Delta}^r_{\textbf{\large .}j^*}$ of primary pivot columns in $\mathcal{P}^{r-1}$, to the left of column $j^*$, and whose primary pivots lie strictly above $j^*-r$ (remember that, by induction, there are no primary pivots to the left of $\widetilde{\Delta}^r_{j^*-r,j^*}$).
 
Express each $\dot{\widetilde{\Delta}}^{r_j}_{\textbf{\large .}j}$ in the first summation in an analogous way, aggregating in a second and separate sum the contribution of primary pivot columns with primary pivot entries strictly above row $j^*-r$. Insert back into its respective place in (\ref{eq:^deltarj*1}) the expression for each $\dot{\widetilde{\Delta}}^{r_j}_{\textbf{\large .}j}$, for $j\in \mathcal{P}^{r,j^*}$, to obtain
\begin{equation}\label{eq:^deltarj*2}
\widetilde{\Delta}^r_{\textbf{\large .}j^*} = \Delta_{\textbf{\large .}j^*} + 
\sum_{j\in \mathcal{P}^{r,j^*}}\alpha_j \Delta_{\textbf{\large .}j}+
\sum_{j\in \mathcal{P}^{r,j^*}} \alpha'_j \dot{\widetilde{\Delta}}^{r_j}_{\textbf{\large .}j} +
\sum_{\scriptsize\begin{array}{c}
j\in \mathcal{P}^{r-1}\backslash \mathcal{P}^{r,j^*} \\[-5pt]
j<j^*\end{array}} \alpha'_j \dot{\widetilde{\Delta}}^{r_j}_{\textbf{\large .}j},
\end{equation}
where the primes in the second summation indicate that the contribution of the `second sum' of the expression of each $\dot{\widetilde{\Delta}}^{r_j}_{\textbf{\large .}j}$ has already been incorporated therein. Furthermore, although this is not explicitly shown in (\ref{eq:^deltarj*2}), $\max\{r_j\mid \alpha'_j\neq 0,~ j\in \mathcal{P}^{r,j^*}\}<\max\{r_j\mid \alpha_j\neq 0,~  j\in \mathcal{P}^{r,j^*}\}$. This follows from the fact that, when one expresses a column, say $\widetilde{\Delta}^{r_j}_{\textbf{\large .}j}$, as a sum of $\Delta_{\textbf{\large .}j}$ and a linear combination of primary pivot columns, these primary pivot columns were identified in iterations previous to the sweeping of the $r_j$-th diagonal and so are associated with lower numbered diagonals. If we repeat this procedure a sufficient (and finite) number of times, eventually all coefficients in the middle sum in the right-hand-side of (\ref{eq:^deltarj*2}) will be zero, and we will have
\begin{equation}\label{eq:^deltarj*3}
\widetilde{\Delta}^r_{\textbf{\large .}j^*} = \Delta_{\textbf{\large .}j^*} + 
 \Delta_{\textbf{\large .}\mathcal{P}^{r,j^*}} \widetilde{\beta}^{r,j^*}+
\sum_{\scriptsize\begin{array}{c}
j\in \mathcal{P}^{r-1}\backslash \mathcal{P}^{r,j^*} \\[-5pt]
j<j^*\end{array}} \gamma_j \dot{\widetilde{\Delta}}^{r_j}_{\textbf{\large .}j}.
\end{equation}

Notice that 
\begin{equation}\label{eq:secondsum}
\dot{\widetilde{\Delta}}^{r_j}_{j^*-r,j} = 0, \quad \mbox{for }j\in \mathcal{P}^{r-1}\backslash \mathcal{P}^{r,j^*},~j<j^*,
\end{equation}
since the primary pivot entries of the columns in $\mathcal{P}^{r-1}\backslash \mathcal{P}^{r,j^*}$, to the left of $j^*$, are strictly above row $j^*-r$ and entries below a primary pivot are zero, see Proposition~\ref{prop:SmaleOverF} (ii).

Using (\ref{eq:delrzero}), (\ref{eq:hatdelrzero}), (\ref{eq:deltarj*}), (\ref{eq:^deltarj*3}) and (\ref{eq:secondsum}), we have
\begin{eqnarray}
0 &=& \Delta_{j^*-r+1..m,j^*} + \Delta_{j^*-r+1..m\,\mathcal{P}^{r,j^*}}\beta^{r,j^*},\label{eq:sis}\\
0 &=& \Delta_{j^*-r+1..m,j^*} + \Delta_{j^*-r+1..m\,\mathcal{P}^{r,j^*}}\widetilde{\beta}^{r,j^*}.\label{eq:sishat}\end{eqnarray}

Since $\Delta_{j^*-r+1..m\,\mathcal{P}^{r,j^*}}$ has full column rank, the solution to the (feasible) linear systems (\ref{eq:sis}) and (\ref{eq:sishat}) is unique, so
\begin{equation}\label{eq:beta=betahat}
\beta^{r,j^*} = \widetilde{\beta}^{r,j^*}.
\end{equation}

But then, (\ref{eq:deltarj*}), (\ref{eq:^deltarj*3}), (\ref{eq:secondsum}) and (\ref{eq:beta=betahat}) imply that
\begin{eqnarray*}
\Delta^r_{j^*-r,j^*} & = & \Delta_{j^*-r,j^*} + \Delta_{j^*-r,\mathcal{P}^{r,j^*}} \beta^{r,j^*}\\
& = & \Delta_{j^*-r,j^*} + \Delta_{j^*-r,\mathcal{P}^{r,j^*}} \widetilde{\beta}^{r,j^*}\\
&=& \widetilde{\Delta}^r_{j^*-r,j^*}.
\end{eqnarray*}
We've proved both entries are equal in value. By induction hypothesis, there are no primary pivots to its left or below it in the corresponding matrices. Thus $\widetilde{\Delta}^r_{j^*-r,j^*}=\Delta^r_{j^*-r,j^*}\neq 0$ will also be marked as a primary pivot.

The proof for the converse is analogous. Thus, by induction, the result follows.\qed

Theorem~\ref{thm:unitpivotTU}, Corollary~\ref{cor:finsteadofz} and Proposition ~\ref{prop:PPcoincidem} imply that the primary pivots of $\widetilde{\Delta}^{m-1}$, the last matrix produced by the application of the Row Cancellation Algorithm over $\F$ to the order $m$ TU connection matrix $\Delta$, are either $1$ or $-1$. This implies the transition matrices (and their inverses) produced by the application of the Row Cancellation Algorithm over $\F$ to a TU connection matrix are all integral. This justifies the application of the Row Cancellation Algorithm over $\F$ to TU connection matrices.  Hence, we have the following result.


\begin{thm}[Equality of primary pivots]\label{thm:smale=incr}
Let $\Delta\in\{-1,0,1\}^m$ be a TU connection matrix with column/row partition $J_0$, \ldots, $J_b$. The primary pivots marked in the application of the Sweeping Algorithm over $\Z$ to $\Delta$ coincide in value and position with the primary pivots marked in the application of Row Cancellation  Algorithm over $\F$ thereto.
\end{thm}

\begin{figure}[htbp]
\begin{center}
Surface connection matrix $\Delta =$\raisebox{-.17\textwidth}{\includegraphics[height=.4\textwidth,viewport=15 0 200 170,clip]{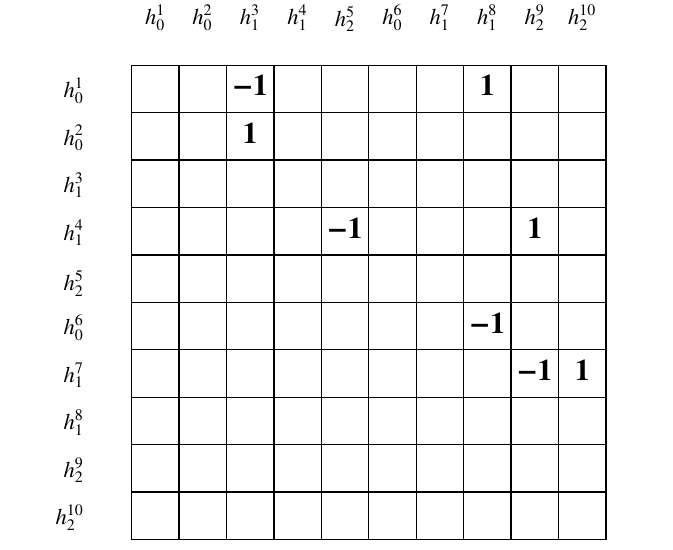}}\\
\begin{tabular}{cc}
Incremental Sweeping Algorithm over $F$ & Smale's Cancellation Algorithm over $\F$\\
\includegraphics[height=.4\textwidth]{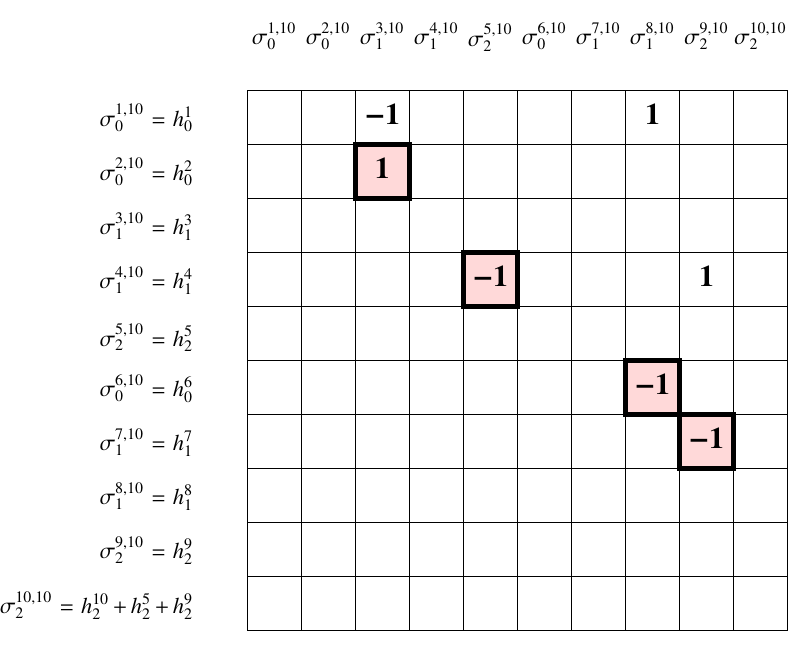}&%
\includegraphics[height=.4\textwidth]{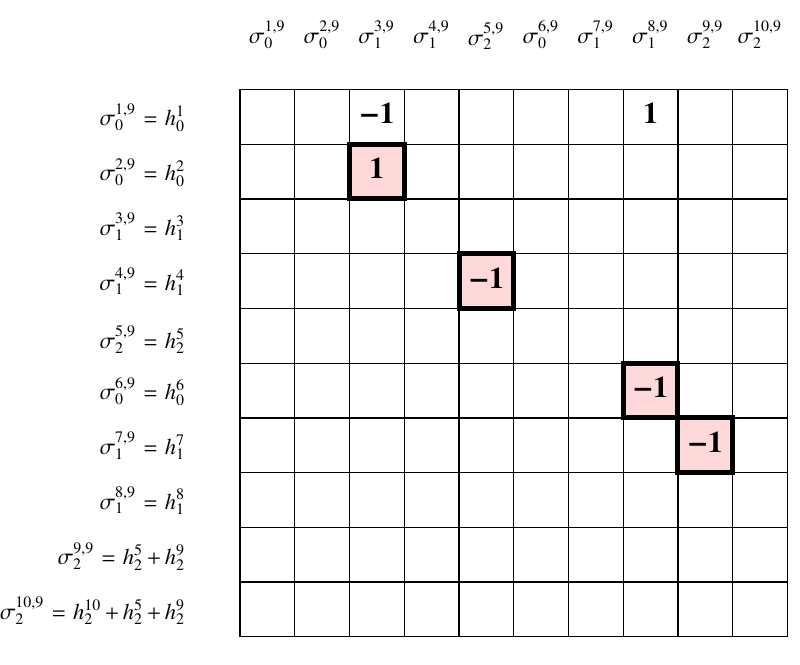}%
\end{tabular}
\end{center}
\caption{Surface connection matrix and the last matrices produced by the application thereto of the Incremental Sweeping Algorithm over $\F$ and Smale's Cancellation Algorithm over $\F$.}\label{fig:IncrXSmale}
\end{figure}

Notice that, although the primary pivots of $\Delta^m$ and $\widetilde{\Delta}^{m-1}$ coincide, the matrices might differ, as exemplified in Figure~\ref{fig:IncrXSmale}. In this figure, cells containing primary entries are highlighted.

\section{ Dynamical interpretation for the algebraic cancellations}\label{surfaces}

\subsection{Totally unimodular connection matrices }\label{subsec:dynamicTU}

In this section we will prove Theorem \ref{ordering-cancellation-TU}.

\vspace{0.3cm}

\dem {\bf (Theorem \ref{ordering-cancellation-TU} )} \ \ 
Recall that $M$  is a closed simply connected manifold of dimension $m>5$ and  $f:M \rightarrow \mathbb{R}$  a Morse-Smale function. Also  $(C,\Delta)$ is a filtered Morse chain complex with finest filtration and   $(E^{r},d^{r})$ is the associated spectral sequence. Assume that $\Delta$ is a totally unimodular matrix.

As  proven in \cite{Cornea:2010}, the non zero differentials  of the spectral sequence  are induced by the pivots. When working with TU connection matrices, it follows from  Theorem \ref{thm:unitpivotTU} (Primary Pivots for TU Connection Matrices)
  that  the primary pivots are always equal to $\pm 1$. Hence, the differentials $d^{r}_p:E^r_{p}\to E^r_{p-r}$ associated to primary pivots are isomorphisms and the ones associated to change of basis pivots always correspond to zero maps. In fact, if a differential $d^{r}_p:E^r_{p}\to E^r_{p-r}$ corresponds to a change of basis pivot, then there is a primary pivot in row $p-r$ and thus $E^r_{p-r}=0$. Consequently, the non zero differentials  are isomorphisms and this
implies that at the next stage of the spectral sequence they produce
{\it algebraic cancellations}, i.e. $E^{r+1}_{p}=E^{r+1}_{p-r}=0$.

The Equality of Primary Pivots Theorem (Theorem \ref{thm:smale=incr})  proves that primary pivots marked when  applying the Sweeping Algorithm over $\Z$ to $\Delta$ coincide in value and position with  primary pivots marked when applying the Row Cancellation Algorithm over $\F$  thereto. Therefore, the algebraic cancellations of the modules $E^r$ of the spectral sequence determined by the SSSA are in one-to-one correspondence with 
the primary pivots determined by the RCA.

In what follows, we analyze the dynamics within   the sequence of matrices obtained by the RCA. Our approach is based on 
Theorem 2.16 in \cite{Franks:1982}, which  states that:
 {\it
Let $M$ be a compact simply connected manifold of dimension $m>5$. Given a finitely generated free chain complex $C$ with $H_{\ast}(C)\approx H_{\ast}(M)$, then there exists a self-indexing Morse function $g$ on $M$ such that if $M_{k}=g^{-1}((-\infty, k+1/2))$ then the chain complex $\{ H_{k}(M_{k},M_{k-1}),\partial_{k}\}${\footnote{The map $\partial_{k}:H_{k}(M_{k},M_{k-1})\rightarrow H_{k-1}(M_{k-1},M_{k-2})$ is the boundary map of the triple $(M_{k},M_{k-1},M_{k-2})$.}} is isomorphic to $C$.}

Let $\{\widetilde{\Delta}^{r}\}_{r=0}^{m}$ be the sequence of matrices produced when one applies the Row Cancellation  to $\Delta$. For each $r$, denote by $\widehat{\Delta}^{r}$ the  matrix obtained from $\widetilde{\Delta}^{r}$ by removing the rows  and columns $(p+1)$ and $(p-\xi+1)$,  for each primary pivot $\widetilde{\Delta}^{r}_{p-\xi+1,p+1}$ in the $\xi$-th diagonal of $\widetilde{\Delta}^{r}$, for $\xi=1,\dots,r$.  Observe that the last matrix $\widehat{\Delta}^{m}$ is null, since all non zero entries of $\widetilde{\Delta}^{m}$ are above a primary pivot.

Now, for each $r=1,\dots,m$, consider the pair  $(C(r),\widehat{\Delta}^{r})$ where $C(r)$ is generated by the subset of $Crit(f)$ consisting of all the critical points of $f$  except for the ones cancelled  in the previous step $(r-1)$. The pair $(C(r), \widehat{\Delta}^{r})$ is a chain complex whose homology coincides with the homology of $M$.  In fact, it follows from  Theorems \ref{thm:unitpivotTU} and \ref{thm:smale=incr} that the primary pivots marked during RCA are $\pm 1$, hence, each  change of basis in RCA
 is a change of basis over $\mathbb{Z}$. Moreover,   in order to construct $C(r)$,    only   pairs of cancelling  critical points are removed and since these do not correspond to generators of $H(M)$, then  the homology of $(C(r), \widehat{\Delta}^{r})$ coincides with the singular homology of $M$.

 By Theorem 2.6 in  \cite{Franks:1982}, there is 
a self-indexing Morse function $g_{r}$ on $M$ such that the chain complex $\{ H_{k}(M_{k}(r),M_{k-1}(r)),\partial_{k}^{r}\}$ is isomorphic to $(C(r),\widehat{\Delta}^{r})$ , where $M_{k}(r)=g_{r}^{-1}((-\infty, k+1/2))$. The chain complex $\{ H_{k}(M_{k}(r),M_{k-1}(r)),\partial_{k}^{r}\}$ is in fact a Morse chain complex for the function $g_{r}$, see \cite{Banyaga:2004}. Consider the Morse flow $\varphi_{r}$  associated to the vector field $-\nabla g_r$. Moreover, since the last matrix $\widetilde{\Delta}^{m}$ produced by the RCA generates  a null matrix $\widehat{\Delta}^{m}$, then the flow $\varphi_{m}$ corresponds to the gradient flow associated to a perfect Morse function $g_{m}$.
\qed

\subsection{Surface connection matrices}\label{subsec:Smale}

In this section, we will consider the case where $M$ is a surface. A surface connection matrix has some additional  properties which enable us to prove stronger dynamical results.  In this case, each algebraic cancellation of the modules of the spectral sequence, codified by the Spectral Sequence Sweeping Algorithm, can be  interpreted  as a dynamical cancellation.

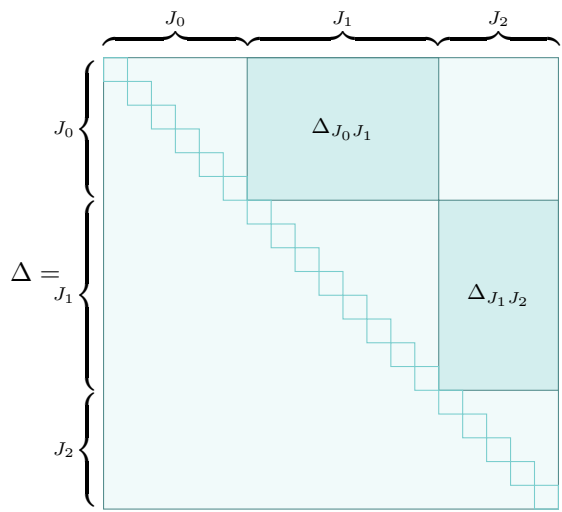
\begin{figure}[htbp]
\centerline{%
\begin{tikzpicture}[x=0.3cm,y=0.3cm,xscale=1.05,yscale=1.05]
\useasboundingbox (0,0) rectangle (19,22);
\draw (3,19) node[above]{$\overbrace{\rule{1.9cm}{0pt}}^{J_0}$};
\draw (10,19) node[above]{$\overbrace{\rule{2.5cm}{0pt}}^{J_1}$};
\draw (16.5,19) node[above]{$\overbrace{\rule{1.6cm}{0pt}}^{J_2}$};\draw (-1.4,10) node[left] {$\Delta =$};
\draw[color=DarkSlateGray4,fill=DarkSlateGray3,fill opacity=0.1] (0,0)--(0,19)--(19,19)--(19,0)--cycle;
\draw[color=DarkSlateGray4,fill=DarkSlateGray3,fill opacity=0.25] (6,19)--++(8,0)--++(0,-14)--++(5,0)--++(0,8)--++(-13,0)--cycle;
\draw (10,16) node{\footnotesize $\Delta_{J_0 J_1}$}
(16.5,9) node{\footnotesize$\Delta_{J_1 J_2}$}
(.1,16) node[rotate=90,above]{$\overbrace{\rule{1.8cm}{0pt}}$}
(-.8,16) node[left]{\scriptsize$J_0$}
(.1,9) node[rotate=90,above]{$\overbrace{\rule{2.5cm}{0pt}}$}
(-.8,9) node[left]{\scriptsize$J_1$}
(.1,2.5) node[rotate=90,above]{$\overbrace{\rule{1.5cm}{0pt}}$}
(-.8,2.5) node[left]{\scriptsize$J_2$};
\foreach \x in {0,...,18}
\draw[shift={(\x,-\x)},color= DarkSlateGray3]
(0,19)--(1,19)--(1,18)--(0,18)--cycle;
\end{tikzpicture}}
\caption{Grouped surface connection matrix with $J_0=\{1,\ldots,6\}$, $J_1=\{7,\ldots,14\}$ and $J_2=\{15,\ldots, 19\}$.}\label{fig:deltatipica}
\end{figure}

Surfaces connection matrices are endowed with special features unique to the two-dimensional case satisfying  the following properties:\begin{itemize}
\item[(i)] $\Delta_{ij}\in \{0,1,-1\}$;
\item[(ii)] columns and rows of $\Delta^0$ may be partitioned into three groups, namely $J_0$, $J_1$ and $J_2$, the first associated with wells ($h_0$'s), the second with saddles ($h_1$'s) and the third with sources ($h_2$'s). Figure~\ref{fig:deltatipica} shows the typical structure of grouped surface connection matrix. Block $\Delta_{J_0 J_1}$ contains the connections from saddles to wells, while block $\Delta_{J_1J_2}$ contains the connections from sources to saddles;
\item[(iii)] each column in $\Delta_{J_0 J_1}$ contains either two nonzero elements, namely $1$ and $-1$, or none;
\item[(iv)] each row in $\Delta_{J_1 J_2}$ contains either two non zero elements, namely $1$ and $-1$, or none;
\end{itemize}
or is obtained from a matrix with the properties above by multiplying a subset of rows and/or columns by $-1$.

The proof of this characterization for surfaces connection matrices can be found in \cite{Bertolim:2013}. This characterization allows us to prove the following result.

\begin{lem} \label{lema:matrizTU} A surface connection matrix is TU.
\end{lem}

\dem It suffices to show that a matrix satisfying (i)--(iv) is TU, since this property is invariant under multiplication of rows and/or columns by $-1$.

First notice that $\Delta_{J_0J_1}$ (resp., $\Delta_{J_1J_2}$) is TU because it is a $0$, $\pm 1$ such that each column (resp., row) contains either none or two non zeros of opposite signs, see~\cite{Schrijver:1986}. Let $M=\Delta_{LC}$ be a square submatrix of $\Delta$. We show by induction on the order of $M$ that its determinant is $0$, $\pm1$. The statement is trivially true when $M$ is $1\times 1$. Assume it is true for $n\times n$ submatrices and consider a submatrix of order $n+1$. Let $L_k=L\cap J_k$ and $C_k=C\cap J_k$, for $k=0$, $1$, $2$.  If $L_2\neq \emptyset$ or $C_0\neq \emptyset$, the determinant of $M$ will be zero, and the statement will be true. So suppose $L_2=C_0=\emptyset$. If $L_0=\emptyset$ (resp., $L_1=\emptyset$, $C_1=\emptyset$, $C_2=\emptyset$), then $M$ is a submatrix of $\Delta_{J_1\mbox{\large\textbf{.}}}$ (resp., $\Delta_{J_0\mbox{\large\textbf{.}}}$, $\Delta_{\mbox{\large\textbf{.}}J_2}$, $\Delta_{\mbox{\large\textbf{.}}J_1}$), a TU matrix, and its determinant is $0$, $\pm1$. So the only case that remains to be analyzed is when $L_0$, $L_1$, $C_1$ and $C_2$ are all nonempty. Notice that $M_{L_1C_1}$ and $M_{L_0C_2}$ are null.

Consider the submatrix $M_{L_0C_1}$. If all its columns contain two non zeros, then the sum of the rows in $M_{L_0\mbox{\large\textbf{.}}}$ is zero and $\det M=0$, proving the assertion. If not, there is a column that is either null, implying $M$ has a null column and thus null determinant, or contains exactly one non zero entry, say $m_{ij}$. Since $M_{L_1C_1}$ is zero, we can apply Laplace's expansion on this column and express $\det M$ as $\pm m_{ij}$ times the determinant of the submatrix $\widetilde{M}$ of $M$ obtained by removing row $i$ and column $j$ therefrom. But the induction hypothesis may be applied to $\widetilde{M}$, since it is of order $n$, and thus prove the assertion for the remaining case.\qed

Lemma~\ref{lema:matrizTU} implies Theorem~\ref{thm:unitpivotTU}  is valid for surface connection matrices.

\begin{cor}[Primary pivots for orientable surfaces]\label{cor:unitpivotsurface} Let $\Delta$ be a surface connection matrix. Then the primary pivots obtained when applying the Incremental Sweeping Algorithm over $\F$ thereto value $\pm 1$.
\end{cor}

As consequence of Corollary \ref{cor:unitpivotsurface}, every non zero differential is an isomorphism and hence it produces a algebraic cancellation  of the modules of the spectral sequence. Our approach  in \cite{Bertolim:2013} was to interpret the algebraic cancellation of the modules of the spectral sequence, which has been coded by the Spectral Sequence Sweeping Algorithm, as dynamical cancellations.

  Whenever a dynamical cancellation occurs, all the connecting orbits to both singularities cancelled must disappear immediately.  The Row Cancellation Algorithm  reflects exactly this situation and hence it is better suited to relate  dynamical interpretation to algebraic cancellations.   
  When the input to the Row Cancellation Algorithm over $\F$ is restricted to the special class of surface connection matrices, it is called \emph{Smale's Cancellation Sweeping Algorithm}.
  
 It follows from Lemma~\ref{lema:matrizTU} and Theorem~ \ref{thm:smale=incr} that, for a surface connection matrix,  the primary pivots marked in  both the Sweeping Algorithm over $\Z$ and  Smale's Cancellation Sweeping  Algorithm coincide in value and position. 

\begin{cor}[Equality of primary pivots for orientable surfaces]\label{cor:smale=incr}
Let $\Delta\in\{-1,0,1\}^m$ be a surface connection matrix with column/row partition $J_0$, $J_1$, $J_2$. The primary pivots marked in the application of the Sweeping Algorithm over $\Z$ to $\Delta$ coincide in value and position with the primary pivots marked in the application of Smale's Cancellation Sweeping  Algorithm thereto.
\end{cor}

%
%

%
%

 We reproduce Theorem $5.1$ in  \cite{Bertolim:2013}, refereed therein as the Ordered Smale's Cancellation Theorem, which asserts:

\begin{thm}\label{ordering-cancellation}Let $(C,\Delta)$ be the Morse chain complex associated to a Morse-Smale function $f$. Let $(E^r,d^r)$ be the associated spectral sequence for the finest filtration $F=\{F_pC\}$ defined by $f$. The algebraic cancellation of the modules $E^r$ of the spectral sequence are in one-to-one correspondence with dynamical cancellations of critical points of $f$. 
\end{thm}

The proof of this theorem, sketched out in \cite{Bertolim:2013}, relies  on Corollary \ref{cor:unitpivotsurface} and Corollary \ref{cor:smale=incr} proved in this article.



In conclusion,  Theorem \ref{thm:smale=incr} is 
crucial in producing the association  of the algebraic cancellation of modules of a spectral sequence of a filtered  chain complex and the dynamical cancellation of the critical points that generates this chain complex.
See  \cite{Bertolim:2013} for more details.

\section*{Final Remarks}

In this article as well as in \cite{Cornea:2010,Rezende:2010}, we have explored the algebraic tool provided by the spectral sequence and its dynamical implications. 
This was a major step in investigating the dynamics associated to the spectral sequence. In particular, we obtined results on  the efect the change of basis of the connection matrices had on the changes in generators of the modules  $E^r_p$ which were coded in the connection matrices determined by the SSSA in terms of a continuation of the initial Morse decomposition.  Many questions arise in this algebraic-dynamical setting.

The work developed in \cite{Cornea:2010} for spectral sequences $(E^r,d^r)$ where each $E^r_p$ is a $\mathbb {Z}$-module admits the possibility of modules with torsion which may disappear through some algebraic cancellation as one calculates the spectral sequence or remain after the sequence stabilizes. What is the dynamical meaning of the torsion in these two different cases? Note that torsion never appears for TU matrices. The following example illustrates this phenmena.

Inspired by  the results in \cite{Bertolim:2013} and in this article, where the algebra has its dynamical correspondence and is mostly determined by the primary pivots, we mean to explore in the case aforementioned, primary and change-of-basis pivots which are not necessarily equal to  $\pm 1$. Our motivation for this investigation is that certain primary and change-of-basis pivots correspond to nonzero differentials of the spectral sequence as proved in \cite{Cornea:2010}.
 

\vspace{0.5cm}

\noindent\hspace{-1cm}\includegraphics[height=7.7cm]{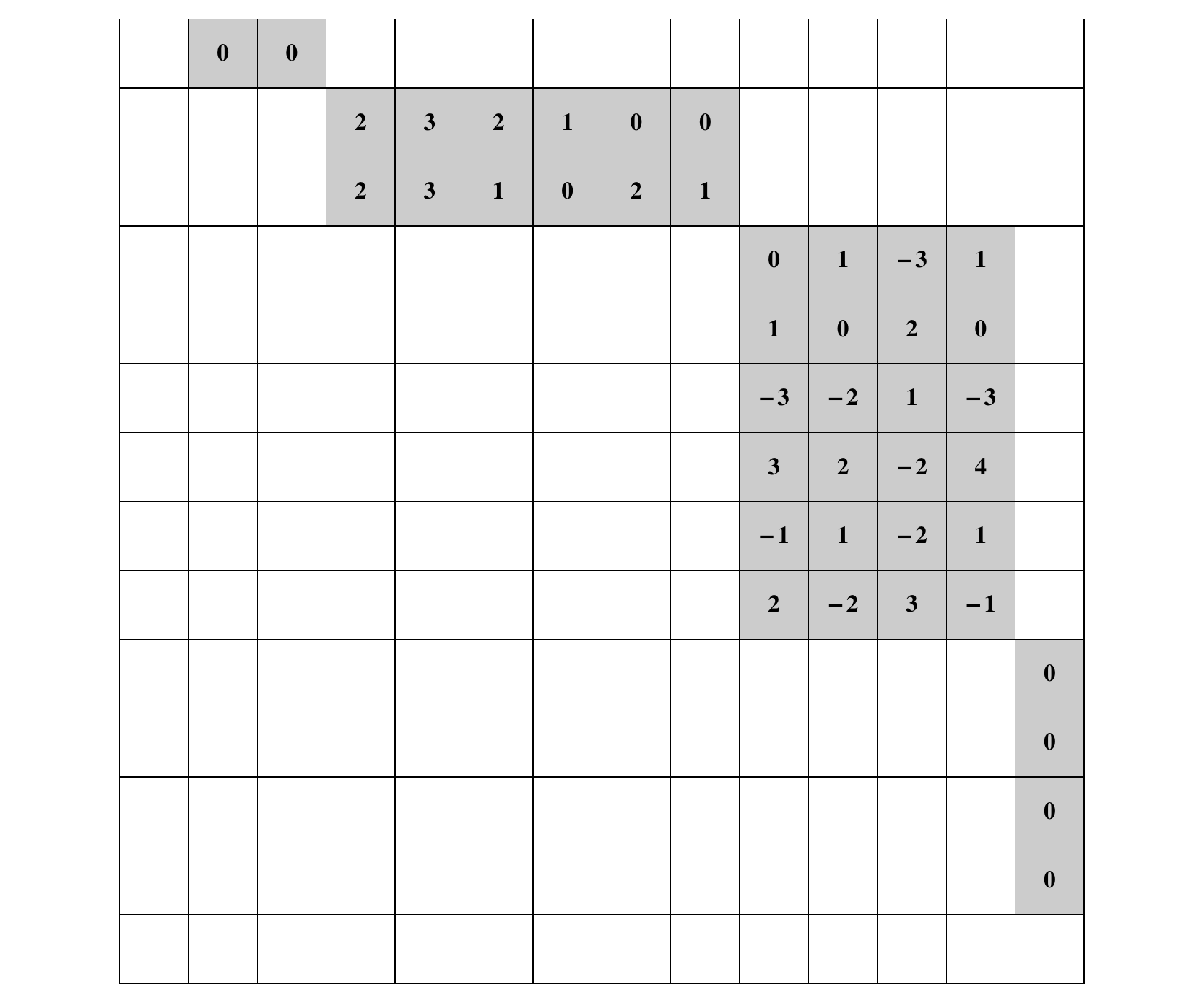}
\includegraphics[height=7.7cm]{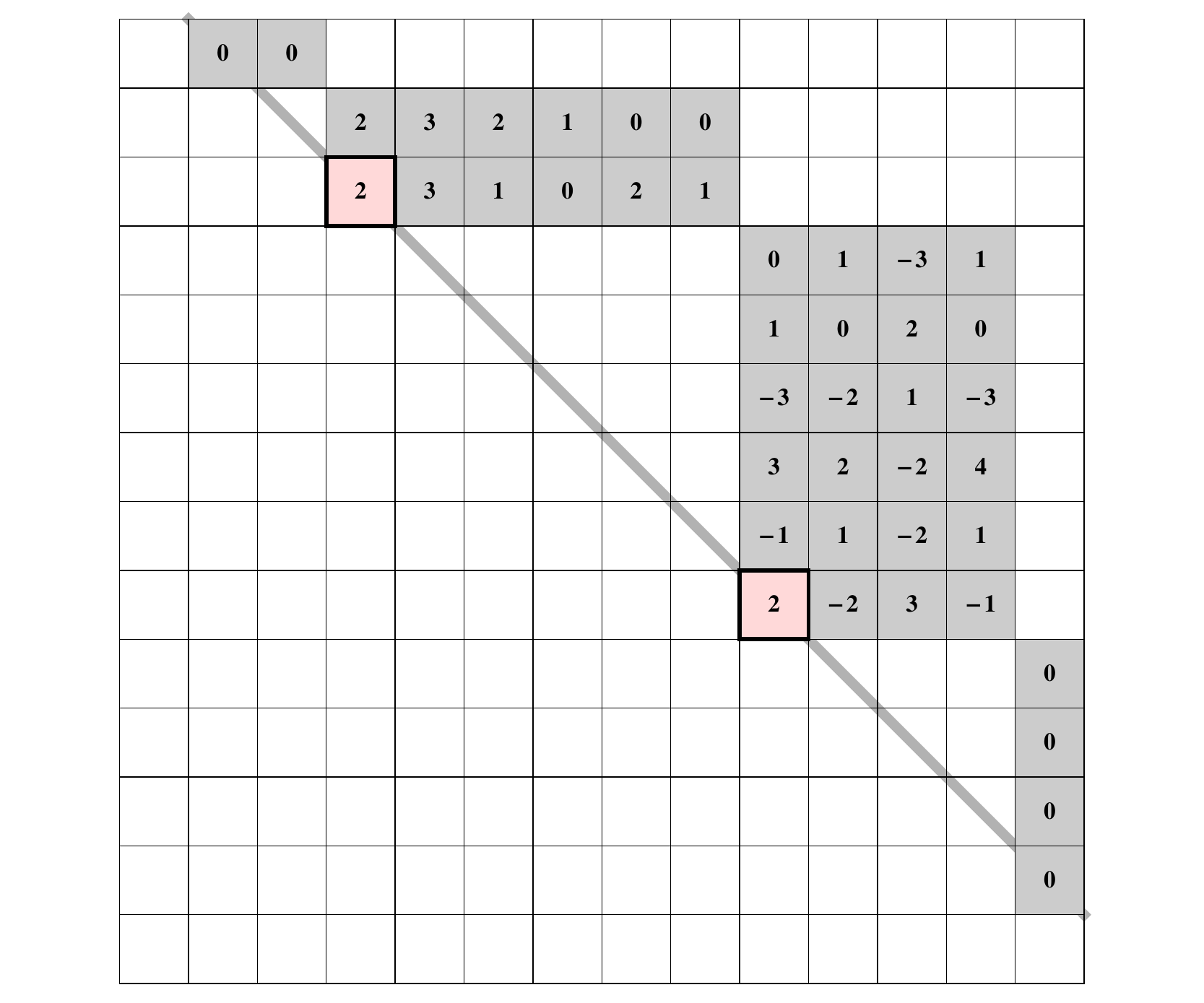}

\noindent\hspace{-1cm}\includegraphics[height=7.7cm]{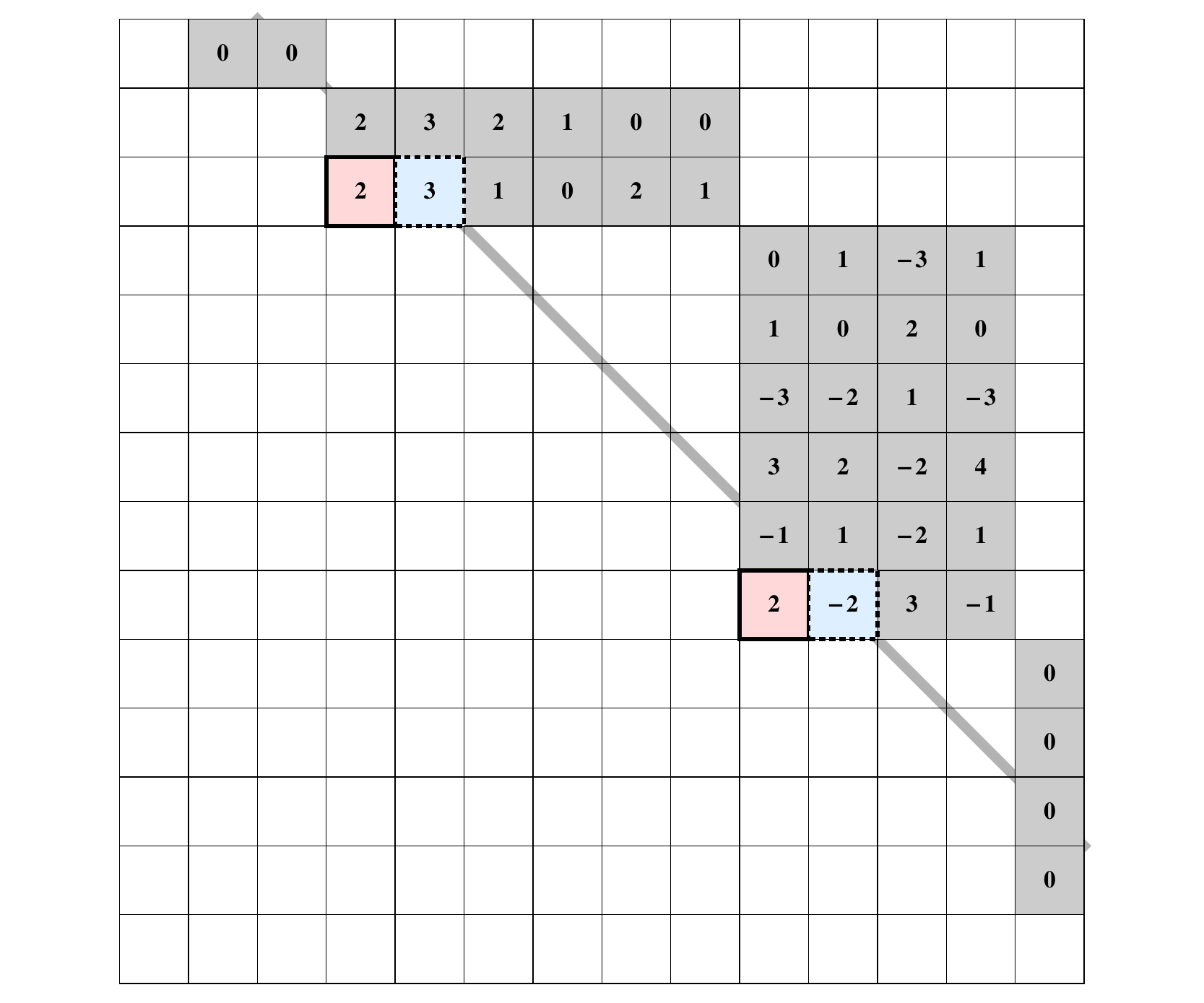}
\includegraphics[height=7.7cm]{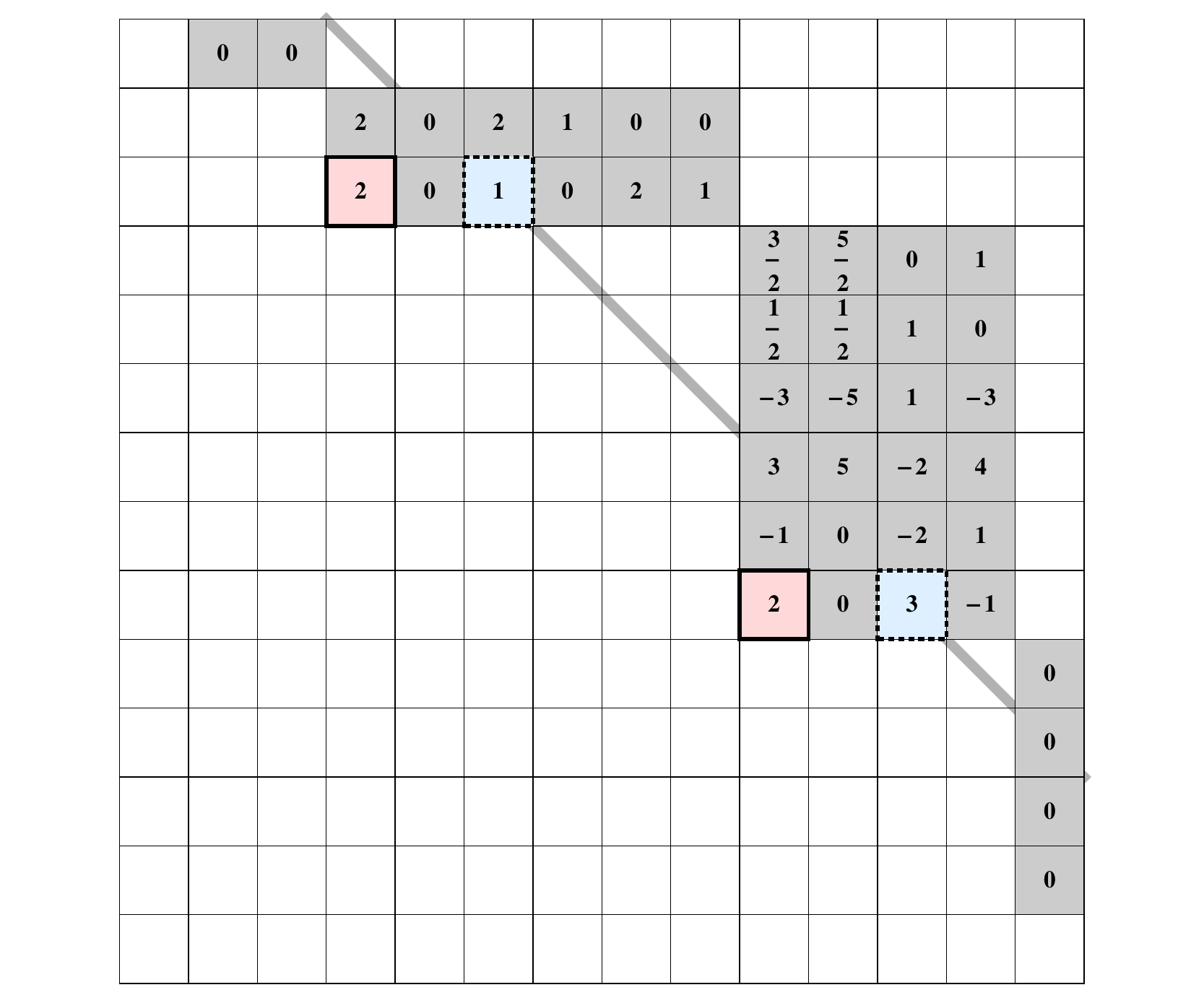}

\vspace{1cm}

\noindent\hspace{-1cm}\includegraphics[height=7.7cm]{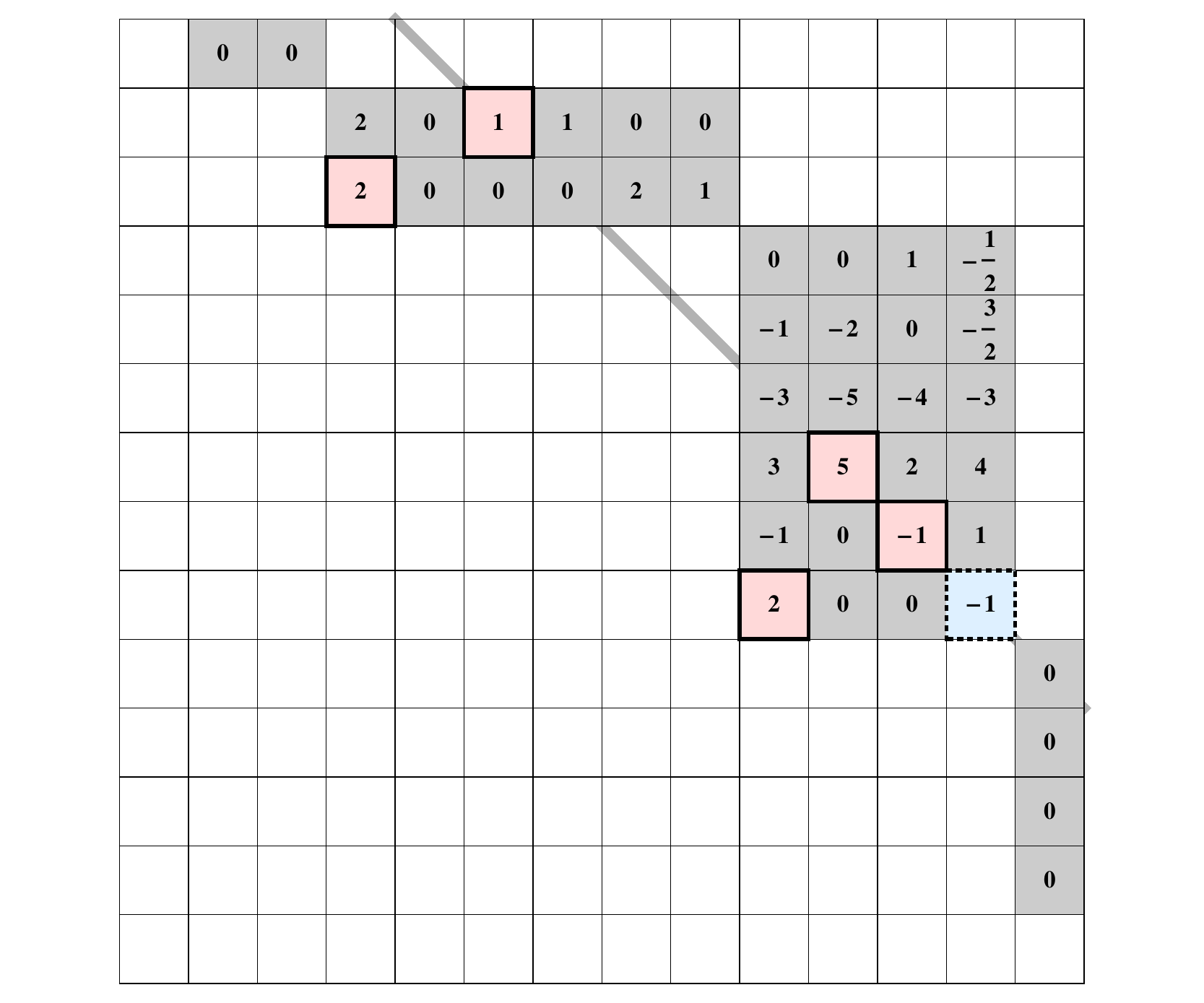}
\includegraphics[height=7.7cm]{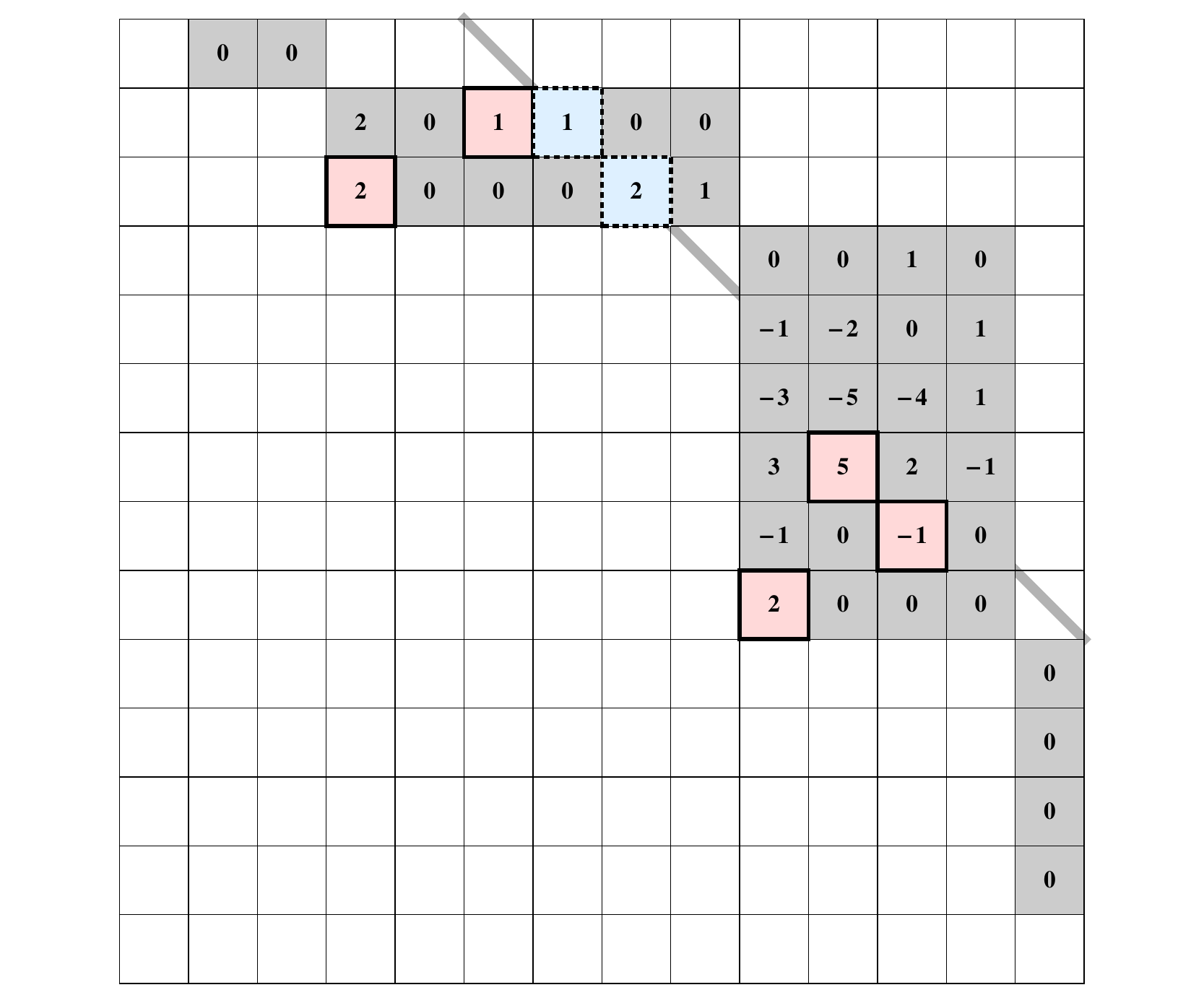}

\vspace{1cm}

\noindent\hspace{-1cm}\includegraphics[height=7.7cm]{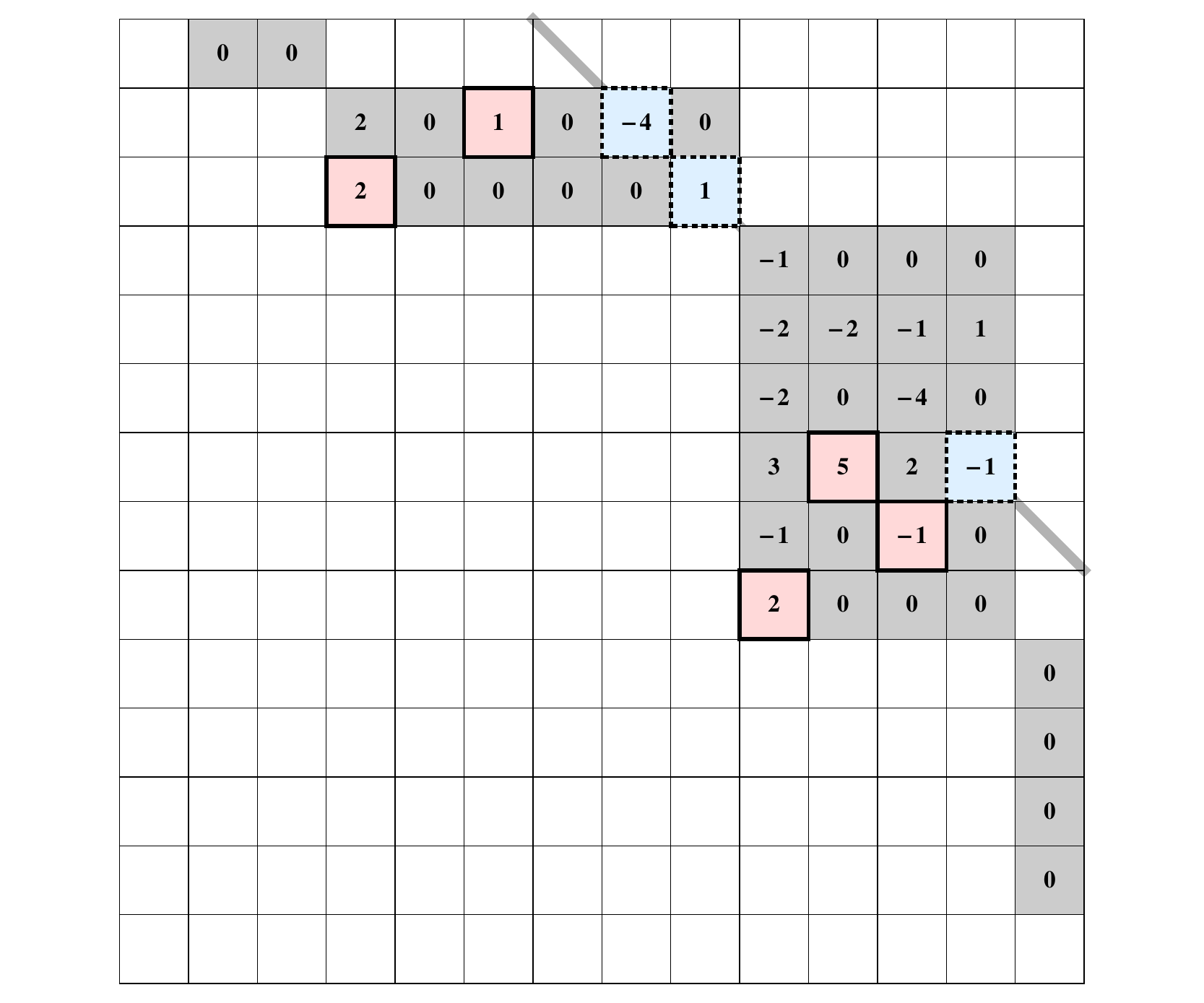}
\includegraphics[height=7.7cm]{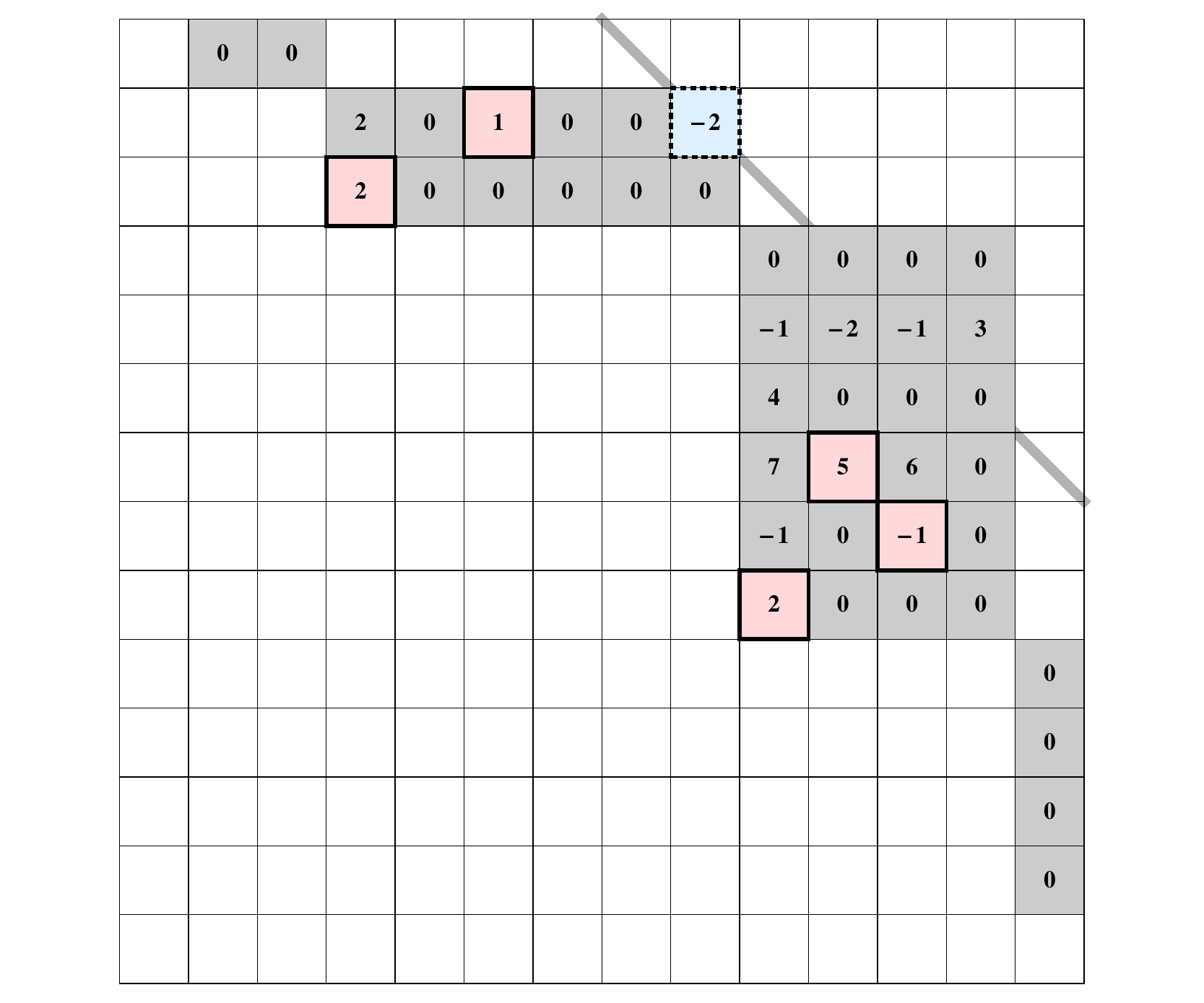}

\vspace{1cm}

\noindent\hspace{-1cm}\includegraphics[height=7.7cm]{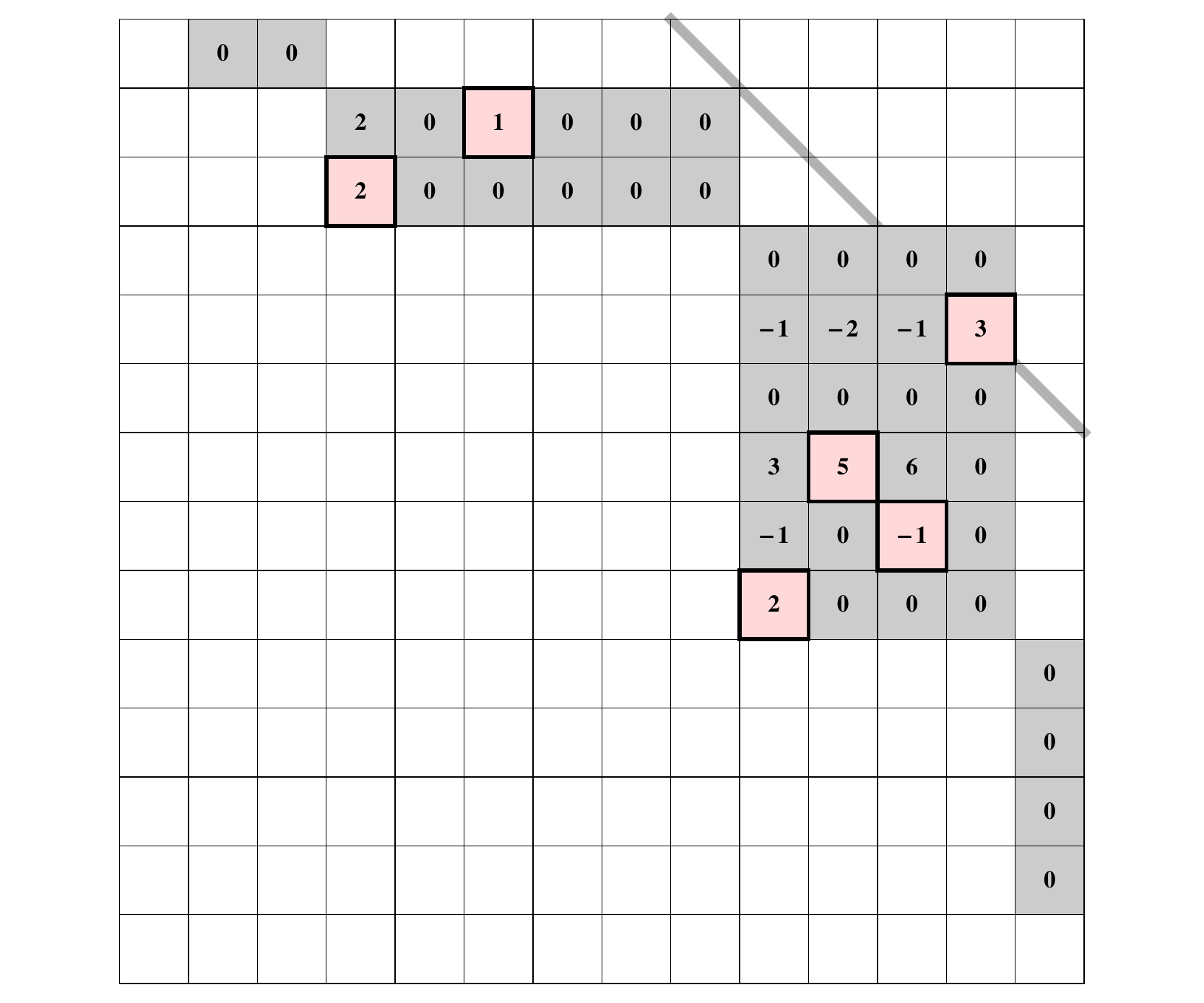}
\includegraphics[height=7.7cm]{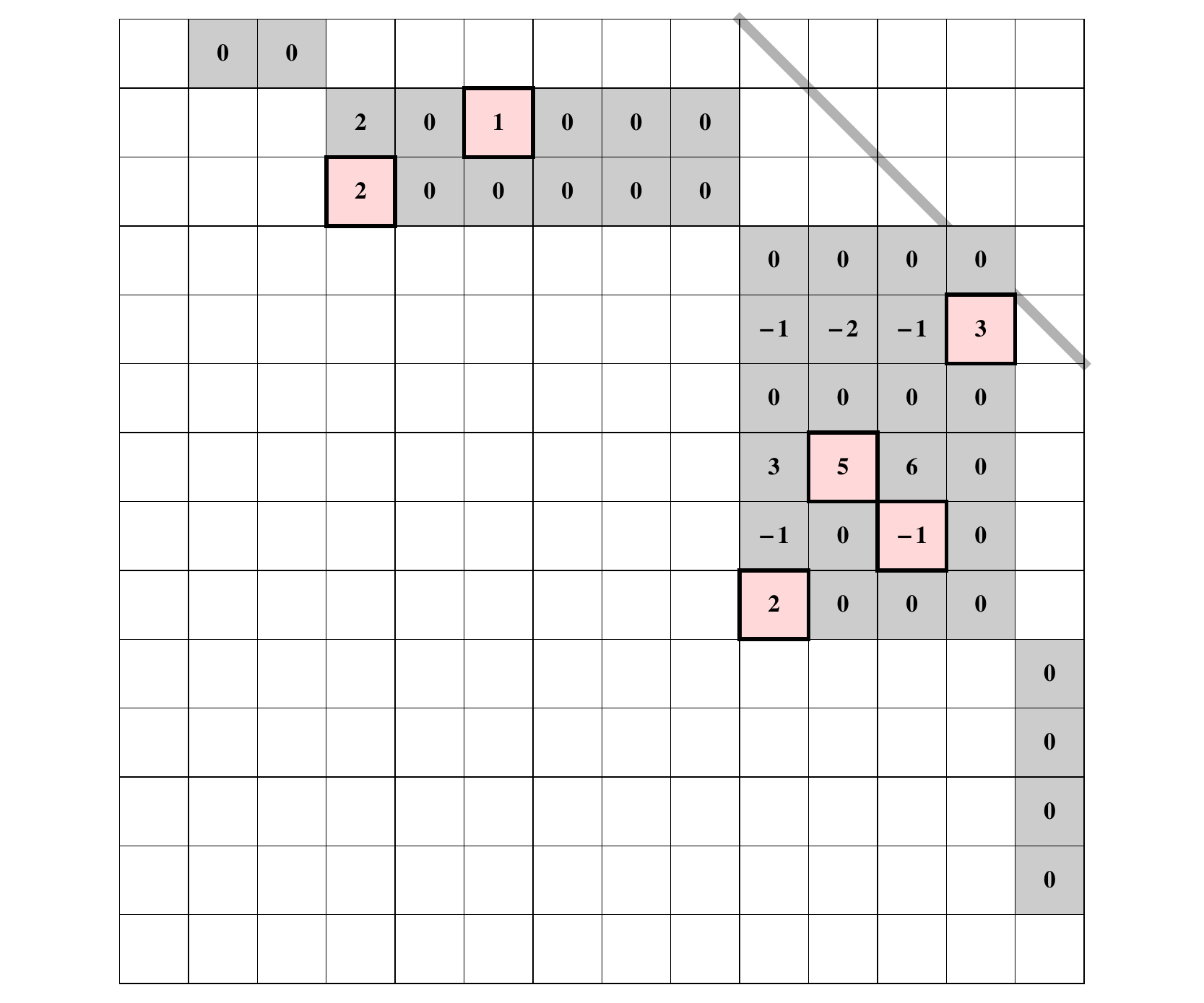}

\vspace{1cm}

As we can see in the example, as $r$ increases, the $\mathbb {Z}$-module $E^{r}_{p}$ changes generators and the SSSA brings this about when a change-of-basis pivot is marked, causing a change of basis over  $\mathbb {Q}$ of the matrix $\Delta^r$ in order to determine the connection matrix $\Delta^{r+1}$ of the next stage. Note that in the example, even though the entries in the intermediary matrices can be fractional, the primary pivots and change-of basis pivots are always integers. Furthermore, the entries of the last matrix are all integers. This is true in general and was proved in  \cite{Rezende:2010}. However, many questions remain.
Up to what point does the SSSA associate a continuation of the initial flow to the unfolding of the algebra? What is the dynamical meaning of the matrices
of the intermediary stages? For instance, what do the rational entries mean?

\includegraphics[height=7.3cm]{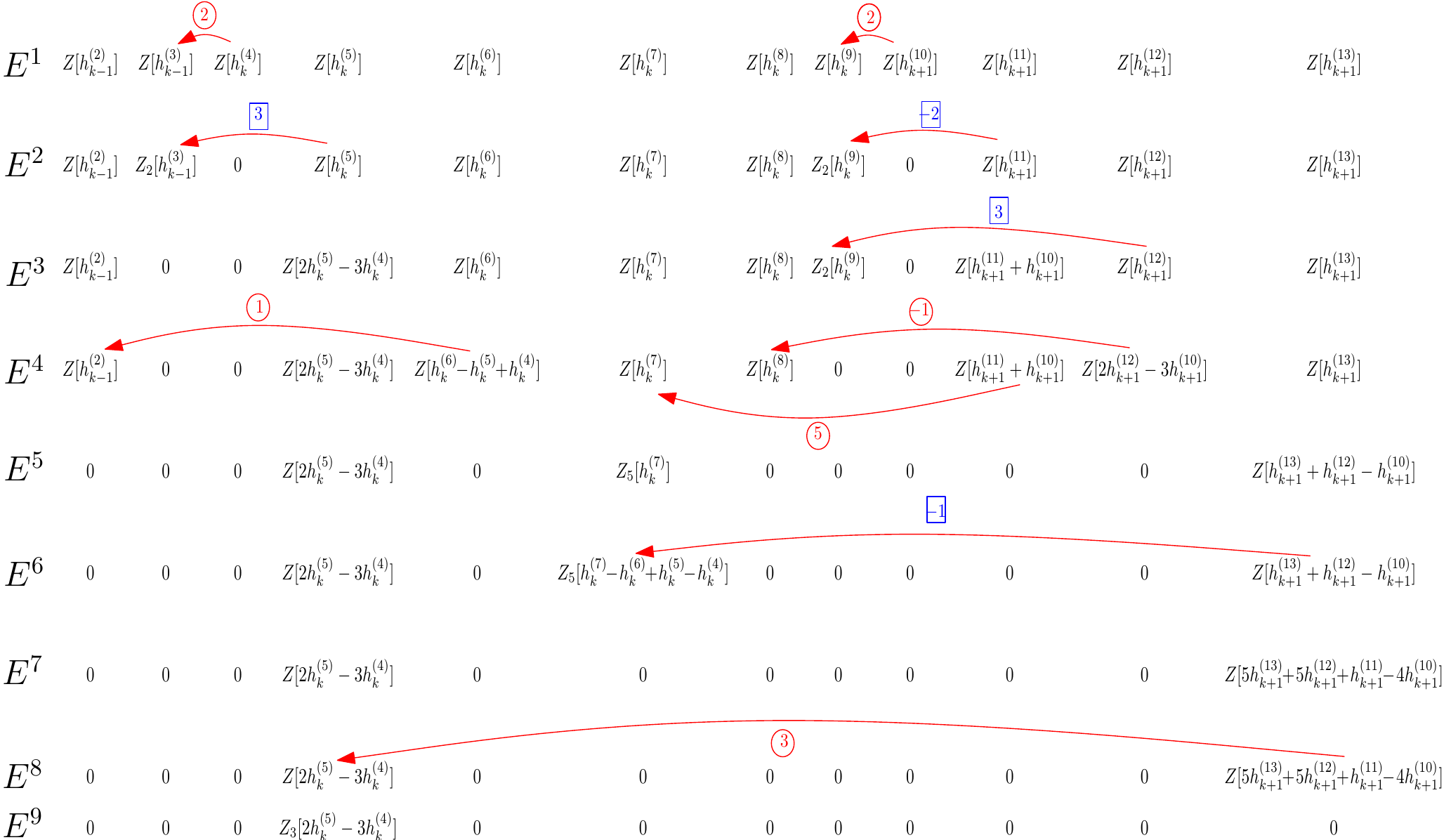}

\vspace{1cm}

{\sc {\footnotesize
\noindent M.A. Bertolim - {Salzburg, Austria.  e-mail: mbertolim@gmail.com} 

\noindent  D.V.S. Lima - IMECC, Universidade Estadual de Campinas, Campinas, SP, Brazil. CEP 13083-859.  \\ e-mail: dahisylima@gmail.com

\noindent  M.P. Mello - IMECC, Universidade Estadual de Campinas, Campinas, SP, Brazil. CEP 13083-859. \\ e-mail: margarid@ime.unicamp.br

\noindent  K.A. de Rezende - IMECC, Universidade Estadual de Campinas, Campinas, SP, Brazil. CEP 13083-859. \\ e-mail: ketty@ime.unicamp.br

 \noindent   M. R. da
Silveira  - CMCC, Universidade Federal do ABC, Santo Andr{\'e}, SP, Brazil. CEP 09210-580.\\ e-mail: mariana.silveira@ufabc.edu.br 
}}

\end{document}